\documentclass[a4paper,11pt]{article}
\usepackage[utf8]{inputenc}

\usepackage{amssymb}
\usepackage{mathtools}
\usepackage[margin=2cm]{geometry}
\usepackage{setspace}
\usepackage{authblk}
\usepackage{graphicx}
\usepackage{bm}   
\usepackage{cite}
\usepackage{multirow}
\usepackage[scriptsize,tight,figbotcap]{subfigure}
\usepackage{lineno}
\usepackage{xcolor}
\usepackage[title]{appendix}
\usepackage{hyperref}
\hypersetup{
    colorlinks,
    linkcolor={blue},
    citecolor={blue},
    urlcolor={blue}
}
\usepackage[font=small,labelfont=bf]{caption}
\usepackage{color}
\usepackage{booktabs}  
\usepackage{framed}
\usepackage{doi}
\usepackage{framed}

\title{A critical analysis of some popular methods for the discretisation of the gradient 
operator in finite volume methods}
\author[1]{Alexandros Syrakos\thanks{Corresponding author. E-mail address: 
\href{mailto:alexandros.syrakos@gmail.com}{alexandros.syrakos@gmail.com}, 
\href{mailto:syrakos@upatras.gr}{syrakos@upatras.gr}}}
\author[1]{Stylianos Varchanis}
\author[1]{Yannis Dimakopoulos}
\author[2]{Apostolos Goulas}
\author[1]{John Tsamopoulos}
\affil[1]{Laboratory of Fluid Mechanics and Rheology, Dept. of Chemical Engineering, University of 
Patras, 26500 Patras, Greece}
\affil[2]{Laboratory of Fluid Mechanics and Turbomachinery, Dept. of Mechanical Engineering, 
Aristotle University of Thessaloniki, 54124 Thessaloniki, Greece}
\date{}

\begin{document}
\newcommand{\vf}[1]{\bm{#1}}
\newcommand{\pd}[2]{\frac{\partial #1}{\partial #2}}
\newcommand{\smashedunderbrace}[2]{\smash{\underbrace{#1}_{#2}}\vphantom{#1}}
\newcommand{\GRAD}[1]{\nabla^{\mathrm{#1}}}

\maketitle

\begin{abstract}
Finite volume methods (FVMs) constitute a popular class of methods for the numerical simulation of 
fluid flows. Among the various components of these methods, the discretisation of the gradient 
operator has received less attention despite its fundamental importance with regards to the accuracy 
of the FVM. The most popular gradient schemes are the divergence theorem (DT) (or Green-Gauss) 
scheme, and the least-squares (LS) scheme. Both are widely believed to be second-order accurate, but 
the present study shows that in fact the common variant of the DT gradient is second-order accurate 
only on structured meshes whereas it is zeroth-order accurate on general unstructured meshes, and 
the LS gradient is second-order and first-order accurate, respectively. This is explained through a 
theoretical analysis and is confirmed by numerical tests. The schemes are then used within a FVM to 
solve a simple diffusion equation on unstructured grids generated by several methods; the results 
reveal that the zeroth-order accuracy of the DT gradient is inherited by the FVM as a whole, and the 
discretisation error does not decrease with grid refinement. On the other hand, use of the LS 
gradient leads to second-order accurate results, as does the use of alternative, consistent, DT 
gradient schemes, including a new iterative scheme that makes the common DT gradient consistent at 
almost no extra cost. The numerical tests are performed using both an in-house code and the popular 
public domain PDE solver OpenFOAM.
\end{abstract}


\begin{framed} \centering
\noindent This is the accepted version of the article published in:

Physics of Fluids 29, 127103 (2017); \doi{10.1063/1.4997682}
\end{framed}


\section{Introduction}
\label{sec: introduction}

Finite volume methods (FVMs) are used widely for the simulation of fluid flows; they are employed 
by several popular general-purpose CFD (Computational Fluid Dynamics) solvers, both commercial 
(e.g.\ ANSYS Fluent, STAR-CD) and open-source (e.g.\ OpenFOAM). One of the key components of FVMs 
is the approximation of the gradient operator. Computing the gradient of the dependent variables is 
needed for the FVM discretisation on non-Cartesian grids, where the fluxes across a face separating 
two finite volumes cannot be expressed as a function of the values of the variables at the centres 
of these two volumes alone, but additional terms involving the gradients must also be included 
\cite{Ferziger_2002}. The gradient operator is even more significant when solving complex flow 
problems such as turbulent flows modelled by the RANS and some LES methodologies 
\cite{Ferziger_2002} or non-Newtonian flows \cite{Oliveira_1998, Afonso_2012, Syrakos_2013, 
Jalali_2016}, where the additional equations solved (turbulence closure, constitutive equation 
etc.)\ may directly involve the velocity gradients. Apart from the main task of solving partial 
differential equations (PDEs), gradient calculation may also be important in auxiliary activities 
such as post-processing \cite{Correa_2011}.

The two most popular methods for calculating the gradient on grids of general geometry are based on 
the use of (a) the divergence theorem (DT) and (b) least-squares (LS) minimisation. They have both 
remained popular for nearly three decades of application of the FVM; for example, the DT method 
(also known as ``Green--Gauss gradient'') has been used in \cite{Barth_1989, Jasak_1996, 
Lilek_1997, Ferziger_2002, Wu_2014, Moukalled_2016}, and the LS method in \cite{Barth_1991, 
Muzaferija_1997, Ollivier_2002, Bramkamp_2004, Wu_2014, Moukalled_2016}. Often, general-purpose 
commercial or open-source CFD codes present the user with the option of choosing between these two 
methods. Their popularity stems from the fact that they are not algorithmically restricted to a 
particular grid cell geometry but can be applied to cells with an arbitrary number of faces. This 
property is important, because in recent years the use of unstructured grids is becoming standard 
practice in simulations of complex engineering processes. The tessellation of the complex geometries 
typically associated with such processes by structured grids is an arduous and extremely 
time-consuming procedure for the modeller, whereas unstructured grid generation is much more 
automated. Therefore, discretisation schemes are sought that are easily applicable on unstructured 
grids of arbitrary geometry whereas early FVMs either used Cartesian grids or relied on coordinate 
transformations which are applicable only on smooth structured grids.

In the literature, usually the gradient discretisation is only briefly discussed within an overall 
presentation of a FVM, with only a relatively limited number of studies devoted specifically to it
(e.g.\ \cite{Mavriplis_2003, Diskin_2008, Shima_2010, Betchen_2010, Correa_2011, Sozer_2014}). This 
suggests that existing gradient schemes are deemed satisfactory, and in fact there seems to be a 
widespread misconception that the DT and LS schemes are second-order accurate on any type of grid. 
For the DT gradient, a couple of studies have shown that it is potentially zeroth-order accurate, 
depending on the grid properties. Syrakos \cite{Syrakos_2006} noticed that it converged to wrong 
values on composite Cartesian grids in the vicinity of the interfaces between patches of different 
fineness. His analysis, given here in Sec.\ \ref{ssec: results refined}, concludes that this is due 
to grid skewness. Later, Sozer et al.\ \cite{Sozer_2014} tested a simpler variant of the DT gradient 
that uses arithmetic averaging instead of linear interpolation and proved that in the 
one-dimensional case it converges to incorrect values if the grid is not uniform. In numerical 
tests they also noticed that the scheme is inconsistent on two-dimensional grids of arbitrary 
topology. In fact, as early as in \cite{Barth_1989} it was briefly mentioned that the DT gradient 
fails to calculate exactly the gradients of linear functions, providing a hint to its inconsistency. 
An important question is whether this inconsistency inhibits convergence of the FVM to the correct 
solution. The results reported in \cite{Wu_2014} concerning the solution of a Poisson equation 
suggest zeroth-order FVM accuracy when the DT gradient is employed on skewed meshes, although the 
authors did not explicitly attribute this to inconsistency of the DT scheme.

The present paper presents mathematical analyses of the orders of accuracy of both methods on 
various types of grid, backed by numerical experiments. The common DT method is proved to be 
zeroth-order accurate on grids of arbitrary geometry and second-order accurate on smooth structured 
meshes. This is shown to be true even if a finite number of ``corrector'' steps is performed (e.g.\ 
\cite{Karimian_2006, Wu_2014}) whereby an iterative procedure uses the gradient calculated at the 
previous iteration to improve the accuracy. Theoretically, this procedure increases the order of 
accuracy only in the limit of infinite iterations. Practically, the desired accuracy can be achieved 
with a finite number of iterations, but this number is a priori unknown and increases with the grid 
fineness. Furthermore, corrector steps make the DT gradient much more expensive than the LS. 
However, since the FVM solution typically proceeds with a number of outer iterations, we exploit 
this fact to interweave the gradient iterations with the FVM outer iterations to obtain a 
first-order accurate DT gradient at a cost that is nearly the same as that of the uncorrected DT 
scheme. On the other hand, the LS method is shown to be first-order accurate on arbitrary grids and 
second-order accurate on structured grids, while for the default distance-based weighting scheme a 
particular non-integer exponent ($-3/2$) enlarges the set of grid types for which the method is 
second-order accurate. The different order of accuracy of the gradient schemes on structured versus 
unstructured grids is attributed to the fact that the former become less skewed and more uniform as 
they are refined, whereas the latter usually do not.

First-order accurate gradients are, for most purposes, compatible with second-order accurate FVMs 
because differentiation of a second-order accurate variable gives a first-order accurate derivative. 
In Section \ref{sec: PDE solution} we proceed beyond the analysis of the gradient schemes 
themselves and test them as components of a FVM for the solution of a Poisson equation on various 
kinds of unstructured grids. Experiments with both an in-house code and OpenFOAM show that the FVMs 
are zeroth-order accurate when they employ the common DT gradient, whereas they are second-order 
accurate if they employ instead the least-squares gradient, or the proposed iteratively corrected 
DT gradient, or an alternative DT gradient scheme, achieving accuracy that is in fact only 
marginally worse than that obtained on Cartesian grids.

\section{Preliminary considerations}
\label{sec: preliminary}

\begin{figure}[t]
 \centering
 \includegraphics[scale=0.85]{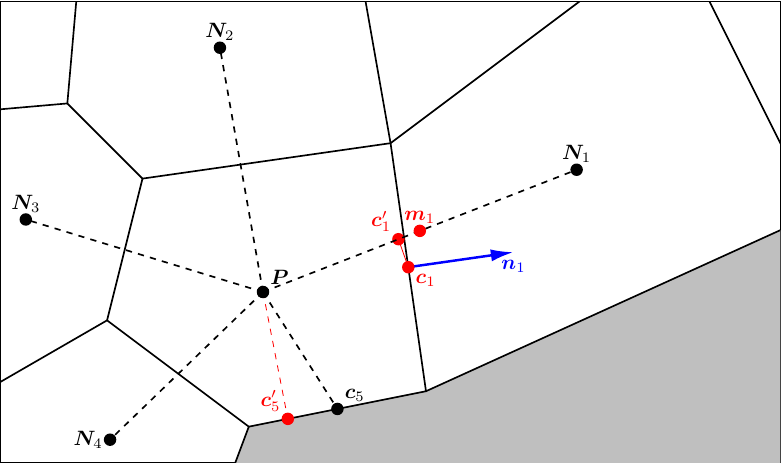}
 \caption{Part of an unstructured grid, showing cell $P$ and its neighbouring cells, each having a 
single common face with $P$. The shaded area lies outside the grid. The faces and neighbours of 
cell 
$P$ are numbered in anticlockwise order, with face $f$ separating $P$ from its neighbour $N_f$. The 
geometric characteristics of its face 1 which separates it from neighbouring cell $N_1$ are 
displayed. The position vectors of the centroids of cells $P$ and $N_f$ are denoted by the same 
characters but in boldface as $\vf{P}$ and $\vf{N}_f$; $\vf{c}_f$ is the centroid of face $f$ and 
$\vf{c}'_f$ is its closest point on the line segment connecting $\vf{P}$ and $\vf{N}_f$; 
$\vf{m}_f$ is the midpoint between $\vf{P}$ and $\vf{N}_f$; $\vf{n}_f$ is the unit vector normal to 
face $f$, pointing outwards of $P$. The shown cell $P$ also has a boundary face (face 5), with no 
neighbour on the other side; in this case $\vf{c}'_5$ is the projection of $\vf{P}$ onto the 
boundary face.}
 \label{fig: unstructured grid}
\end{figure}

We will focus on two-dimensional problems, although the one- and three-dimensional cases are also 
discussed where appropriate and the conclusions are roughly the same in all dimensions. Let $x$, 
$y$ denote the usual Cartesian coordinates, and let $\phi(x,y)$ be a function defined over a 
domain, whose gradient $\nabla\phi$ we wish to calculate. It is convenient to introduce a convention 
where subscripts beginning with a dot (.) denote differentiation with respect to the ensuing 
variable(s), e.g.\ $\phi_{\!.x} \equiv \partial \phi / \partial x$, $\phi_{\!.xy} \equiv \partial^2 
\phi / \partial x \partial y$ etc. If the variables $x$, $y$ etc.\ are used as subscripts without a 
leading dot then they are used simply as indices, without any differentiation implied. Therefore, we 
seek the gradient $\nabla \phi = (\phi_{\!.x}, \phi_{\!.y})$ of the function $\phi$. The domain 
over which the function is defined is discretised by a grid into a number of non-overlapping finite 
volumes, or \textit{cells}. A cell can be arbitrarily shaped, but its boundary must consist of a 
number of straight faces, as in Fig.\ \ref{fig: unstructured grid}, each separating it from a 
single other cell or from the exterior of the domain (the latter are called boundary faces). We 
assume a cell-centred finite volume method, meaning that the values of $\phi$ are known only at the 
geometric centres of the cells and at the centres of the boundary faces. The notation that is 
adopted in order to describe the geometry of the grid is presented in Fig.\ \ref{fig: unstructured 
grid}. Also, we will denote vectors by boldface characters.

Our goal is to derive approximate algebraic gradient operators $\nabla^\mathrm{a}$ which return 
values $\nabla^{\mathrm{a}} \phi(\vf{P}) \approx \nabla \phi(\vf{P})$ at each cell centroid 
$\vf{P}$, using information only from the immediate neighbouring cells and boundary faces. The 
components of the approximate gradient are denoted as $\nabla^{\mathrm{a}}\phi = 
(\phi_{\!.x}^{\mathrm{a}}, \phi_{\!.y}^{\mathrm{a}})$. The operators $\nabla^{\mathrm{a}}$ must be 
capable of approximating the derivative on grids of arbitrary geometry. To study the effect of grid 
geometry on the accuracy, we must define some indicators of grid irregularity. With the present grid 
arrangement, we find it useful to define three kinds of such grid irregularities, which we shall 
call here ``non-orthogonality'', ``unevenness'' and ``skewness'' (other possibilities exist, see 
e.g.\ \cite{Kallinderis_2009}). We will define these terms with the aid of Fig.\ \ref{fig: 
unstructured grid}. With the nomenclature defined in that figure, we will say that face $f$ of cell 
$P$ exhibits non-orthogonality if $\vf{N}_f - \vf{P}$ is not parallel to $\vf{n}_f$; a measure of 
non-orthogonality is the angle between the vectors $(\vf{N}_f - \vf{P})$ and $\vf{n}_f$. Also, we 
will say that face $f$ exhibits unevenness if the midpoint of the line segment joining $\vf{P}$ and 
$\vf{N}_f$, $\vf{m}_f = (\vf{P}+\vf{N}_f)/2$, does not coincide with $\vf{c}'_f$ (i.e.\ the cell 
centres are unequally spaced on either side of the face); a measure of unevenness is $\|\vf{c}'_f - 
\vf{m}_f\| / \|\vf{N}_f - \vf{P}\|$. Finally, we will say that face $f$ exhibits skewness if 
$\vf{c}'_f$ does not coincide with $\vf{c}_f$ (i.e. the line joining the cell centres does not pass 
through the centre of the face); a measure of the skewness is $\|\vf{c}_f - \vf{c}'_f\| / \|\vf{N}_f 
- \vf{P}\|$. 

\begin{figure}[t]
 \centering
 \includegraphics[scale=0.80]{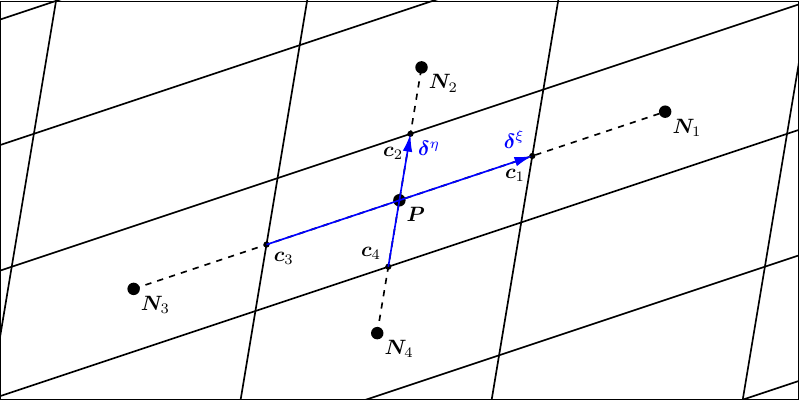}
 \caption{Part of a grid formed by equispaced parallel grid lines. See the text for details.}
 \label{fig: parallelogram grid}
\end{figure}

Before discussing the DT and LS gradient schemes, it is useful to examine a simpler scheme which is 
applicable on very plain grids that are formed from two families of equispaced parallel straight 
lines intersecting at a constant angle, as in Fig.\ \ref{fig: parallelogram grid}. In this case all 
the cells of the grid are identical parallelograms (Cartesian grids with constant spacing are a 
special case of this category). Figure \ref{fig: parallelogram grid} shows a cell $P$ belonging to 
such a grid, and its four neighbouring cells. The vectors $\vf{\delta}^{\xi}$ and 
$\vf{\delta}^{\eta}$ are parallel to the grid lines and span the size of the cells. Due to the grid 
properties it holds that $\vf{\delta}^{\xi} = \vf{c}_1 - \vf{c}_3 = \vf{N}_1 - \vf{P} = \vf{P} - 
\vf{N}_3$ and $\vf{\delta}^{\eta} = \vf{c}_2 - \vf{c}_4 = \vf{N}_2 - \vf{P} = \vf{P} - \vf{N}_4$. It 
can be assumed that two variables, $\xi$ and $\eta$, are distributed in the domain, such that in 
the direction of $\vf{\delta}^{\xi}$ the variable $\xi$ varies linearly while $\eta$ is constant, 
and in the direction of $\vf{\delta}^{\eta}$ the variable $\eta$ varies linearly while $\xi$ is 
constant. Then the grid can be considered to be constructed by drawing lines of constant $\xi$ and 
of constant $\eta$, equispaced (in the $\xi, \eta$ space) by $\Delta \xi$ and $\Delta \eta$, 
respectively. Let us assume also that the grid density can be increased by adjusting the spacings 
$\Delta \xi$ and $\Delta \eta$, but their ratio must be kept constant, e.g.\ if $\Delta \xi = h$ 
then $\Delta \eta = \alpha h$ with $\alpha$ being a constant, independent of $h$. The variable $h$ 
determines the grid fineness. Therefore, the direction of the grid vectors $\vf{\delta}^{\xi}$ and 
$\vf{\delta}^{\eta}$ remains constant with grid refinement, but their lengths are proportional to 
the grid parameter $h$.

This idealized grid exhibits no unevenness or skewness. It possibly exhibits non-orthogonality, but 
this poses no problem as far as the gradient calculation is concerned. Since points $\vf{N}_3$, 
$\vf{P}$ and $\vf{N}_1$ are collinear and equidistant, the rate of change of any quantity $\phi$ in 
the direction of $\vf{\delta}^{\xi}$ can be approximated at point $\vf{P}$ from the values at 
$\vf{N}_3$ and $\vf{N}_1$ using second-order accurate central differencing. In the same manner, the 
rate of change in the direction of $\vf{\delta}^{\eta}$ at $\vf{P}$ can be approximated from the 
values at $\vf{N}_4$ and $\vf{N}_2$. So, let the grid vectors be written in Cartesian coordinates as 
$\vf{\delta}^{\xi} = (\delta^{\xi}_x, \delta^{\xi}_y)$ and $\vf{\delta}^{\eta} = (\delta^{\eta}_x, 
\delta^{\eta}_y)$, respectively. Then, by expanding the function $\phi$ in a two-dimensional Taylor 
series along the Cartesian directions, centred at point $\vf{P}$, and using that to express the 
values at the points $\vf{N}_1 = \vf{P} + \vf{\delta}^{\xi}$, $\vf{N}_3 = \vf{P} - 
\vf{\delta}^{\xi}$, $\vf{N}_2 = \vf{P} + \vf{\delta}^{\eta}$ and $\vf{N}_4 = \vf{P} - 
\vf{\delta}^{\eta}$ we get
\begin{align*}
  \phi(\vf{N}_1) - \phi(\vf{N}_3) &= \nabla \phi(\vf{P}) \cdot (2\vf{\delta}^{\xi})  \;+\; O(h^3) 
\\
  \phi(\vf{N}_2) - \phi(\vf{N}_4) &= \nabla \phi(\vf{P}) \cdot (2\vf{\delta}^{\eta}) \;+\; O(h^3)
\end{align*}
This system can be solved for \(\nabla\phi(\vf{P}) \equiv (\phi_{\!.x},\phi_{\!.y})\) to get
\begin{equation*}
 \begin{bmatrix}
   \phi_{\!.x}(\vf{P}) \\[0.15cm]
   \phi_{\!.y}(\vf{P})
 \end{bmatrix}
 \; = \quad
 \frac{1}{2\Omega_P}
 \left[ \begin{array}{rr}
   \delta_y^{\eta}  & -\delta_y^{\xi}  \\[0.15cm]
  -\delta_x^{\eta}  &  \delta_x^{\xi}
 \end{array} \right]
 \!\cdot\!
 \begin{bmatrix}
  \phi(\vf{N}_1) - \phi(\vf{N}_3) \\[0.15cm]
  \phi(\vf{N}_2) - \phi(\vf{N}_4)
 \end{bmatrix}
 \quad + \quad
 \frac{1}{2\Omega_P}
 \left[ \begin{array}{rr}
   \delta_y^{\eta}  & -\delta_y^{\xi}  \\[0.15cm]
  -\delta_x^{\eta}  &  \delta_x^{\xi}
 \end{array} \right]
 \!\cdot\!
 \begin{bmatrix}
  O(h^3) \\[0.15cm]
  O(h^3)
 \end{bmatrix}
\end{equation*}
where $\Omega_P = \delta_x^{\xi} \delta_y^{\eta} - \delta_y^{\xi} \delta_x^{\eta} = \| 
\vf{\delta}^{\xi} \times \vf{\delta}^{\eta}\|$ is the volume (area in 2D) of cell $P$. The last term 
in the above equation, involving the unknown $O(h^3)$ terms, is of order $O(h^2)$ because $\Omega_P 
= O(h^2)$ and all the $\delta$'s are $O(h)$. So, carrying out the matrix multiplications we arrive 
at
\begin{equation} \label{eq: grad parallelogram}
 \nabla\phi(\vf{P}) \;\equiv\;
 \begin{bmatrix}
   \phi_{\!.x}(\vf{P}) \\[0.15cm]
   \phi_{\!.y}(\vf{P})
 \end{bmatrix}
 \;=\;
 \underbrace{
 \frac{1}{2\Omega_P}
 \begin{bmatrix}
   \delta_y^{\eta} (\phi(\vf{N}_1) \!-\! \phi(\vf{N}_3))  \,-\,  \delta_y^{\xi}  (\phi(\vf{N}_2) 
\!-\! \phi(\vf{N}_4)) \\[0.15cm]
   \delta_x^{\xi}  (\phi(\vf{N}_2) \!-\! \phi(\vf{N}_4))  \,-\,  \delta_x^{\eta} (\phi(\vf{N}_1) 
\!-\! \phi(\vf{N}_3))
 \end{bmatrix}
 }_{\nabla^{\mathrm{s}} \phi (\vf{P})}
 \;+\;
 \begin{bmatrix}
  O(h^2) \\[0.15cm]
  O(h^2)
 \end{bmatrix}
\end{equation}
Finally, we drop the unknown $O(h^2)$ terms in Eq.\ \eqref{eq: grad parallelogram} and are left with 
a second-order approximation to the gradient, $\nabla^{\mathrm{s}}$; the dropped terms are called 
the \textit{truncation error} of the operator $\nabla^{\mathrm{s}}$. The fact that this formula has 
been derived for a grid constructed from equidistant parallel lines may seem too restrictive, but in 
fact the utility of Eq.\ \eqref{eq: grad parallelogram} goes beyond this narrow context, as will now 
be explained.

\begin{figure}[t]
 \centering
 \includegraphics[scale=0.90]{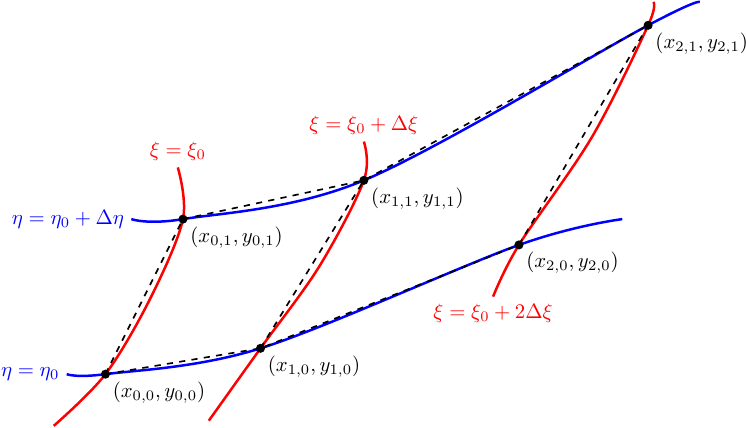}
 \caption{Part of a grid (dashed straight lines) constructed from lines of constant $\xi$ (red 
curves) and $\eta$ (blue curves), where $\xi$, $\eta$ are variables distributed smoothly in the 
domain. The lines are equispaced with constant spacings $\Delta \xi$ and $\Delta \eta$, 
respectively. Grid node $(i,j)$ is located at $(\xi_i,\eta_j) = (\xi_0 + i\,\Delta \xi, \eta_0 + 
j\,\Delta \eta)$ in the computational space, where $(\xi_0,\eta_0)$ is a predefined point, and at 
$(x_{i,j},y_{i,j})$ in the physical space.}
 \label{fig: structured grid}
\end{figure}

Consider again structured grid generation based on a pair of variables $\xi$, $\eta$ distributed 
smoothly in the domain, where curves of constant $\xi$ and $\eta$ are drawn at equal intervals of 
$\Delta \xi$ and $\Delta \eta$, respectively. This time $\xi$ and $\eta$ are not required to vary 
linearly nor to be constant along straight lines. Therefore, a curvilinear grid such as that shown 
in Fig.\ \ref{fig: structured grid} may result, constructed by joining the points of intersection 
of these two families of curves by straight line segments (the dashed lines in Fig.\ \ref{fig: 
structured grid}).

There is a one-to-one correspondence of coordinates $(x,y)$ of the physical space to coordinates 
$(\xi,\eta)$ of the computational space. Therefore, not only are the computational coordinates 
functions of the physical coordinates ($\xi = \xi(x,y)$ and $\eta = \eta(x,y)$), but the physical 
coordinates $(x,y)$ can also be regarded as functions of the computational coordinates ($x = 
x(\xi,\eta)$ and $y = y(\xi,\eta)$). Since the latter vary smoothly in the domain, $(x,y)$ can be 
expanded in Taylor series around a reference point $(\xi_0,\eta_0)$. For example, for the $x$ 
coordinate,
\begin{equation} \label{eq: taylor x}
 x(\xi_0+\delta\xi, \eta_0+\delta\eta) \;=\; 
 x(\xi_0,\eta_0) \;+\;
 x_{\!.\xi}\delta\xi \;+\; x_{\!.\eta}\delta\eta
 \;+\;
 O(h^2)
\end{equation}
where the derivatives are evaluated at point $(\xi_0,\eta_0)$. In this way, we can express the 
coordinates of all the grid vertices $(x_{i,j}, y_{i,j})$ shown in Fig.\ \ref{fig: structured grid} 
as functions of $(x_{0,0},y_{0,0}) \equiv (x(\xi_0,\eta_0), y(\xi_0,\eta_0))$ and of the 
derivatives 
$x_{.\xi}$ etc.\ there. Using these expansions, it is easy to show that 
\begin{equation} \label{eq: structured skewness}
 \lim_{\Delta\xi,\Delta\eta \rightarrow 0} \left( \frac{x_{1,1}-x_{1,0}}{x_{0,1}-x_{0,0}} \right) 
\;=\;
 \lim_{\Delta\xi,\Delta\eta \rightarrow 0} \left( \frac{y_{1,1}-y_{1,0}}{y_{0,1}-y_{0,0}} \right) 
\;=\; 1
\end{equation}
and
\begin{equation} \label{eq: structured unevenness}
 \lim_{\Delta\xi,\Delta\eta \rightarrow 0} \left( \frac{x_{2,0}-x_{1,0}}{x_{1,0}-x_{0,0}} \right) 
\;=\;
 \lim_{\Delta\xi,\Delta\eta \rightarrow 0} \left( \frac{y_{2,0}-y_{1,0}}{y_{1,0}-y_{0,0}} \right) 
\;=\; 1
\end{equation}

Equation \eqref{eq: structured skewness} implies that, as the grid is refined by reducing the 
$\Delta \xi$, $\Delta \eta$ spacings, neighbouring grid lines of the same family become more and 
more parallel, so that grid skewness tends to zero. Equation \eqref{eq: structured unevenness} 
implies that, as the grid is refined, neighbouring cells tend to become of equal size, so that grid 
unevenness tends to zero. The conclusion is that on structured grids which are constructed from 
smooth distributions of auxiliary variables $(\xi,\eta)$, grid refinement\footnote{It is stressed 
that grid refinement must be performed in the \textit{computational} space $(\xi,\eta)$, by 
simultaneously reducing the spacings $\Delta \xi$ and $\Delta \eta$. Otherwise, if refinement is 
performed directly in the physical space $(x,y)$, for example by joining the centroid of each cell 
to the centroids of its faces thus splitting it into four child cells, then equations \eqref{eq: 
structured skewness} -- \eqref{eq: structured unevenness} do not necessarily hold.} causes the 
geometry of a cell and its neighbours to locally approach that depicted in Fig.\ \ref{fig: 
parallelogram grid}. This has an important consequence: any numerical scheme for computing the 
gradient that reduces to Eq.\ \eqref{eq: grad parallelogram} on parallelogram grids (such as that of 
Fig.\ \ref{fig: parallelogram grid}) is of second-order accuracy on smooth structured grids (such as 
that of Fig.\ \ref{fig: structured grid})\footnote{To be precise, this depends on the skewness 
decreasing fast enough with grid refinement -- see Section \ref{sec: gauss}. For structured grids, 
using Taylor series such as Eq.\ \eqref{eq: taylor x}, it can be shown that the skewness is 
$O(h)$.}. It so happens that both the DT scheme and the LS scheme belong to this category. 
Unfortunately, grid refinement does not engender such a quality improvement when it comes to 
unstructured meshes.

\section{Gradient calculation using the divergence theorem}
\label{sec: gauss}

For grids of more general geometry, like that of Fig.\ \ref{fig: unstructured grid}, a more general 
method is needed. Let $\Omega_P$ denote the volume of cell $P$ and $S_P$ its bounding surface. The 
DT gradient scheme is based on a derivative of the divergence theorem, which can be expressed as 
follows:
\begin{equation*}
 \int_{\Omega_P} \nabla \phi \, \mathrm{d}\Omega \;=\; \int_{S_P} \phi \, \vf{n} \, \mathrm{d}s
\end{equation*}
where $\vf{n}$ is the unit vector perpendicular to $S_P$ at each point, pointing outwards of the 
cell, while $\mathrm{d}\Omega$ and $\mathrm{d}s$ are infinitesimal elements of the volume and 
surface, respectively.   The bounding surface $S_P$ can be decomposed into $F$ faces which are 
denoted by $S_f$, $f = 1, \ldots, F$ ($F = 5$ in Fig.\ \ref{fig: unstructured grid}). These faces 
are assumed to be straight (planar, in 3 dimensions), as in Fig.\ \ref{fig: unstructured grid}, so 
that the normal unit vector $\vf{n}$ has a constant value $\vf{n}_f$ along each face $f$. Therefore, 
the above equation can be written as
\begin{equation} \label{eq: gauss theorem}
 \int_{\Omega_P} \nabla \phi \, \mathrm{d}\Omega \;=\; \sum_{f=1}^F \left( \vf{n}_f \int_{S_f} \phi 
\, \mathrm{d}s \right)
\end{equation}

According to the midpoint integration rule \cite{Ferziger_2002, Burden_2011}, the mean value of a 
quantity over cell $P$ (or face $f$) is equal to its value at the centroid $\vf{P}$ of the cell (or 
$\vf{c}_f$ of the face), plus a second-order correction term. Applying this to the mean values of 
$\nabla \phi$ and $\phi$ over $\Omega_P$ and $S_f$ we get, respectively:
\begin{equation} \label{eq: midpoint rule cell}
 \frac{1}{\Omega_P} \int_{\Omega_P} \nabla \phi \, \mathrm{d}\Omega \;=\; \nabla \phi (\vf{P}) 
\;+\; \vf{O}(h^2)
 \; \Rightarrow \quad 
 \int_{\Omega_P} \nabla \phi \, \mathrm{d}\Omega \;=\; \nabla \phi (\vf{P}) \, \Omega_P \;+\; 
\vf{O}(h^4)
\end{equation}
\begin{equation} \label{eq: midpoint rule faces}
 \frac{1}{S_f} \int_{S_f} \phi \, \mathrm{d}s \;=\; \phi(\vf{c}_f) \;+\; O(h^2) 
 \; \Rightarrow \quad 
 \int_{S_f} \phi \, \mathrm{d}s \;=\; \phi(\vf{c}_f) \, S_f \;+\; O(h^3)
\end{equation}
where $h$ is a characteristic grid spacing (we used that $\Omega_P = O(h^2)$ and $S_f = O(h)$). 
Substituting these expressions into the divergence theorem, Eq.\ \eqref{eq: gauss theorem}, we get
\begin{equation} \label{eq: gauss gradient halfway}
 \nabla \phi(\vf{P}) \;=\; \frac{1}{\Omega_P} \sum_{f=1}^F \phi(\vf{c}_f) \, S_f\, \vf{n}_f \;+\; 
\vf{O}(h)
\end{equation}
The above formula is exact as long as the unknown $\vf{O}(h)$ term is retained, which consists 
mostly of face contributions arising from Eq.\ \eqref{eq: midpoint rule faces}, whereas the volume 
contribution from Eq.\ \eqref{eq: midpoint rule cell} is only $\vf{O}(h^2)$. Dropping this term 
would leave us with a first-order accurate approximation. However, this approximation cannot be 
the final formula for the gradient because Eq.\ \eqref{eq: gauss gradient halfway} contains 
$\phi(\vf{c}_f)$, the $\phi$ values at the face centres, whereas we need a formula that uses only 
the values at the \textit{cell} centres. The common practice is to approximate $\phi(\vf{c}_f)$ by 
$\phi(\vf{c}'_f)$, the exact values of $\phi$ at points $\vf{c}'_f$ rather than $\vf{c}_f$ (see 
Fig.\ \ref{fig: unstructured grid}); these values also do not belong to the set of cell-centre 
values, but since $\vf{c}'_f$ lies on the line segment joining cell centres $\vf{P}$ and 
$\vf{N}_f$, the value $\phi(\vf{c}'_f)$ can in turn be approximated by linear interpolation between 
$\phi(\vf{P})$ and $\phi(\vf{N}_f)$, say $\bar{\phi}(\vf{c}'_f)$ (the overbar denotes linear 
interpolation):
\begin{equation} \label{eq: linear interpolation}
 \bar{\phi}(\vf{c}'_f) \;\equiv\; \frac{\| \vf{c}'_f - \vf{N}_f \|}{\| \vf{N}_f - \vf{P} \|} \, 
\phi(\vf{P})
                       \;+\;      \frac{\| \vf{c}'_f - \vf{P} \|}{\| \vf{N}_f - \vf{P} \|} \, 
\phi(\vf{N}_f)
                       \;=\; \phi(\vf{c}'_f) \;+\; O(h^2)
\end{equation}
Linear interpolation is known to be second-order accurate, hence the $O(h^2)$ term in the above 
equation. Thus, by using $\bar{\phi}(\vf{c}'_f)$ instead of $\phi(\vf{c}_f)$ in the right hand side 
of Eq.\ \eqref{eq: gauss gradient halfway} and dropping the unknown $\vf{O}(h)$ term we obtain an 
approximation to the gradient which depends only on cell-centre values of $\phi$:
\begin{equation} \label{eq: gauss gradient 0}
 \nabla^{\mathrm{d}0} \phi (\vf{P}) \;\equiv\; \frac{1}{\Omega_P} \sum_{f=1}^F 
\bar{\phi}(\vf{c}'_f) \, S_f \, \vf{n}_f 
\end{equation}
This is called the ``divergence theorem'' (DT) gradient. It applies in both two and three 
dimensions.

An important question is whether and how the replacement of $\phi(\vf{c}_f)$ by 
$\bar{\phi}(\vf{c}'_f)$, that led from Eq.\ \eqref{eq: gauss gradient halfway} to the formula 
\eqref{eq: gauss gradient 0}, has affected the accuracy. To answer this question we first need to 
deduce how much $\bar{\phi}(\vf{c}'_f)$ differs from $\phi(\vf{c}_f)$. This can be done by using 
the exact $\phi(\vf{c}'_f)$ as an intermediate value. We will consider structured grids and 
unstructured grids separately.

\subsubsection*{Structured grids}

As discussed in the previous Section, the skewness of smooth structured grids diminishes with 
refinement, and in fact expressing the points involved in its definition as Taylor series of the 
form \eqref{eq: taylor x} it can be shown that $\|\vf{c}_f - \vf{c}'_f\| / \|\vf{N}_f - \vf{P}\| = 
O(h)$, or $\vf{c}_f - \vf{c}'_f = \vf{O}(h^2)$. Therefore, expanding $\phi(\vf{c}_f)$ in a Taylor 
series about point $\vf{c}'_f$ gives $\phi(\vf{c}_f) = \phi(\vf{c}'_f) + \nabla \phi(\vf{c}'_f) 
\cdot (\vf{c}_f - \vf{c}'_f) + O(h^2) = \phi(\vf{c}'_f) + O(h^2)$. Next, $\phi(\vf{c}'_f)$ can be 
expressed in terms of $\bar{\phi}(\vf{c}'_f)$ according to Eq.\ \eqref{eq: linear interpolation}. 
Putting everything together results in the sought relationship:
\begin{equation} \label{eq: phi_c vs phi_c_bar structured}
 \phi(\vf{c}_f) \;=\; \phi(\vf{c}'_f) + O(h^2) \;=\; (\bar{\phi}(\vf{c}'_f) + O(h^2)) + O(h^2) 
\;=\; 
\bar{\phi}(\vf{c}'_f) + O(h^2) 
\end{equation}
Substituting this into Eq.\ \eqref{eq: gauss gradient halfway}, and considering that $\Omega_P = 
O(h^2)$ and $S_f = O(h)$, we arrive at
\begin{equation} \label{eq: gauss gradient exact structured}
 \nabla \phi (\vf{P}) \;=\; \frac{1}{\Omega_P} \sum_{f=1}^F \bar{\phi}(\vf{c}'_f) \, S_f \, 
\vf{n}_f \;+\; \vf{O}(h)
\end{equation}
The first term on its right-hand side is just the approximate gradient $\nabla^{\mathrm{d0}}$, Eq.\ 
\eqref{eq: gauss gradient 0}. Therefore, the truncation error of $\nabla^{\mathrm{d0}}$ is 
$\vf{\tau} \equiv \nabla^{\mathrm{d}0}\phi(\vf{P}) - \nabla\phi(\vf{P})  = \vf{O}(h)$. Actually, the 
situation is even better because the leading terms of the contributions of opposite faces to the 
$\vf{O}(h)$ term in Eq.\ \eqref{eq: gauss gradient exact structured} cancel out leaving a net 
$\vf{O}(h^2)$ truncation error, making the method second-order accurate. The proof is tedious and 
involves writing analytic expressions for the truncation error contributions of each face; it can be 
performed more easily for a grid such as that shown in Fig.\ \ref{fig: parallelogram grid} where the 
divergence theorem procedure described in the present Section is an alternative path to arrive at 
the exact same formula, Eq.\ \eqref{eq: grad parallelogram}, $\nabla^{\mathrm{d}0} \equiv 
\nabla^{\mathrm{s}}$, with its $\vf{O}(h^2)$ truncation error.

\begin{figure}[tb]
 \centering
 \includegraphics[scale=0.90]{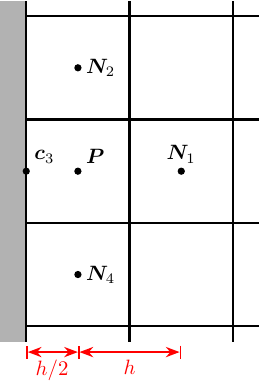}
 \caption{A boundary cell $P$ belonging to a Cartesian grid.}
 \label{fig: boundary grid}
\end{figure}

The error cancellation between opposite faces does not always occur; sometimes the worst-case 
scenario of $\vf{O}(h)$ truncation error predicted by Eq.\ \eqref{eq: gauss gradient exact 
structured} holds. An example is at boundary cells. Fig.\ \ref{fig: boundary grid} shows a boundary 
cell $P$ belonging to a Cartesian grid which exhibits neither skewness nor unevenness (in fact, it 
does not even exhibit non-orthogonality). Yet cell $P$ has no neighbour on the boundary side, and 
so the centre of its boundary face, $\vf{c}_3$, is used instead of a neighbouring cell centre. This 
introduces ``unevenness'' in the $x$- (horizontal) direction because the distances $\|\vf{N}_1 - 
\vf{P}\| = h$ and $\|\vf{c}_3 - \vf{P}\| = h/2$ are not equal. The $x$- component of the divergence 
theorem gradient \eqref{eq: gauss gradient 0} reduces to
\begin{equation} \label{eq: gradient at boundary}
 \phi_{\!.x}^{\mathrm{d}0}(\vf{P}) \;=\;
 \frac{1}{2h} \left( \phi(\vf{N}_1) + \phi(\vf{P}) - 2\phi(\vf{c}_3) \right)
\end{equation}
which is only first-order accurate, since expanding $\phi(\vf{N}_1)$ and $\phi(\vf{c}_3)$ in Taylor 
series about $\vf{P}$ gives
\begin{equation*}
 \phi_{\!.x}^{\mathrm{d}0}(\vf{P}) \;=\;
 \phi_{\!.x}(\vf{P}) \;+\; \frac{1}{8} \phi_{\!.xx}(\vf{P}) \, h \;+\; O(h^2)
\end{equation*}
This offers a nice demonstration of the effect of error cancellation between opposite faces. 
Equation \eqref{eq: gauss gradient 0} is obtained from Eq.\ \eqref{eq: gauss gradient halfway} by 
using interpolated values $\bar{\phi}(\vf{c}'_f)$ instead of the exact but unknown values 
$\phi(\vf{c}_f)$ at the face centres. Therefore, one might expect that since the cell of Fig.\ 
\ref{fig: boundary grid} has a boundary face and the exact value at its centre, $\phi(\vf{c}_3)$, 
is used rather than an interpolated value, the result would be more accurate; but the above 
analysis shows exactly the opposite: the error increases from $O(h^2)$ to $O(h)$. This is due to 
the fact that by dropping the interpolation error on the boundary face the corresponding error on 
the opposite face $1$ is no longer counterbalanced.

\begin{figure}[thb]
 \centering
 \includegraphics[scale=0.85]{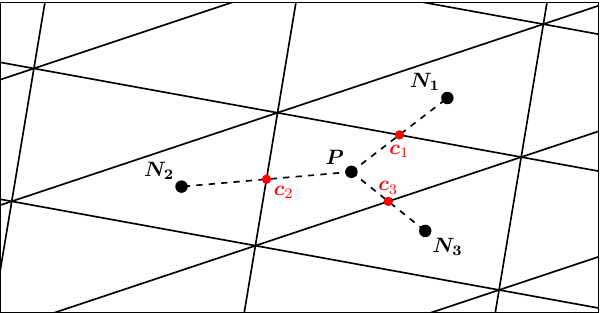}
 \caption{Grid of triangles constructed through bisection of the cells of the structured grid of 
Fig.\ \ref{fig: parallelogram grid} along their same diagonal.}
 \label{fig: triangular structured grid}
\end{figure}

Another example where error cancellation does not occur and the formal $\vf{O}(h)$ accuracy 
predicted by Eq.\ \eqref{eq: gauss gradient exact structured} holds is the common case of grid that 
consists of triangles that come from dividing each cell of a smooth structured grid along the same 
diagonal. For example, application of this procedure to the grid shown in Fig.\ \ref{fig: 
parallelogram grid} results in the grid of Fig.\ \ref{fig: triangular structured grid}. The latter 
may be seen to exhibit neither unevenness nor skewness as the face centres $\vf{c}_f$ coincide with 
the midpoints of the line segments joining cell centre $\vf{P}$ to its neighbours $\vf{N}_f$. Thus 
Eq.\ \eqref{eq: phi_c vs phi_c_bar structured} holds, leading to Eq.\ \eqref{eq: gauss gradient 
exact structured}. This time, however, a tedious but straightforward calculation where neighbouring 
$\phi(\vf{N}_f)$ values are expressed in Taylor series about $\vf{P}$ and substituted in \eqref{eq: 
gauss gradient exact structured} shows that there is no cancellation and the truncation error 
remains $\vf{O}(h)$. If the grid comes from triangulation of a curvilinear structured grid (Fig.\ 
\ref{fig: structured grid}) then the skewness is not zero but it diminishes with refinement at an 
$O(h)$ rate (Eq.\ \eqref{eq: phi_c vs phi_c_bar structured}) and exactly the same conclusions hold.

In fact, from the considerations leading to Eq.\ \eqref{eq: phi_c vs phi_c_bar structured} it 
follows that if a grid generation algorithm is such that skewness diminishes as $O(h^p)$ then the 
order of accuracy of the DT gradient is at least $\min\{p,1\}$.

\subsubsection*{Unstructured grids}

Unstructured grids usually consist of triangles (or tetrahedra in 3D), or of polygonal (polyhedral 
in 3D) cells which are also formed by a triangulation process. They are typically constructed using 
algorithms that are based on geometrical principles that do not depend on grid fineness, so that 
there is some self-similarity between coarse and fine grids and refinement does not reduce the 
skewness, i.e.\ $\|\vf{c}_f - \vf{c}'_f\| / \|\vf{N}_f - \vf{P}\| = O(1)$. This means that 
$\vf{c}_f - \vf{c}'_f = \vf{O}(h)$ and therefore instead of Eq.\ \eqref{eq: phi_c vs phi_c_bar 
structured} we now have
\begin{equation} \label{eq: phi_c vs phi_c_bar unstructured}
 \phi(\vf{c}_f) \;=\; \phi(\vf{c}'_f) + O(h) \;=\; (\bar{\phi}(\vf{c}'_f) + O(h^2)) + O(h) \;=\; 
\bar{\phi}(\vf{c}'_f) + O(h) 
\end{equation}
Substituting this into the exact equation \eqref{eq: gauss gradient halfway} we get
\begin{align*}
 \nabla \phi (\vf{P}) \;&=\; \frac{1}{\Omega_P} \sum_{f=1}^F \left( \bar{\phi}(\vf{c}'_f) + O(h) 
\right)\, S_f\, \vf{n}_f \;+\; 
                             \vf{O}(h)
\\[0.2cm] 
                        &=\; \frac{1}{\Omega_P} \sum_{f=1}^F \bar{\phi}(\vf{c}'_f) \, S_f \, 
\vf{n}_f \;+\; 
                             \frac{1}{\Omega_P} \sum_{f=1}^F  O(h) \, S_f \, \vf{n}_f \;+\; 
\vf{O}(h)
\end{align*}
which, considering that $\Omega_P = O(h^2)$ and $S_f = O(h)$, means that
\begin{equation} \label{eq: gauss gradient exact unstructured}
 \nabla \phi (\vf{P}) \;=\; \frac{1}{\Omega_P} \sum_{f=1}^F \bar{\phi}(\vf{c}'_f) \, S_f \, 
\vf{n}_f \;+\; \vf{O}(1) \;+\; \vf{O}(h)
\end{equation}
Again, the first term on its right-hand side is the approximate gradient $\nabla^{\mathrm{d0}}$, 
Eq.\ \eqref{eq: gauss gradient 0}. But this time, unlike in the structured grid case, the 
truncation error of $\nabla^{\mathrm{d0}}$ is $\vf{\tau} \equiv \nabla^{\mathrm{d}0}\phi(\vf{P}) - 
\nabla\phi(\vf{P}) = \vf{O}(1)$. Nor do the leading terms of the face contributions to the 
truncation error cancel out to increase its order, because on unstructured grids the orientations 
and sizes of the faces are unrelated. This is an unfortunate result, because it means that the 
approximation \eqref{eq: gauss gradient 0} is zeroth-order accurate, as the error $\vf{O}(1)$ does 
not decrease with grid refinement, but instead $\nabla^{\mathrm{d}0} \phi (\vf{P})$ converges to a 
value that is not equal to $\nabla\phi(\vf{P})$ (see Section \ref{ssec: results refined} for an 
example).

Acknowledging that the lack of accuracy is due to the bad representation of the $\phi(\vf{c}_f)$ 
values in Eq.\ \eqref{eq: gauss gradient halfway} by $\bar{\phi}(\vf{c}'_f)$, application of the 
formula \eqref{eq: gauss gradient 0} is often followed by a ``corrector step'' where, instead of 
the values $\bar{\phi}(\vf{c}'_f)$, the ``improved'' values $\hat{\phi}(\vf{c}_f)$ are used, defined 
as
\begin{equation} \label{eq: phi cor}
 \hat{\phi}(\vf{c}_f) \;\equiv\; \bar{\phi}(\vf{c}'_f)
                      \;+\;      \overline{\nabla^{\mathrm{d}0} \phi}(\vf{c}'_f) \!\cdot\! 
(\vf{c}_f - \vf{c}'_f) 
\end{equation}
where $\overline{\nabla^{\mathrm{d}0} \phi}(\vf{c}'_f)$ is obtained by linear interpolation (Eq.\ 
\eqref{eq: linear interpolation}) between $\nabla^{\mathrm{d}0} \phi (\vf{P})$ and 
$\nabla^{\mathrm{d}0} \phi (\vf{N}_f)$ at point $\vf{c}'$ (these were calculated in the previous 
step \eqref{eq: gauss gradient 0}, the ``predictor'' step). This results in an approximation 
which is hopefully more accurate than \eqref{eq: gauss gradient 0}:
\begin{equation} \label{eq: gauss gradient 1}
 \nabla^{\mathrm{d}1} \phi (\vf{P}) \;\equiv\; \frac{1}{\Omega_P} \sum_{f=1}^F \hat{\phi}(\vf{c}_f) 
\, S_f \, \vf{n}_f 
\end{equation}
The values $\hat{\phi}(\vf{c}_f)$ are expected to be better approximations to $\phi(\vf{c}_f)$ than 
$\bar{\phi}(\vf{c}'_f)$ are, because Eq.\ \eqref{eq: phi cor} tries to account for skewness by 
mimicking the Taylor series expansion
\begin{equation} \label{eq: phi_c from phi_c'}
\phi(\vf{c}_f) \;=\; \phi(\vf{c}'_f) \;+\; \nabla\phi(\vf{c}'_f) \cdot (\vf{c}_f - \vf{c}'_f) \;+\; 
                     O(h^2)
\end{equation}
Unfortunately, since in Eq.\ \eqref{eq: phi cor} only a crude approximation of 
$\nabla\phi(\vf{c}'_f)$ is used, namely $\overline{\nabla^{\mathrm{d}0} \phi}(\vf{c}'_f) = 
\nabla\phi(\vf{c}'_f) + \vf{O}(1)$, what we get by subtracting Eq.\ \eqref{eq: phi cor} from Eq.\ 
\eqref{eq: phi_c from phi_c'} is $\hat{\phi}(\vf{c}_f) = \phi(\vf{c}_f) + O(h)$ which may have 
greater accuracy but not greater \textit{order} of accuracy than the previous estimate 
$\bar{\phi}(\vf{c}'_f) = \phi(\vf{c}_f) + O(h)$, Eq.\ \eqref{eq: phi_c vs phi_c_bar unstructured}. 
Substituting this into Eq.\ \eqref{eq: gauss gradient halfway} we arrive again at an equation 
similar to \eqref{eq: gauss gradient exact unstructured} which shows that the error of the 
approximation \eqref{eq: gauss gradient 1} is also of order $\vf{O}(1)$.

Further correction steps may be applied in the same manner; a fixed finite number of such steps may 
increase the accuracy, but the order of accuracy with respect to grid refinement will remain zero.
But if this procedure is repeated until convergence to an operator $\nabla^{{\mathrm{d}}\infty}$, 
say, then $\nabla^{{\mathrm{d}}\infty}$ would simultaneously satisfy both Eqs.\ \eqref{eq: phi cor} 
and \eqref{eq: gauss gradient 1}, or combined in a single equation:
\begin{equation} \label{eq: gauss gradient infinity}
 \nabla^{\mathrm{d}\infty} \phi(\vf{P}) \;-\; 
 \frac{1}{\Omega_P} \sum_{f=1}^F \overline{\nabla^{\mathrm{d}\infty} \phi}(\vf{c}'_f) \!\cdot\! 
   (\vf{c}_f - \vf{c}'_f) \, S_f \, \vf{n}_f
 \;=\;
 \frac{1}{\Omega_P} \sum_{f=1}^F \bar{\phi}(\vf{c}'_f) \, S_f \, \vf{n}_f
\end{equation}
where we have moved all terms involving the gradient to the left-hand side. An analogous equation 
can be derived for the exact gradient, from Eqs.\ \eqref{eq: gauss gradient halfway}, \eqref{eq: 
phi_c from phi_c'}, and \eqref{eq: linear interpolation}:
\begin{equation} \label{eq: gauss gradient infinity exact}
 \nabla \phi(\vf{P}) \;-\; 
 \frac{1}{\Omega_P} \sum_{f=1}^F \overline{\nabla \phi}(\vf{c}'_f) \!\cdot\! 
   (\vf{c}_f - \vf{c}'_f) \, S_f \, \vf{n}_f
 \;=\;
 \frac{1}{\Omega_P} \sum_{f=1}^F \bar{\phi}(\vf{c}'_f) \, S_f \, \vf{n}_f \;+\;
 \vf{O}(h)
\end{equation}
By subtracting Eq.\ \eqref{eq: gauss gradient infinity exact} from Eq.\ \eqref{eq: gauss gradient 
infinity} we get an expression for the truncation error $\vf{\tau} \equiv 
\nabla^{\mathrm{d}\infty}\phi - \nabla \phi$:
\begin{equation} \label{eq: gauss gradient infinity tau}
 \vf{\tau}(\vf{P}) \;-\; 
 \frac{1}{\Omega_P} \sum_{f=1}^F \overline{\vf{\tau}}(\vf{c}'_f) \!\cdot\! 
   (\vf{c}_f - \vf{c}'_f) \, S_f \, \vf{n}_f
 \;=\;
 \vf{O}(h)
\end{equation}
The left-hand side is a linear combination of not only $\vf{\tau}(\vf{P})$ but also 
$\vf{\tau}(\vf{N}_f)$ at all neighbour points. Since $\Omega_P = O(h^2)$, $S_f = O(h)$ and 
$\vf{c}_f - \vf{c}'_f = O(h)$, the coefficients of this linear combination are $O(1)$ i.e.\ they do 
not depend on the grid fineness but only on the grid geometry (skewness, unevenness etc.). 
Expression \eqref{eq: gauss gradient infinity tau} suggests that $\vf{\tau} = \vf{O}(h)$, i.e.\ 
that $\nabla^{\mathrm{d}\infty}$ is first-order accurate, although to ascertain this one would have 
to assemble Eqs.\ \eqref{eq: gauss gradient infinity tau} for all grid cells in a large linear 
system $A\tau = b \Rightarrow \tau = A^{-1} b$ (where $\tau$ stores the $x$- and $y$- (and $z$-, in 
3D) components of $\vf{\tau}$ at all cell centres and $b = O(h)$), select a suitable matrix norm 
$\|\cdot\|$, and ensure that $\|A^{-1}\|$ remains bounded as $h \rightarrow 0$.

So, iterating Eqs.\ \eqref{eq: phi cor} and \eqref{eq: gauss gradient 1} until convergence one 
would hope that a first-order accurate gradient $\nabla^{\mathrm{d}\infty}$ would be obtained. As 
reported in \cite{Betchen_2010}, convergence of this procedure is not guaranteed and 
underrelaxation may be needed, leading to a large number of required iterations and great 
computational cost. In practice, the desired accuracy will be reached with a finite number of 
iterations, but this number is not known a priori and must be determined by trials; furthermore, in 
order to maintain the property that grid refinement improves the accuracy, the number of iterations 
should increase with grid refinement in order to avoid accuracy stagnation. Considering that the 
main rival of the DT gradient, the LS gradient, costs approximately as much as the DT gradient with 
a single corrector step (see Sec.\ \ref{ssec: in-house tests}), it can be seen that this iterative 
procedure is impractical and may end up consuming most of the computational time of a FVM solver. 
Alternatively, one could directly solve all the equations \eqref{eq: gauss gradient infinity} 
simultaneously in a large linear system; this is also proposed and tested in \cite{Betchen_2010}, 
along with a second-order accurate method that additionally computes the Hessian matrix. Obviously, 
this approach also involves a very large computational cost.

However, we would like to propose here a way to avoid the extra cost of the iterative procedure 
for the gradient. The FVM generally employs outer iterations to solve a PDE (e.g.\ SIMPLE 
iterations, as opposed to the inner iterations of the linear solvers employed). Outer iterations are 
necessary when the equations solved are non-linear, but may be used also when solving linear 
problems if some terms are treated with the deferred-correction approach. FVMs that use gradient 
schemes such as those presently discussed are almost always iterative, with the gradients computed 
using the values of the dependent variable from the previous outer iteration. The idea is then to 
exploit these outer iterations, dividing the ``gradient iterations'' among them: at each outer 
iteration $n$ a single ``gradient iteration'' is performed according to
\begin{equation} \label{eq: gauss outer iterative scheme}
 \nabla^{\mathrm{d}n} \phi(\vf{P}) \;=\; 
 \frac{1}{\Omega_P} \sum_{f=1}^F \left( 
 \overline{\phi^{n-1}}(\vf{c}'_f) \;+\; 
 \overline{\nabla^{\mathrm{d}(n-1)} \phi}(\vf{c}'_f) \!\cdot\! (\vf{c}_f - \vf{c}'_f)
 \right) \, S_f \, \vf{n}_f
\end{equation}
which means that one has to store not only the solution of the previous outer iteration 
$\phi^{n-1}$ but also the gradient of that iteration $\nabla^{\mathrm{d}(n-1)}\phi$, which many 
codes do already. Since only one calculation is performed per outer iteration, the cost of this 
procedure is almost as low as that of the uncorrected DT gradient \eqref{eq: gauss gradient 0}, but 
it converges to the first-order accurate $\nabla^{\mathrm{d}\infty}$ gradient operator instead of 
the zeroth-order accurate $\nabla^{\mathrm{d}0}$. For explicit time-dependent FVM methods, Eq.\ 
\eqref{eq: gauss outer iterative scheme} could be applied with $n-1$ denoting the previous time 
step; in the very first time step it may be necessary to perform several ``gradient iterations'' to 
obtain a consistent gradient to begin with.

Iterative schemes that are different from \eqref{eq: gauss outer iterative scheme} can also be 
devised; a similar iterative technique is used in \cite{Wang_2017} in order to determine the 
coefficients of reconstruction polynomials without solving large linear systems. There, Jacobi, 
Gauss-Seidel, and SOR-type iterations are applied, with the latter found to be the most efficient. 
The present scheme \eqref{eq: gauss outer iterative scheme} is reminiscent of the Jacobi iterative 
procedure (although not exactly equivalent, since $\nabla \phi(\vf{P})$ appears also at the 
right-hand side, within the averaged gradient, evaluated from the previous iteration). A scheme 
reminiscent of the Gauss-Seidel iterative procedure would be one where in the right-hand side of 
Eq.\ \eqref{eq: gauss outer iterative scheme} new values of the gradient at neighbour cells 
(calculated at the current outer iteration) are used whenever available instead of using the old 
values. However, in the present work only the scheme \eqref{eq: gauss outer iterative scheme} is 
tested in Section \ref{ssec: in-house tests}.

At this point we would like to make a brief comment concerning the calculation of the gradient at 
boundary cells. The above analysis has assumed that at boundary faces the values of $\phi$ at the 
face centres are available. However, this is not always the case; in situations where these values 
are not available, they must be approximated to at least second-order accuracy as otherwise the DT 
gradient will be inconsistent, as the preceding analysis has shown. For example, in problems with 
Neumann boundary conditions the directional derivative normal to the boundary, $g$ say, is given 
rather than the boundary values. If face 5 of Fig.\ \ref{fig: unstructured grid} belongs to such a 
boundary, then expressing $\phi(\vf{P})$ in a Taylor series about point $\vf{c}'_5$ gives
\begin{equation} \label{eq: neumann bc}
 \phi(\vf{c}'_5) \;=\; \phi(\vf{P}) \;+\; g \|\vf{c}'_5 - \vf{P}\| \;+\; O(h^2)
\end{equation}
(where $g$ is measured in the direction pointing out of the domain). Equation \eqref{eq: neumann 
bc} provides a second-order accurate approximation for $\phi(\vf{c}'_5)$, but only a first-order 
accurate approximation for $\phi(\vf{c}_5)$. Therefore, if we just use the value 
$\bar{\phi}(\vf{c}'_5) \equiv \phi(\vf{P}) \;+\; g \|\vf{c}'_5 - \vf{P}\|$ in the DT gradient 
formula then we will get a zeroth-order accurate gradient. Note that this will hold even for 
structured grids if the grid lines intersect the boundary at an angle ($\vf{c}'_5 \neq \vf{c}_5$). 
It is not difficult to show that a second order approximation is $\phi(\vf{c}_5) \approx 
\bar{\phi}(\vf{c}'_5) + \nabla\phi(\vf{P}) \cdot (\vf{c}_5 - \vf{c}'_5)$. This introduces $\nabla 
\phi(\vf{P})$ also in the right-hand side of the gradient expression \eqref{eq: gauss gradient 
halfway}, which poses no problem for the iterative procedure \eqref{eq: gauss outer iterative 
scheme}.

\begin{figure}[tb]
 \centering
 \includegraphics[scale=0.75]{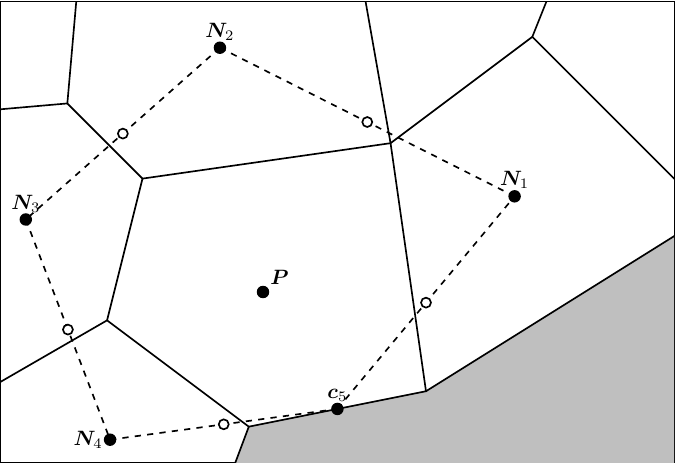}
 \caption{The divergence theorem method can approximate $\nabla \phi(\vf{P})$ to first-order 
accuracy on arbitrarily irregular grids if applied to the auxiliary cell bounded by dashed lines, 
rather than on the actual cell $P$. This cell is formed by joining the centroids of the 
neighbouring cells and of the boundary faces of cell $P$ by straight line segments. The open 
circles denote the centroids of these segments.}
 \label{fig: dual cell}
\end{figure}

Finally, we note that alternative schemes can be derived from the Gauss divergence theorem that are 
consistent on unstructured grids and are worth mentioning, although the present work focuses on the 
standard method. For example, consider the application of the divergence theorem method not to the 
actual cell $P$ but to the auxiliary cell marked by dashed lines in Fig.\ \ref{fig: dual cell}. The 
endpoints of each face of this cell are either cell centroids or boundary face centroids, and so 
the value of $\phi$ can be computed at any point on the face to second-order accuracy, by linear 
interpolation between the two endpoints. Therefore, in Eq.\ \eqref{eq: gauss gradient halfway} the 
values $\phi(\vf{c}_f)$ ($\vf{c}_f$ being the face centroids of the \textit{auxiliary} cell, marked 
by empty circles in Fig.\ \ref{fig: dual cell}) are calculated to second-order accuracy rather than 
first-order, leading to first-order accuracy of the computed gradient. Second-order accuracy is 
inhibited also by the fact that the midpoint integration rule (Eq.\ \eqref{eq: midpoint rule cell}) 
requires that the gradient be computed at the centroid of the auxiliary cell, which now does not, 
in general, coincide with point $\vf{P}$ (the situation changes if $\vf{P}$ and the centroid of the 
auxiliary cell tend to coincide with grid refinement). This method, along with a more complex 
variant, is mentioned in \cite{Barth_1989}; it is also tested in Section \ref{ssec: in-house tests}.
Yet another method would be to use cell $P$ itself rather than an auxiliary cell, but, similarly to 
the previous method, to calculate the values $\phi(\vf{c}_f)$ in Eq.\ \eqref{eq: gauss gradient 
halfway} by linear interpolation from the values at the face endpoints (vertices) rather than from 
the values at the cell centroids straddling the face. An extra step must therefore precede where 
$\phi$ is approximated at the cell vertices to second-order accuracy from its values at the cell 
centroids. This adds to the computational cost and furthermore requires of the grid data structures 
to contain lists relating each vertex to its surrounding cells. More information on this method and 
further references can be found in \cite{Diskin_2011, Moukalled_2016}.

\section{Gradient calculation using least squares minimisation}
\label{sec: least squares}


The starting point of the ``least-squares'' (LS) method for calculating the gradient of a quantity 
$\phi$ at the centroid of a cell $P$ is the expression of the values of $\phi$ at the centroids of 
all neighbouring cells as Taylor series expansions about the centroid of $P$. For convenience, we 
decompose the position vectors in Cartesian coordinates as  $\vf{P} = (x_0,y_0)$, $\vf{N}_f = 
(x_f,y_f)$ and denote $\Delta x_f \equiv x_f - x_0$, $\Delta y_f \equiv y_f - y_0$. Then, for each 
neighbour $f$ it follows from the Taylor expansion that
\begin{equation} \label{eq: ls Taylor}
 \phi(\vf{N}_f) - \phi(\vf{P}) \;=\; \nabla\phi(\vf{P})\cdot(\vf{N}_f-\vf{P}) \;+\; \varsigma_f
\end{equation}
where
\begin{equation} \label{eq: varsigma_f}
 \varsigma_f \;=\; \frac{1}{2} \phi_{.xx} (\Delta x_f)^2 \;+\; 
                   \phi_{.xy} \, \Delta x_f \, \Delta y_f \;+\;
                   \frac{1}{2} \phi_{.yy} (\Delta y_f)^2 \;+\; O(h^3)
\end{equation}
is an associated truncation error (the second derivatives are evaluated at $\vf{P}$). The basic 
idea of the method is to drop the unknown $\varsigma_f$ terms and solve the remaining linear system 
by least squares since, in general, the number of equations ($F$ = number of neighbouring cells) 
will be greater than the number of unknowns (two, $\phi_{\!.x}(\vf{P})$ and $\phi_{\!.y}(\vf{P})$). 
Of course some of the ``neighbours'' may be boundary face points such as the centroid $\vf{c}_5$ in 
Fig.\ \ref{fig: unstructured grid}, or, in the case of Neumann boundary conditions, the projection 
point $\vf{c}'_5$. In the latter case $\phi(\vf{c}'_5)$ can be approximated by Eq.\ \eqref{eq: 
neumann bc}. Furthermore, unlike for the DT method, it is now easy to incorporate additional, more 
distant cells as neighbours, a practice that may prove advantageous in some cases\footnote{It may 
even be necessary in some cases. For example, some of the triangular cells at the corners of the 
grid in Fig.\ \ref{sfig: grid diffusion netgen} have only one neighbour cell; if the values of 
$\phi$ are not available at the boundary faces (e.g.\ if $\phi$ is the pressure in incompressible 
flows) then at least one more distant cell must be used, because the LS calculation requires at 
least two neighbour points.} (see e.g.\ \cite{Diskin_2011}), although in the present work we will 
restrict ourselves to using only the cells that share a face with $P$.

It is advantageous to first multiply each equation $f$ by a suitably selected weight $w_f$ 
(weighted least squares). Then the system of all these equations can be written in matrix form as
\newcommand*{\weightsM}{
 \begin{bmatrix}
    w_1     &  \multicolumn{3}{c}{\text{\kern 1em\smash{\raisebox{-2.5ex}{\Huge 0}}}} \\
            &  w_2  &          &       \\
            &       &  \ddots  &       \\
  \multicolumn{3}{c}{\text{\kern -1em\smash{\raisebox{0ex}{\Huge 0}}}}  &  w_F
 \end{bmatrix}}
\begin{equation} \label{eq: ls system weighted}
 \underbrace{
 \weightsM
 }_W
 \!\cdot\!
 \underbrace{
 \begin{bmatrix}
  \Delta\phi_1 \\
  \Delta\phi_2\\
  \vdots \\
  \Delta\phi_F
 \end{bmatrix}
 }_b
 \;=\;
 \underbrace{
 \weightsM
 }_W
 \!\cdot\!
 \left(
 \smashedunderbrace{
 \begin{bmatrix}
  \Delta x_1 & \Delta y_1 \\
  \Delta x_2 & \Delta y_2 \\
  \vdots     & \vdots     \\
  \Delta x_F & \Delta y_F
 \end{bmatrix}
 }{A}
 \!\cdot\!
 \smash{\underbrace{
 \begin{bmatrix}
  \phi_{\!.x}(\vf{P}) \\[0.15cm]
  \phi_{\!.y}(\vf{P})
 \end{bmatrix}
 \vphantom{\weightsM}}_z}
 \;+\;
 \smashedunderbrace{
 \begin{bmatrix}
  \varsigma_1 \\
  \varsigma_2 \\
  \vdots \\
  \varsigma_F
 \end{bmatrix}
 }{\varsigma}
 \right)
\end{equation}
where $\Delta \phi_f \equiv \phi(\vf{N}_f) - \phi(\vf{P})$. The alternative matrix notation $Wb = 
W(Az+\varsigma)$ is also introduced above, to facilitate the discussions that follow. The above 
equations are exact, with the vector $\varsigma$ ensuring that the system \eqref{eq: ls system 
weighted} has a unique solution, which is independent of the weights, despite having more equations 
than unknowns. The redundant equations are just linear combinations of the non-redundant ones. But 
once $\varsigma$ is dropped, the system will, in general, have no solution.

Systems with no solution can be solved in the ``least squares'' sense \cite{Trefethen_1997, 
Strang_2006}. In matrix notation, if a linear system $Az=b$ with more equations than unknowns cannot 
be solved because the vector $b$ does not lie in the column space of the matrix $A$ then the best 
that can be done is to find the vector $z$ that minimises the error $b-Az$. This error will be 
minimised when its projection onto the column space of $A$ is zero, i.e.\ when $b-Az$ is 
perpendicular to each column of $A$, or $A^{\mathrm{T}}(b-Az) = 0 \Rightarrow A^{\mathrm{T}}Az = 
A^{\mathrm{T}}b$ (where $\mathrm{T}$ denotes the transpose). This latter system (called the 
\textit{normal equations}) has a unique solution provided that $A$ has independent columns. The 
solution $z$ is called the least squares solution because it minimises the $L_2$ norm $\|b - Az\|$, 
i.e.\ the sum of the squares of the individual components of the error $e \equiv b - Az$ (note that 
each individual error component $e_i$ is the error of the corresponding equation $i$, $e_i = b_i - 
\sum_j a_{ij} z_j$).

In the weighted case $WAz = Wb$ the matrix $A$ is replaced by the product $WA$ and $b$ by $Wb$, so 
that the corresponding solution to the normal equations is
\begin{equation} \label{eq: WLS solution}
z \;=\; (A^T W^{\mathrm{T}} W A)^{-1} A^{\mathrm{T}} W^{\mathrm{T}} W b
\end{equation}
Now the quantity minimised is the norm of the error of the weighted system, $W(b-Az) = We$. Thus, 
if equation $i$ is assigned a larger weight than equation $j$ then the method will prefer to make 
$e_i$ small at the expense of $e_j$.

For now, we proceed with the least squares methodology on Eq.\ \eqref{eq: ls system weighted} 
without discarding $\varsigma$ in order to assess its impact on the error. We therefore first 
left-multiply the system \eqref{eq: ls system weighted} by $(WA)^\mathrm{T}$, and then solve 
it to obtain
\begin{equation} \label{eq: WLS matrix notation}
 z \;=\; (A^{\mathrm{T}} W^{\mathrm{T}} W A)^{-1} A^{\mathrm{T}} W^{\mathrm{T}} W b
   \;+\; (A^{\mathrm{T}} W^{\mathrm{T}} W A)^{-1} A^{\mathrm{T}} W^{\mathrm{T}} W \varsigma
\end{equation}
The matrix $A^{\mathrm{T}} W^{\mathrm{T}} W A$ is invertible provided that $WA$ has independent 
columns, which requires that $A$ have independent columns because $W$ is diagonal. The two columns 
of $A$ will be independent if there are also two independent rows, i.e.\ vectors $\vf{N}_f - 
\vf{P}$, so we need at least two neighbours $\vf{N}_f$ to lie at different directions with respect 
to $\vf{P}$, which is normally the case (in three dimensions $A$ has three columns and three 
independent neighbour directions are needed). Then, substituting for $A$, $W$, $z$, $b$ and 
$\varsigma$ from Eq.\ \eqref{eq: ls system weighted} into Eq.\ \eqref{eq: WLS matrix notation} and 
performing the matrix multiplications and inversions, we arrive at
\begin{equation} \label{eq: WLS final exact}
 z \;\equiv\; \nabla \phi(\vf{P})
 \;=\; \nabla^{\mathrm{ls}} \phi(\vf{P}) \;+\; \tau
\end{equation}
where $\nabla^{\mathrm{ls}} \phi(\vf{P})$ is the first term on the right side of Eq.\ \eqref{eq: 
WLS matrix notation},
\begin{equation} \label{eq: WLS final}
 \nabla^{\mathrm{ls}} \phi(\vf{P})
 \;\equiv\;
 \frac{1}{\mathrm{D}}
 \underbrace{
 \begin{bmatrix}
   \sum\limits_{f=1}^F (\Delta y_f)^2 w_f^2         &  -\sum\limits_{f=1}^F \Delta x_f \Delta y_f 
w_f^2 \\[0.4cm]
  -\sum\limits_{f=1}^F \Delta x_f \Delta y_f w_f^2  &   \sum\limits_{f=1}^F (\Delta x_f)^2 w_f^2
 \end{bmatrix}
 }_M
 \cdot
 \underbrace{
 \begin{bmatrix}
  \sum\limits_{f=1}^F \Delta x_f \Delta \phi_f w_f^2 \\[0.4cm]
  \sum\limits_{f=1}^F \Delta y_f \Delta \phi_f w_f^2 
 \end{bmatrix}
 }_{\beta_b}
\end{equation}
with $(1/D)M = (A^{\mathrm{T}} W^{\mathrm{T}} W A)^{-1}$ and $\beta_b = A^{\mathrm{T}} 
W^{\mathrm{T}} W b$, and $\tau$ is the second term:
\begin{equation} \label{eq: WLS e}
 \tau
 \;=\;
 \frac{1}{\mathrm{D}} M
 \cdot
 \begin{bmatrix}
  \sum\limits_{f=1}^F \Delta x_f \varsigma_f w_f^2   \;&\;
  \sum\limits_{f=1}^F \Delta y_f \varsigma_f w_f^2 
 \end{bmatrix}^{\mathrm{T}}
 \;=\; \frac{1}{D} M \, \beta_{\varsigma}
\end{equation}
with $\beta_{\varsigma} = A^{\mathrm{T}} W^{\mathrm{T}} W \varsigma$. In the 2D case, $D = |M|$, 
the determinant of $M$. Equation \eqref{eq: WLS final exact} gives the exact gradient $\nabla 
\phi(\vf{P})$, but since the truncation error $\tau$ is unknown we drop it and use the expression 
\eqref{eq: WLS final} alone as the approximate ``least squares'' gradient. In the 3D case the 
explicit formula is more involved than Eq.\ \eqref{eq: WLS final}, but again the least squares 
gradient is obtained by solving Eq.\ \eqref{eq: WLS solution}. If, just for the purposes of 
discussing the 3D case, we denote the Cartesian components of the displacement vectors as $\vf{N}_f 
- \vf{P} = (\Delta x_f^1, \Delta x_f^2, \Delta x_f^3)$, then the $(i,j)$ entry of the matrix 
$A^{\mathrm{T}} W^{\mathrm{T}} W A$ is $\sum_f \Delta x_f^i \Delta x_f^j w_f^2$ (it is a symmetric 
matrix) and the $i$-th component of the vector $A^{\mathrm{T}} W^{\mathrm{T}} W b$ is $\sum_f \Delta 
x_f^i \Delta \phi_f w_f^2$. The normal equations \eqref{eq: WLS solution} are sometimes 
ill-conditioned, and in such cases it is better to solve the least-squares system by $QR$ 
factorisation of $A$ \cite{Trefethen_1997}.

The above derivation has produced the explicit expression \eqref{eq: WLS e} for the truncation 
error $\tau$. To analyse it further, we can assume that all the weights share the same 
dependency on the grid spacing, namely $w_f = O(h^q)$ for some real number $q$ (independent of 
$f$), as is the usual practice. Then the factors of Eq.\ \eqref{eq: WLS e} have the following 
magnitudes: Since $\Delta x_f$ and $\Delta y_f$ are $O(h)$ the coefficients of $M$ have magnitude 
$O(h^{2+2q})$. Consequently, $M$ being $2\times2$, its determinant will have magnitude 
$(O(h^{2+2q}))^2 \Rightarrow 1/D = O(h^{-2(2+2q)})$. Finally, considering that $\varsigma_f = 
O(h^2)$, the components of $\beta_{\varsigma}$ are of $O(h^{2q+3})$. Multiplying all these 
together, Eq.\ \eqref{eq: WLS e} shows that $\tau = O(h)$, independently of $q$. This is not 
surprising, given that the approximation is based on Eq.\ \eqref{eq: ls Taylor} which assumes a 
linear variation of $\phi$ in the vicinity of point $\vf{P}$. So, $\nabla^{\mathrm{ls}}$ is at least 
first order accurate, even on grids of arbitrary geometry.

The order of accuracy of $\nabla^{\mathrm{ls}}$ may be higher than one if some cancellation occurs 
between the components of $\tau$ for certain grid configurations, similarly to the DT operator. In 
particular, a tedious but straightforward calculation shows that when applied to a parallelogram 
grid such as that shown in Fig.\ \ref{fig: parallelogram grid}, $\nabla^{\mathrm{ls}}$ again reduces 
to the second-order accurate formula \eqref{eq: grad parallelogram}, provided only that the weights 
of parallel faces are equal: $w_1 = w_3$ and $w_2 = w_4$. This will hold due to symmetry if the 
weights are dependent only on the grid geometry. Of course, as for the DT gradient, this has the 
consequence that the LS gradient is second-order accurate also on smooth curvilinear structured 
grids.

\subsubsection*{The choice of weights}

According to the preceding analysis, the least squares method is first-order accurate on grids of 
arbitrary geometry and second-order accurate on smooth structured grids. This holds irrespective of 
the choice of weights, and in fact it holds even for the unweighted method ($w_f = 1$). The question 
then arises of whether a suitable choice of weights can offer some advantage.

The weights commonly used are of the form $w_f = (\Delta r_f)^{-q}$ where $\Delta r_f = \| \vf{N}_f 
- \vf{P} \|$ is the distance between the two cell centres and $q$ is an integer, usually chosen as
$q=1$ \cite{Barth_1991, Mavriplis_2003, Moukalled_2016} or $q=2$ \cite{Ollivier_2002, Correa_2011}. 
The unweighted method amounts to $q = 0$. As noted, the least squares method finds the approximate 
gradient that minimises $\|W(b-Az)\|$, which for the various choices of $q$ amounts to minimising:
\begin{align}
 \label{eq: minimised q=0}
 q = 0:   \qquad &   \sum_f \left[ \Delta r_f \left(
                     \frac{\Delta \phi_f}{\Delta r_f} \;-\;
                     \nabla^{\mathrm{ls}} \phi(\vf{P}) \cdot \vf{d}_f
                     \right) \right]^2
\\[0.2cm]
 \label{eq: minimised q=1}
 q = 1:   \qquad &   \sum_f \left(
                     \frac{\Delta \phi_f}{\Delta r_f} \;-\;
                     \nabla^{\mathrm{ls}} \phi(\vf{P}) \cdot \vf{d}_f
                     \right)^2
\\[0.2cm]
 \label{eq: minimised q=2}
 q = 2:   \qquad &   \sum_f \left[ \frac{1}{\Delta r_f} \left(
                     \frac{\Delta \phi_f}{\Delta r_f} \;-\;
                     \nabla^{\mathrm{ls}} \phi(\vf{P}) \cdot \vf{d}_f
                     \right) \right]^2
\end{align}
where $\vf{d}_f = (\vf{N}_f - \vf{P}) / \Delta r_f$ is the unit vector in the direction from 
$\vf{P}$ to $\vf{N}_f$, and so $\nabla^{\mathrm{ls}} \phi \cdot \vf{d}_f$ is the least squares 
directional derivative $\partial \phi / \partial r_f$ in that direction. In the $q = 1$ case, 
expression \eqref{eq: minimised q=1} shows that the least squares procedure shows no preference in 
trying to set the directional derivative in each neighbour direction $f$ equal to the finite 
difference $\Delta \phi_f / \Delta r_f$. On the other hand, in the unweighted method \eqref{eq: 
minimised q=0} the discrepancies between the directional derivatives and the finite differences are 
weighted by the distances $\Delta r_f$ so that the method prefers to reduce $\Delta \phi / \Delta 
r_f - \nabla^{\mathrm{ls}} \phi \cdot \vf{d}_f$ along the directions of the distant neighbours at 
the expense of the directions of the closer neighbours. The exact opposite holds for the $q = 2$ 
case \eqref{eq: minimised q=2} where the discrepancies are weighted by $1 / \Delta r_f$ so that the 
result is determined mostly by the close neighbours. Intuitively, this latter choice seems more 
reasonable as the linearity of the variation of $\phi$ is lost as one moves away from $\vf{P}$ and 
thus at distant neighbours the finite differences $\Delta \phi / \Delta r_f$ are less accurate 
approximations of the directional derivatives. Thus it is not surprising that usually the $q > 0$ 
methods outperform the unweighted method.

These are the common weight choices, but it so happens that the particular non-integer exponent $q = 
3/2$ confers enhanced accuracy compared to the other $q$ choices under special but not uncommon 
circumstances. This fact does not appear to be widely known in the literature; we have seen it only 
briefly mentioned in \cite{Muzaferija_1994, Shima_2010}. So, consider the vector $\beta_{\varsigma}$ 
in the expression \eqref{eq: WLS e} for the error, and substitute for $\varsigma$ from Eq.\ 
\eqref{eq: varsigma_f}. Then the first component of $\beta_{\tau}$ becomes 
\begin{align}
\nonumber
 \sum_f \Delta x_f \, \varsigma_f \, w_f^2 \;&=\; 
 \frac{\phi_{.xx}}{2} \sum_f (\Delta x_f)^3 w_f^2 \;+\;
 \phi_{.xy} \sum_f (\Delta x_f)^2 \Delta y_f \, w_f^2 \;+\;
 \frac{\phi_{.yy}}{2} \sum_f \Delta x_f \, (\Delta y_f)^2 w_f^2
\\[0.2cm]
\label{eq: beta_tau component 1}
 &+\; \sum_f O(h^4) \, w_f^2
\end{align}
Now, if two neighbours, $i$ and $j$ say, lie at opposite directions to point $\vf{P}$ but at the 
same distance then $w_i = w_j$, $\Delta x_i = -\Delta x_j$ and $\Delta y_i = -\Delta y_j$ so that 
their contributions in each of the first three sums in the right hand side of Eq.\ \eqref{eq: 
beta_tau component 1} cancel out. If all neighbour points are arranged in such pairs, like in the 
grid of Fig.\ \ref{fig: parallelogram grid}, then these three sums become zero leaving only the 
$\sum O(h^4) w_f^2$ term. The same holds for the second component of $\beta_{\varsigma}$, so that 
$\beta_{\varsigma} = O(h^{4-2q})$ overall because $w_f = O(h^{-q})$. Then Eq.\ \eqref{eq: WLS e} 
gives $\tau = (1/D) M \beta_{\varsigma} = O(h^2)$  i.e.\ the method is second-order accurate for 
any exponent $q$.

The particular choice $w_f = (\Delta r_f)^{-3/2}$ amounts to dropping the $w_f$'s from the first 
three sums of the right-hand side of Eq.\ \eqref{eq: beta_tau component 1} and replacing therein 
every instance of $\Delta x_f$ by $\Delta x_f / \Delta r_f$ and every instance of $\Delta y_f$ by 
$\Delta y_f / \Delta r_f$. These ratios are precisely $\cos \theta_f$ and $\sin \theta_f$, 
respectively, where $\theta_f$ is the angle that the direction vector $\vf{N}_f - \vf{P}$ makes 
with the horizontal direction. If two neighbours, $i$ and $j$, lie at opposite directions then 
$\theta_i = \theta_j + \pi$ so that $\cos \theta_i = -\cos \theta_j$ and $\sin \theta_i = - \sin 
\theta_j$, and their contributions cancel out in the aforementioned three sums, \textit{irrespective 
of whether these two neighbours lie at equal distances to $\vf{P}$ or not}. The same holds for the 
second component of $\beta_{\varsigma}$. Thus, if all neighbour points are arranged in such pairs 
then again what remains of the right-hand side of Eq.\ \eqref{eq: beta_tau component 1} is only the 
last term and so $\tau = O(h^2)$. For example, the LS gradient with $q=3/2$ is second-order 
accurate at the boundary volume $P$ of Fig.\ \ref{fig: boundary grid}, whereas it is only 
first-order accurate with any other choice of $q$. Furthermore, the same result will hold if the 
neighbours are not arranged in pairs at opposite directions but tend to become so with grid 
refinement. This is the case with smooth structured grids, as shown in Section \ref{sec: 
preliminary}, and therefore the LS gradient with $q = 3/2$ is second-order accurate at boundary 
cells of all smooth structured grids.

Another property of the $q = 3/2$ LS gradient is that it is second order accurate if all neighbour 
points are arranged at equal angles. Unfortunately, this property holds only with more than three 
neighbour points, which limits its usefulness. A proof is provided in Appendix \ref{sec: appendix ls 
angles}. 

Finally, a question that arises naturally is whether full 2nd-order accuracy can be achieved on 
arbitrary grids by allowing non-diagonal entries in the weights matrix. It turns out that this is 
indeed possible, but yields a method that is equivalent to the least squares solution of a system of 
Taylor expansions \eqref{eq: ls Taylor} with terms higher than first-order included. The procedure 
is sketched in Appendix \ref{sec: appendix ls O(2)}. We do not advocate it, because direct solution 
of the system of higher-order Taylor expansions would have the added advantage of solving also for 
the second derivatives -- see e.g.\ \cite{Liu_2017}. An alternative second-order accurate method is 
described in \cite{Vaassen_2008}. The second-order accurate methods are much more expensive than the 
present method.

\section{Numerical tests on the accuracy of the gradient schemes}
\label{sec: tests}

\subsection{One-dimensional tests}
\label{ssec: results 1-d}

The methods are first tested on a one-dimensional problem so as to examine the effect of unevenness, 
isolated from skewness. The derivative of the single-variable function $\phi(x) = \tanh x$ is 
calculated at 101 equispaced points spanning the $x \in [0,2]$ interval. The results are compared 
against the exact solution $\phi_{\!.x} = 1 - (\tanh x)^2$ and the mean absolute error $\sum_i 
|\phi_{\!.x}(x_{(i)}) - \phi_{\!.x}^{\mathrm{ls}} (x_{(i)})| / 101$ is recorded for each method. In 
order to introduce unevenness, the neighbours of point $x_i$ are not chosen from this set of 
equispaced points but are set at $x_{i,f} = x_i + \Delta x_f / 2^r$ where the $\Delta x_f$ belong to 
a predetermined set of displacements $\{\Delta x_f\}_{f=1}^F$ (common to all $x_i$ points) and the 
integer $r$ is the level of grid refinement. For example, if the chosen set of displacements is 
$\{-0.05, 0.1\}$ ($F$ = 2 neighbours), then the derivative at point $x_i$ will be calculated using 
the values of $\phi$ at the three points $\{x_i - 0.05 / 2^r, \; x_i, \; x_i + 0.1 / 2^r\}$. The 
order of accuracy of each method is determined by incrementing the level of refinement $r$.

\begin{figure}[tb]
 \subfigure[$\{\Delta x_f^0\} = \{0.05, 0.10\}$] {\label{sfig: Dx 0.05 0.10} 
    \includegraphics[scale=0.97]{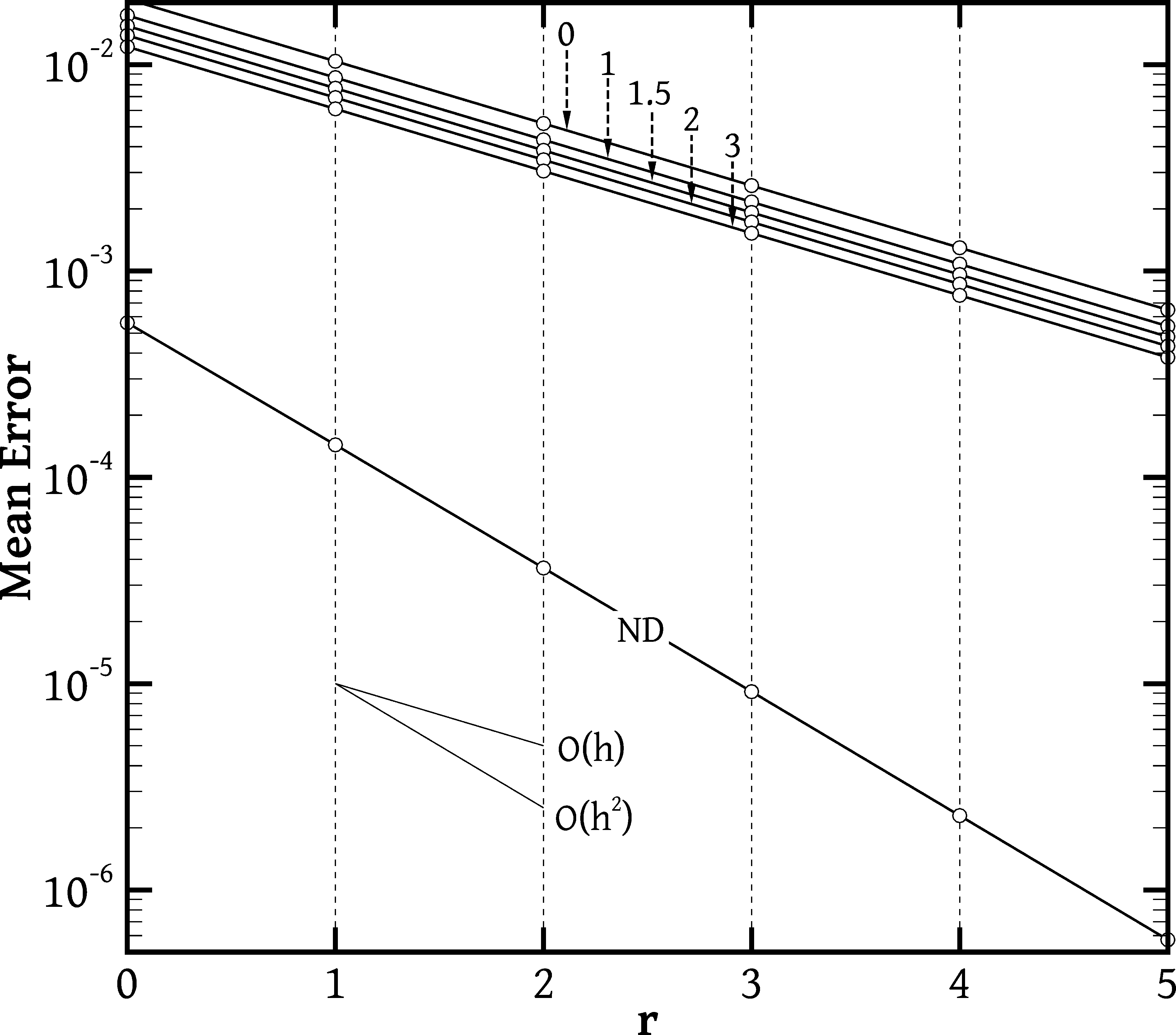}}
 \subfigure[$\{\Delta x_f^0\} = \{-0.10, 0.05\}$] {\label{sfig: Dx -0.10 0.05} 
    \includegraphics[scale=0.97]{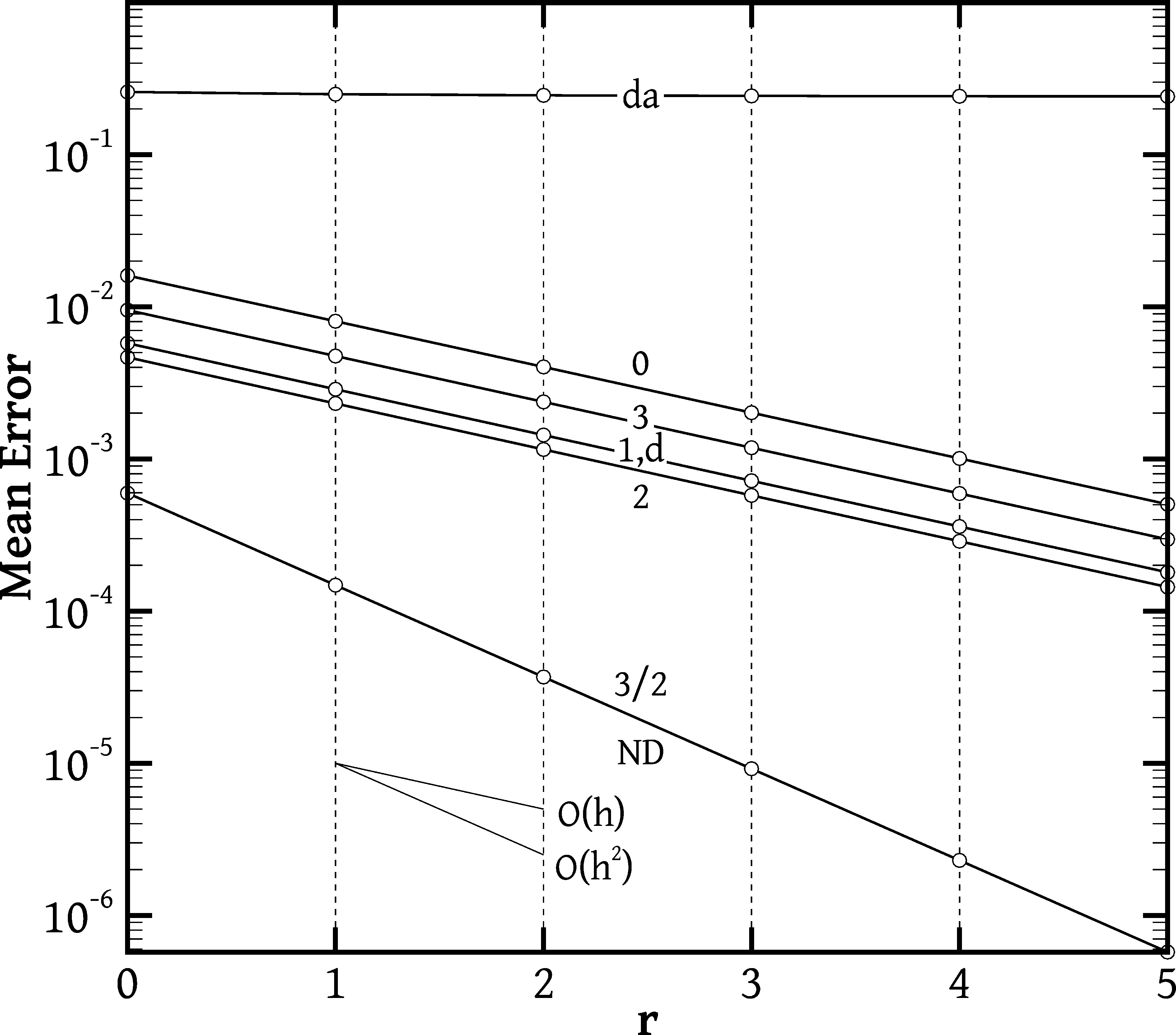}}
\\
 \subfigure[$\{\Delta x_f^0\} = \{-0.10, 0.05, 0.15\}$] {\label{sfig: Dx -0.10 0.05 0.15} 
    \includegraphics[scale=0.97]{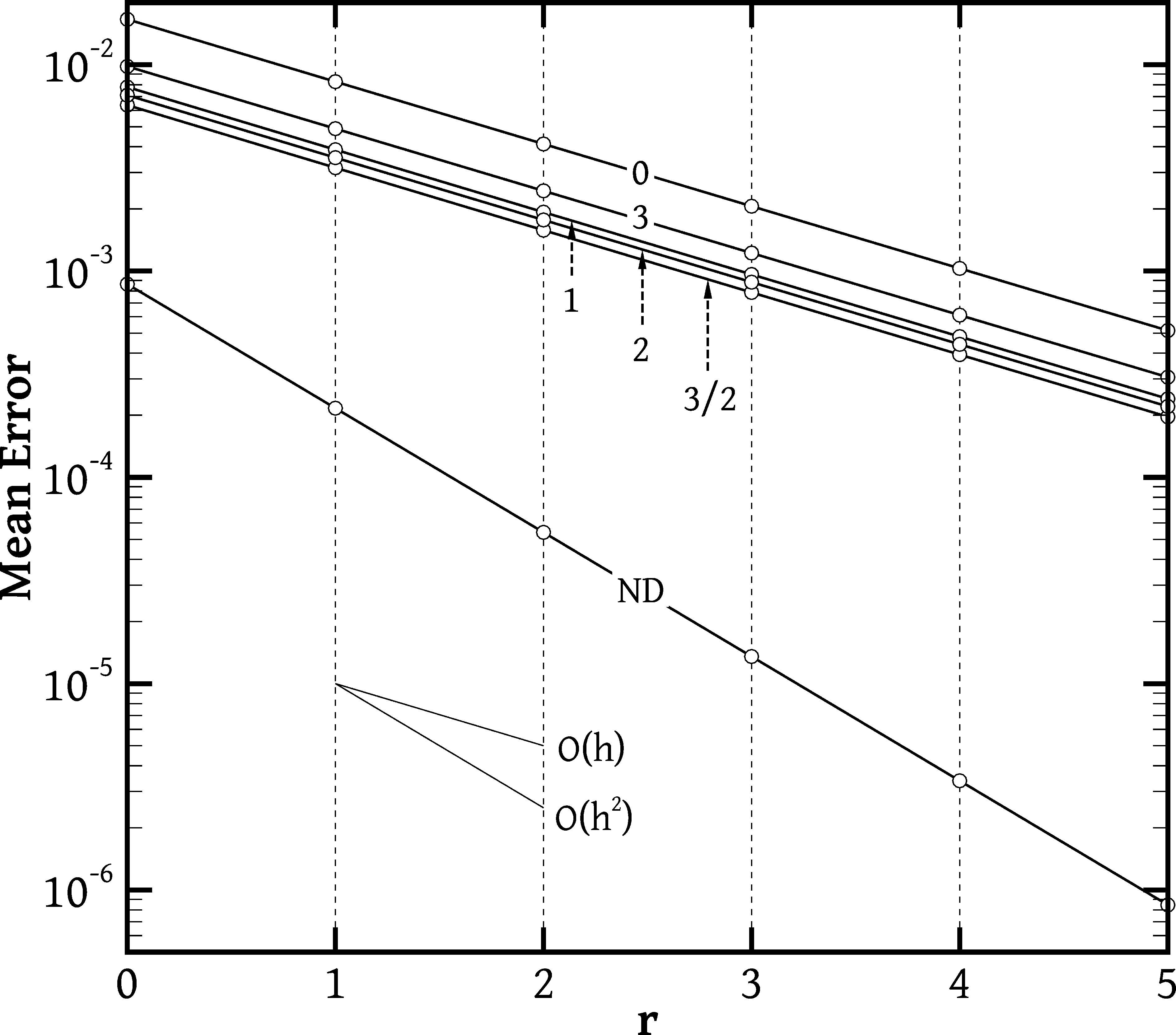}}
 \subfigure[$\{\Delta x_f^0\} = \{-0.20, -0.10, 0.05, 0.15\}$] {\label{sfig: Dx -0.20 -0.10 0.05 
0.15} 
    \includegraphics[scale=0.97]{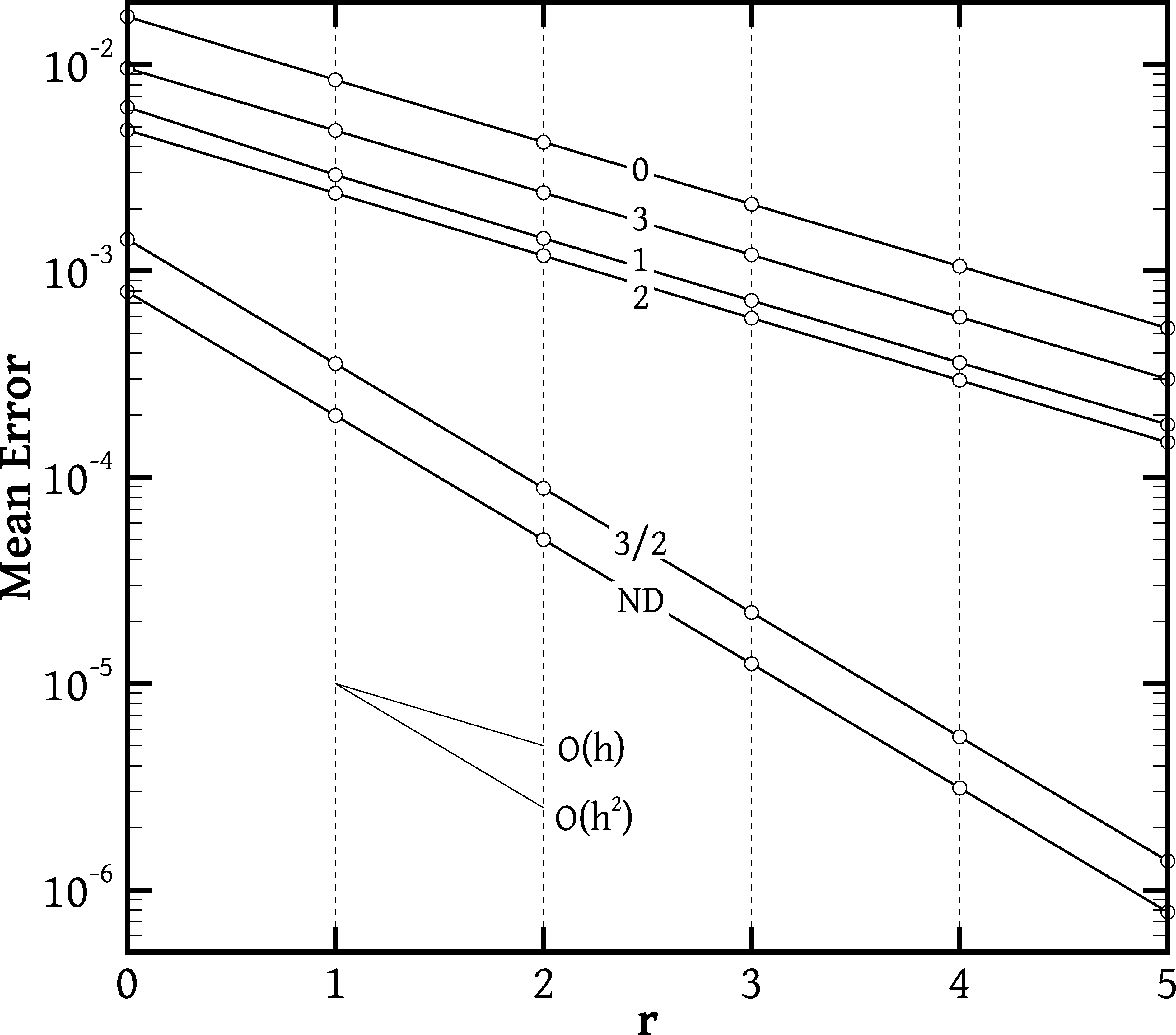}}
 \caption{Mean errors of the calculation of the derivative of $\phi(x) = \tanh(x)$ with various 
methods -- see the text for details.}
 \label{fig: 1d ls results}
\end{figure}

In order to test the methods thoroughly, several sets of initial displacements $\{\Delta 
x_f^0\}_{f=1}^F$ were used. In Fig.\  \ref{fig: 1d ls results}, each diagram corresponds to a 
different such set and the mean error is plotted as a function of the number of displacement 
halvings, $r$. The slope of each curve reveals the order of accuracy of the corresponding method.
The methods tested are: (a) the DT method, denoted ``d'', (b) the LS methods with weight exponents 
$q = 0, 1, 1.5, 2$ and $3$, indicated on each curve, (c) the non-diagonal weights LS method of 
Appendix \ref{sec: appendix ls O(2)} denoted ``ND'', and (d) a simpler variant of the DT method 
which is sometimes used \cite{Sozer_2014}, denoted ``da'', where the values at face centres are 
calculated by arithmetic averaging, $\phi(\vf{c}'_f) = (\phi(\vf{P})+\phi(\vf{N}))/2$, instead of 
linear interpolation \eqref{eq: linear interpolation}. The DT methods are only applicable in the 
case plotted in Fig.\ \ref{sfig: Dx -0.10 0.05} because exactly two neighbours are required, one on 
each side of $x_i$. To apply the DT methods $x_i$ is regarded as the centroid of a cell of size 
equal to the minimum distance between $x_i$ and any of its neighbours.

A result not shown is that when the displacements stencil is symmetric (e.g.\ $\{\Delta x_f^0\} = 
\{-0.1, 0.1\}$) all methods produce identical, second-order accurate results. This is consistent 
with the two-dimensional case where on symmetric grids like that of Fig.\ \ref{fig: parallelogram 
grid} all methods reduce to the formula \eqref{eq: grad parallelogram}. When the stencil is not 
symmetric, Fig.\ \ref{fig: 1d ls results} shows that the LS methods with a diagonal weights matrix 
become first-order accurate, except for the $q = 3/2$ method which retains second-order accuracy 
when there are equal numbers of neighbours on either side (Figs.\ \ref{sfig: Dx -0.10 0.05} and 
\ref{sfig: Dx -0.20 -0.10 0.05 0.15}). The unweighted method ($q = 0$) is always the least accurate; 
the optimum accuracy is achieved with $1 \leq q \leq 2$, while a further increase of $q$ is 
unprofitable ($q = 3$). The method of Appendix \ref{sec: appendix ls O(2)} is always second-order 
accurate. Concerning the DT methods (Fig.\ \ref{sfig: Dx -0.10 0.05}), the method of Sec.\ \ref{sec: 
gauss}, indicated by ``d'' in the figure, gives identical results with the $q = 1$ least squares 
method. On the other hand, the simplified method, indicated by ``da'', is zeroth-order accurate, 
which agrees with the findings reported in \cite{Sozer_2014}.

\subsection{Uniform Cartesian grids}
\label{ssec: results cartesian}

Next, the two-dimensional methods are used to calculate the gradient of the function $\phi(x,y) = 
\tanh(x) \!\cdot\! \tanh(y)$ on the unit square $(x,y) \in [0,1] \times [0,1]$ using uniform 
Cartesian grids of different fineness. The exact gradient is $\phi_{\!.x} = (1 - (\tanh x)^2) \tanh 
y$ and $\phi_{\!.y} = (1 - (\tanh y)^2) \tanh x$. All grid cells are geometrically identical squares 
of side $h = 0.25/2^r$ where $r$ is the level of refinement; however, boundary cells are 
topologically different from interior cells because they have one or more boundary faces where the 
function value at the face centre has to be used (Fig.\ \ref{fig: boundary grid}). It so happens 
that the function $\tanh$ has zero second derivative at the boundaries $x = 0$ and $y = 0$, which 
may artificially increase the order of accuracy of the methods there. However, the general behaviour 
of the methods at boundary cells can be observed at the $x = 1$ and $y = 1$ boundaries where no such 
special behaviour of the $\tanh$ function applies.

\begin{figure}[tb]
 \centering
 \subfigure[mean error, $\tau_{\mathrm{mean}}$]
   {\label{sfig: e_mean cartesian} \includegraphics[scale=1.0]{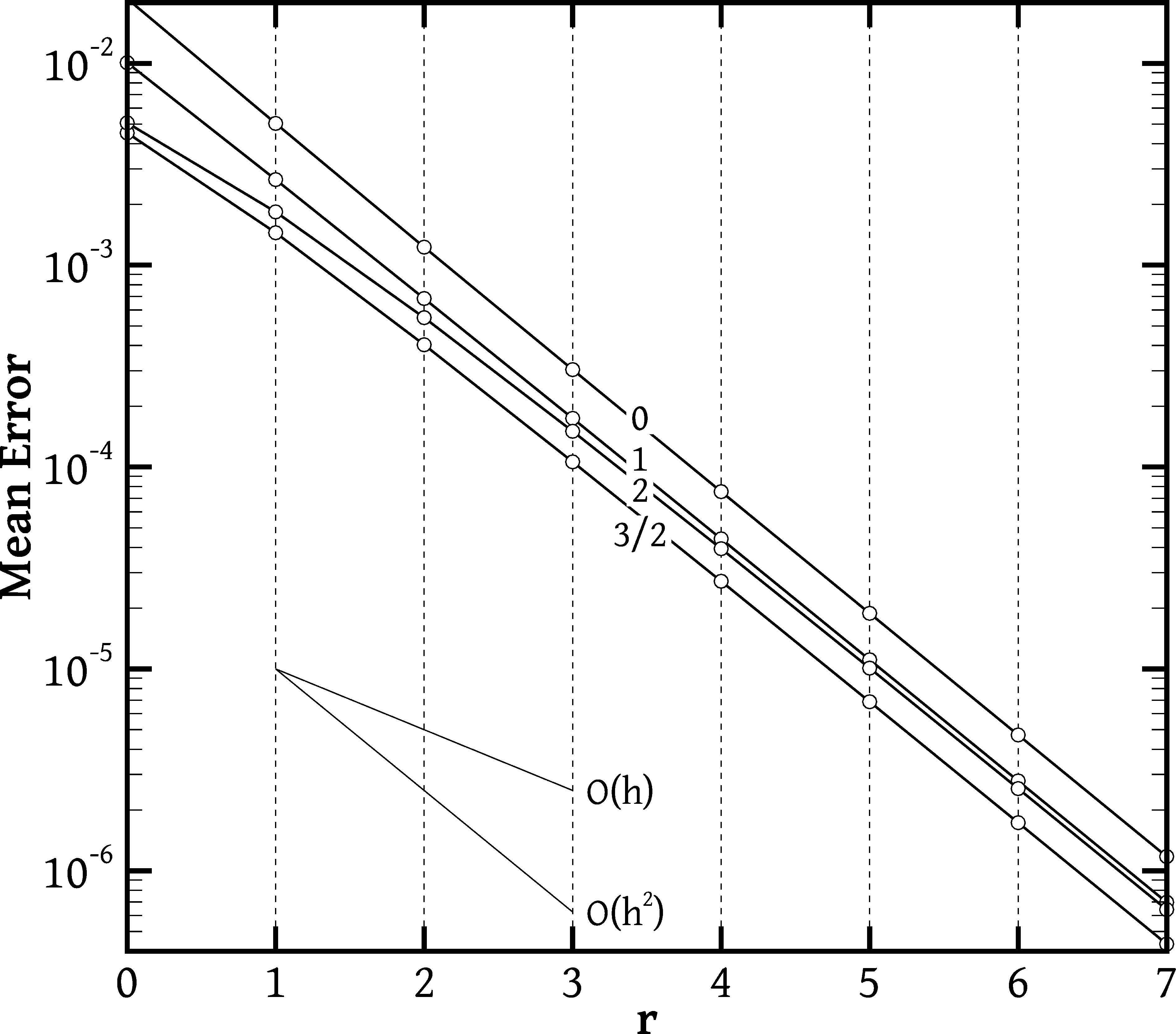}}
 \subfigure[maximum error, $\tau_{\mathrm{max}}$]
   {\label{sfig: e_max cartesian}  \includegraphics[scale=1.0]{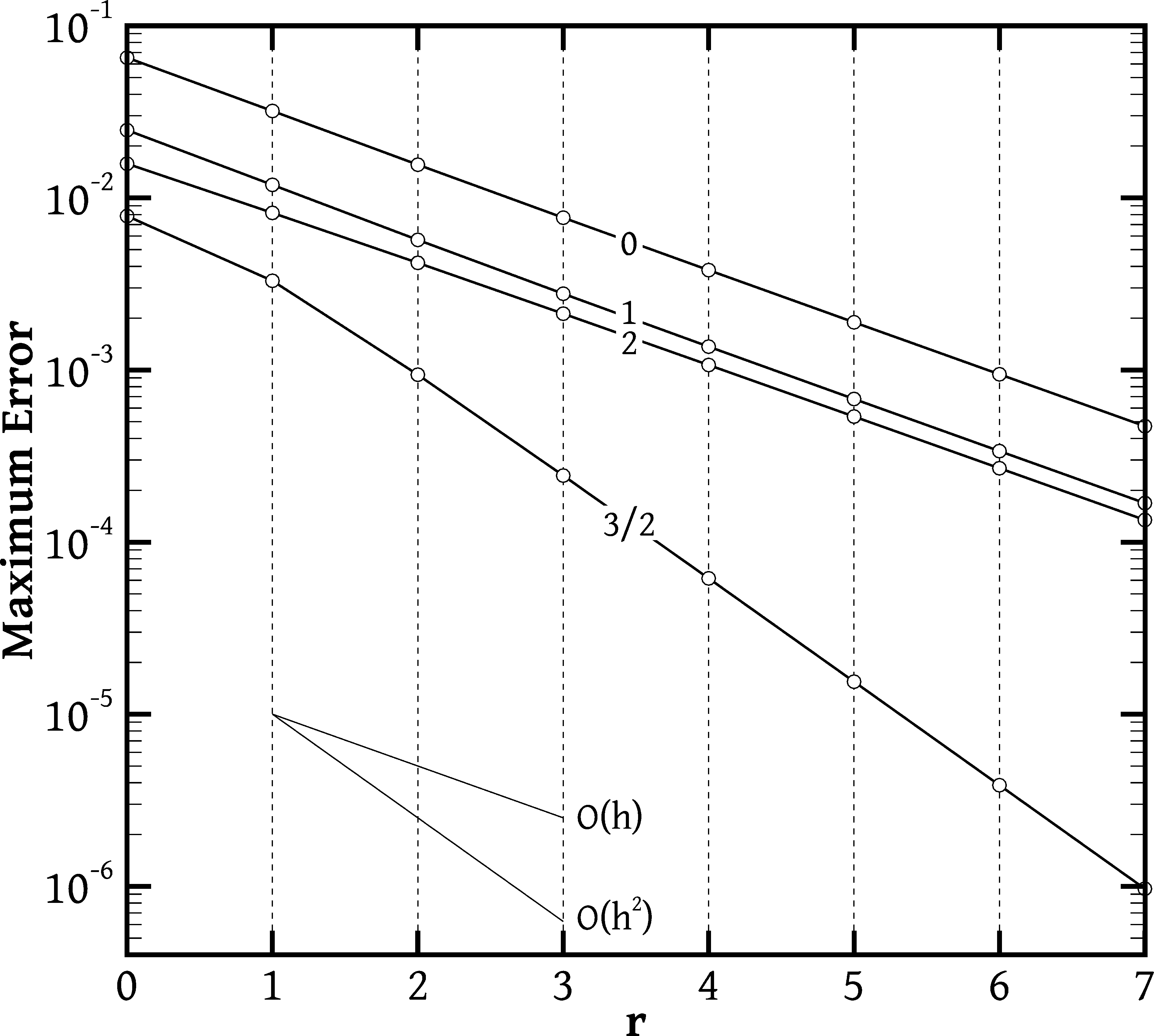}}
 \caption{The mean \subref{sfig: e_mean cartesian} and maximum \subref{sfig: e_max cartesian} 
errors (defined by Eqs.\ \eqref{eq: error mean} and \eqref{eq: error max}, respectively) of the 
various methods for calculating the gradient of the function $\phi = \tanh(x) \tanh(y)$ on uniform 
Cartesian grids. The abscissa $r$ designates the grid; grid $r$ has a uniform spacing of $h = 
0.25/2^r$. The exponent $q$ used for the LS weights is shown on each curve. The DT method produces 
identical results with the $q$ = 1 LS method.}
 \label{fig: errors cartesian}
\end{figure}

On each grid the gradient of $\phi$ is calculated at all cell centres using the DT method 
$\nabla^{\mathrm{d}0}$ (Eq.\ \eqref{eq: gauss gradient 0}), and the LS methods 
$\nabla^{\mathrm{ls}}$ (Eq.\ \eqref{eq: WLS final}) with weight exponents $q = 0, 1, 1.5$ 
and $2$. Since there is no skewness ($\vf{c}' = \vf{c}$ in Eq.\ \eqref{eq: phi cor}), the 
application of corrector steps is meaningless. The methods are evaluated by comparing the mean and 
maximum truncation errors, defined as
\begin{align}
\label{eq: error mean}
 \tau_{\mathrm{mean}} \;&\equiv\; \frac{1}{M_r} \sum_{P=1}^{M_r} \| 
\nabla^{\mathrm{a}}\phi(\vf{P}) - \nabla \phi(\vf{P}) \|
\\[0.2cm]
\label{eq: error max}
 \tau_{\mathrm{max}}  \;&\equiv\; \max_{P=1}^{M_r} \| \nabla^{\mathrm{a}}\phi(\vf{P}) - \nabla 
\phi(\vf{P}) \|
\end{align}
where $\|\cdot\|$ denotes the $L_2$ norm of a vector in a single cell, $M_r$ is the number of cells 
of grid $r$, and $\nabla^{\mathrm{a}}$ is any of the gradient schemes considered. These errors are 
plotted in Figs.\ \ref{sfig: e_mean cartesian} and \ref{sfig: e_max cartesian}, respectively.

The theory predicts that all methods should be second-order accurate at interior cells, i.e.\ 
$\vf{\tau} = O(h^2)$, because they all reduce to formula \eqref{eq: grad parallelogram} there. At 
boundary cells (Fig.\ \ref{fig: boundary grid}) the favourable conditions that are responsible for 
second-order accuracy are lost, and the methods should revert to first-order accuracy except for 
the LS method with $q = 3/2$ of which second-order accuracy is still expected. Figure \ref{sfig: 
e_max cartesian} confirms that the maximum error, which occurs at some boundary cell, is 
$\tau_{\mathrm{max}} = O(h)$ for all methods, except the $q = 3/2$ method for which 
$\tau_{\mathrm{max}} = O(h^2)$. Of the other methods, the unweighted LS method ($q = 0$) is the 
least accurate and the $q = 2$ method is slightly better than the $q = 1$ method. The DT method 
gives identical results with the $q = 1$ LS method, as they both revert to Eq.\ \eqref{eq: gradient 
at boundary} at boundary cells.

Figure \ref{sfig: e_mean cartesian} shows the mean errors. Since all methods give the same results 
at interior cells any differences are due to the different errors at the boundary cells. The slope 
of the curves is the same, and corresponds to $\tau_{\mathrm{mean}} = O(h^2)$. This does not 
contradict with the $O(h)$ errors at the boundary cells: the number of such cells along each side of 
the domain equals $1/h$, so that in total there are $O(1/h)$ boundary cells, each contributing an 
$O(h)$ error. Their total contribution to $\tau_{\mathrm{mean}}$ in Eq.\ \eqref{eq: error mean} 
is therefore $O(1/h) \!\cdot\! O(h) / M_r = O(h^2)$ since $M_r = 1/h^2$. The interior cell 
contribution is $O(h^2)$ also, since $\tau = O(h^2)$ at each individual interior cell.

\subsection{Smooth curvilinear grids}
\label{ssec: results curvilinear}

\begin{figure}[tb]
 \centering
 \subfigure[$r=0$] {\label{sfig: grid elliptic 0} 
\includegraphics[scale=0.75]{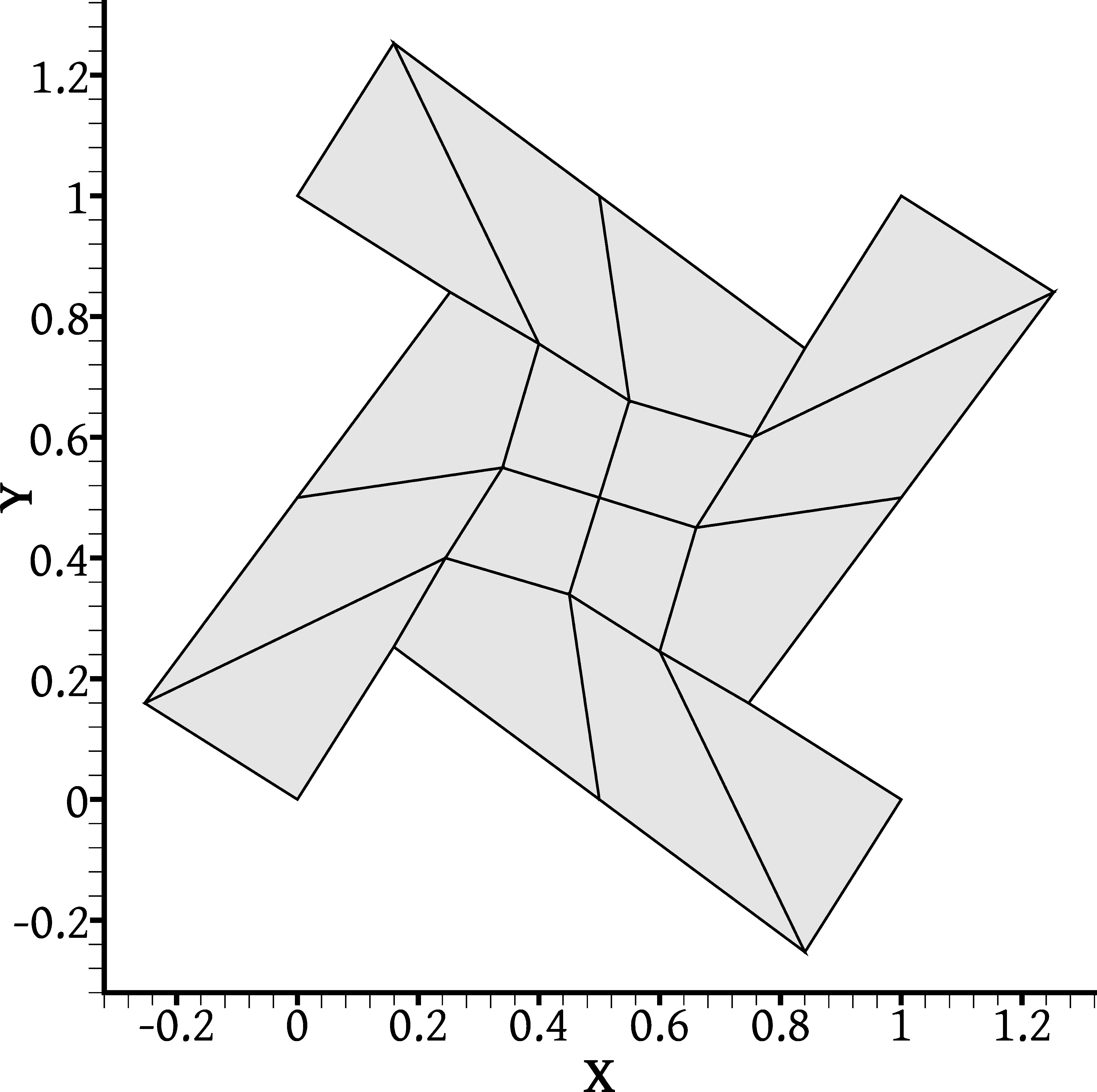}}
 \subfigure[$r=1$] {\label{sfig: grid elliptic 1} 
\includegraphics[scale=0.75]{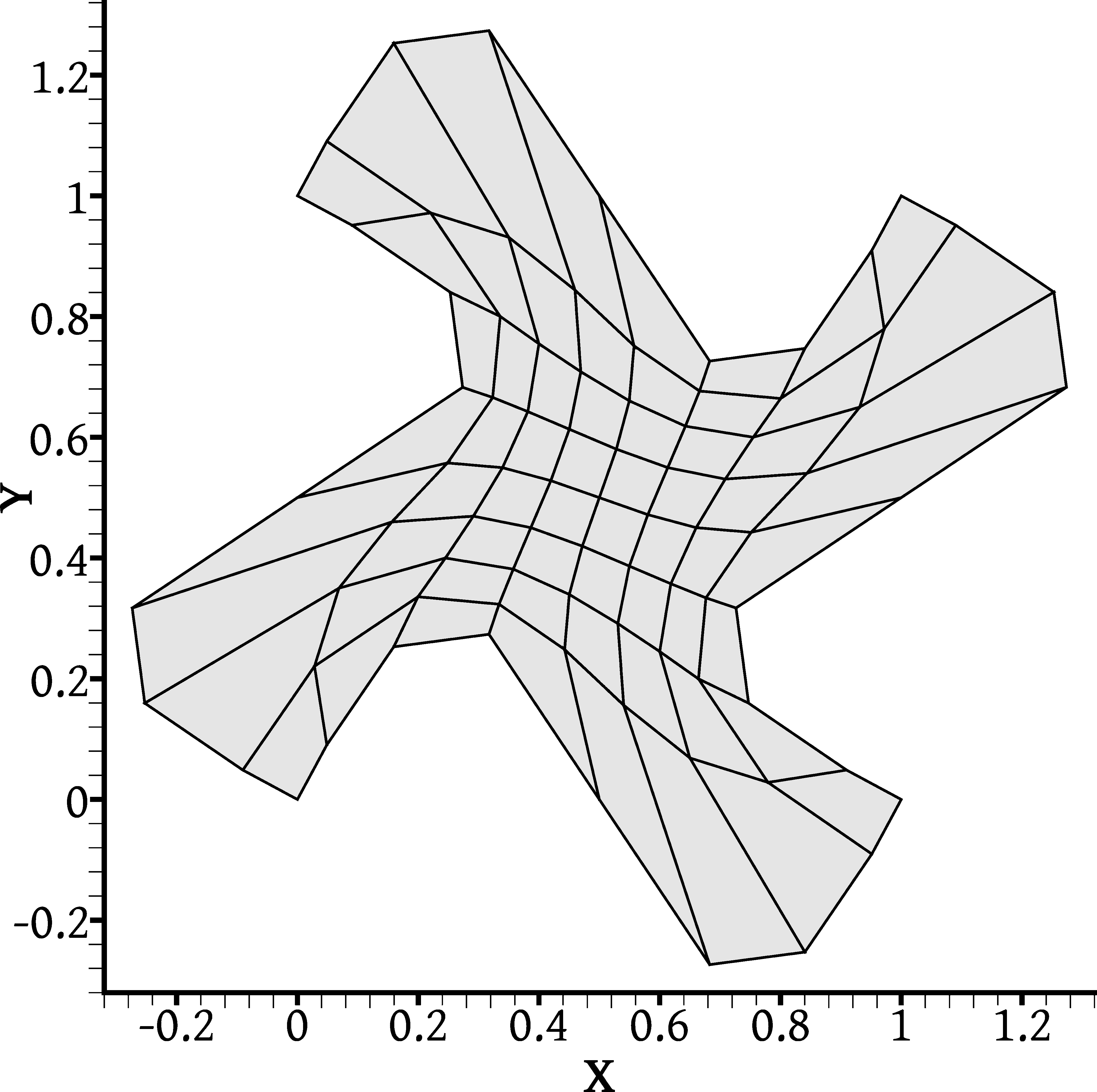}}
 \subfigure[$r=2$] {\label{sfig: grid elliptic 2} 
\includegraphics[scale=0.75]{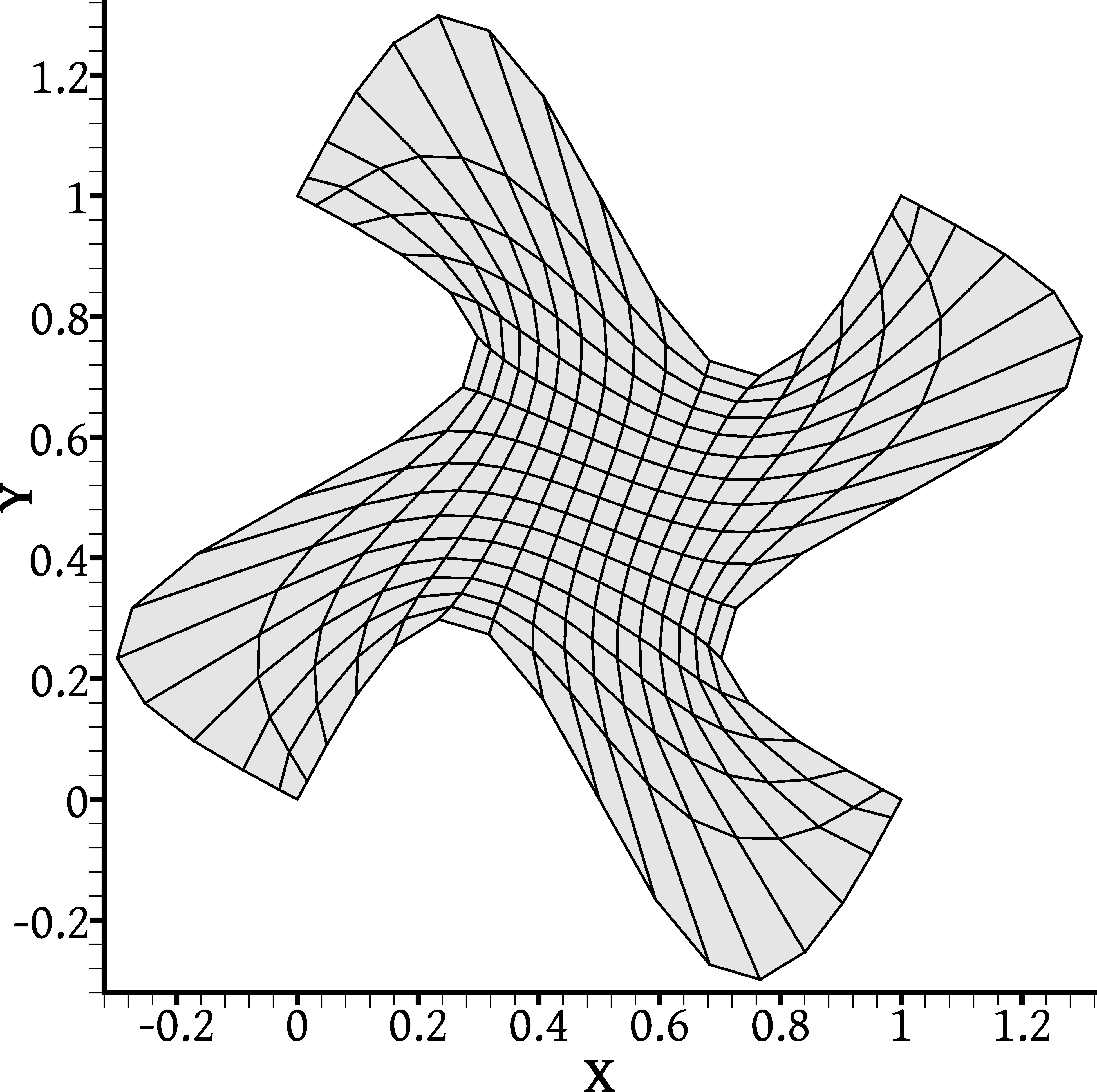}}
 \subfigure[$r=3$] {\label{sfig: grid elliptic 3} 
\includegraphics[scale=0.75]{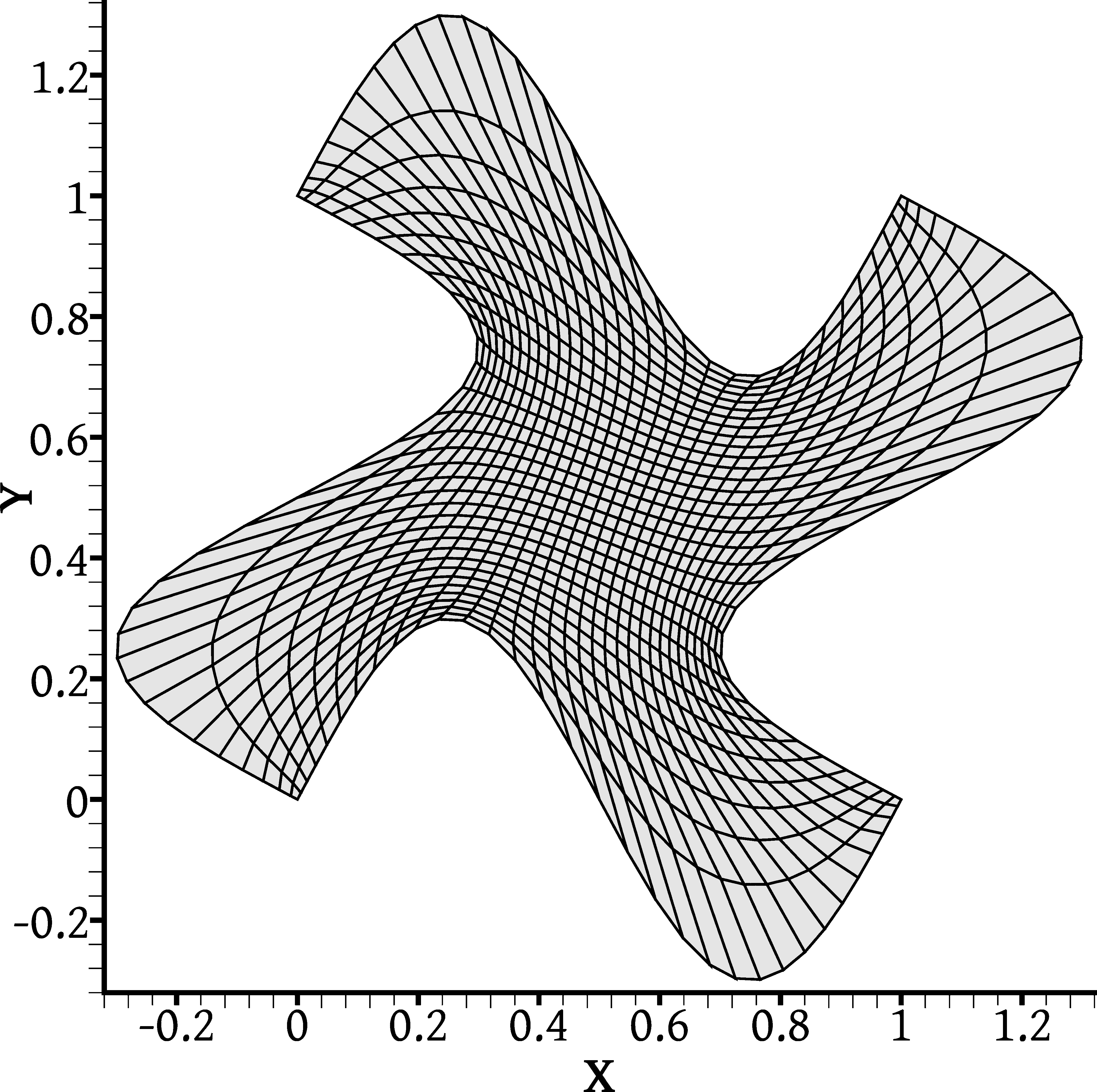}}
 \caption{The first four of the series of grids constructed on a domain with sinusoidal boundaries 
via elliptic grid generation. See the text for more details.}
 \label{fig: grids elliptic}
\end{figure}

Next, we try the methods on smooth curvilinear grids. The same function $\phi = \tanh(x) \!\cdot\! 
\tanh(y)$ is differentiated, but the domain boundaries now have the shapes of two horizontal and 
two vertical sinusoidal waves (Fig.\ \ref{fig: grids elliptic}), beginning and ending at the points 
$(0,0)$, $(1,0)$, $(1,1)$ and $(0,1)$. The grid is generated using a very basic elliptic grid 
generation method \cite{Thompson_1985}. In particular, smoothly varying functions $\xi(x,y)$ and 
$\eta(x,y)$ are assumed in the domain, and the grid consists of lines of constant $\xi$ and of 
constant $\eta$. The left, right, bottom and top boundaries correspond to $\xi = 0$, $\xi = 1$, 
$\eta = 0$ and $\eta = 1$, respectively. In the interior of the domain $\xi$ and $\eta$ are assumed 
to vary according to the following Laplace equations:
\begin{align*}
 \xi_{.xx} \;+\; \xi_{.yy} \;&=\; 0
\\
 \eta_{.xx} \;+\; \eta_{.yy} \;&=\; 0
\end{align*}

These equations guarantee that $\xi$ and $\eta$ vary smoothly in the domain, but their solution  
$\xi = \xi(x,y)$, $\eta = \eta(x,y)$ is not much help in constructing the grid. Instead, we need the 
inverse functions $x=x(\xi,\eta)$, $y=y(\xi,\eta)$ which explicitly set the locations of all grid 
nodes; node $(i,j)$ is located at $(x_{i,j}, y_{i,j}) \equiv (x(\xi = i\,\Delta\xi, 
\eta=j\,\Delta\eta), y(\xi = i\,\Delta\xi, \eta=j\,\Delta\eta))$. With $\xi,\eta \in [0,1]$, the 
constant spacings $\Delta \xi$ and $\Delta \eta$ are adjusted according to the desired grid 
fineness. Therefore, using the chain rule of partial differentiation it can be shown that the above 
equations can be expressed in inverse form as
\begin{align*}
  g_{22} \, x_{\!.\xi\xi} \;-\; 2 g_{12} \, x_{\!.\xi\eta} \;+\; g_{11} \, x_{\!.\eta\eta} \;&=\; 0
\\
  g_{22} \, y_{\!.\xi\xi} \;-\; 2 g_{12} \, y_{\!.\xi\eta} \;+\; g_{11} \, y_{\!.\eta\eta} \;&=\; 0
\end{align*}
where
\begin{align*}
 g_{11} \;&=\; x^2_{\!.\xi} +y^2_{\!.\xi}
\\
 g_{22} \;&=\; x^2_{\!.\eta} +y^2_{\!.\eta}
\\
 g_{12} \;&=\; x_{\!.\xi} x_{\!.\eta} + y_{\!.\xi} y_{\!.\eta}
\end{align*}
In order to cluster the points near the boundaries in the physical domain, we accompany the above 
equations with the following boundary conditions: at the bottom boundary we set $x = 0.5 + 0.5 
\sin(\pi (\xi-0.5))$ and $y = \sin(2\pi x)$, and at the left boundary we set $y = 0.5 + 0.5 
\sin(\pi (\eta-0.5))$ and $x = -\sin(2\pi y)$. At the top and right boundaries we set the same 
conditions, respectively, adding 1 to $x$ at the right boundary and 1 to $y$ at the top boundary. 
Better results can be obtained by using a more elaborate method such as described in 
\cite{Dimakopoulos_2003}, but this suffices for the present purposes.

Now $(x,y)$ have become the dependent variables while $(\xi,\eta)$ are the independent variables, 
which acquire values in the unit square $[0,1]\times[0,1]$. The grid equations were solved 
numerically with a finite difference method on a $513 \times 513$ point uniform Cartesian grid. Note 
that the dependent variables $(x,y)$ are stored at the grid nodes, i.e.\ at the intersection points 
of the grid lines, instead of at the cell centres. The derivatives are approximated by second-order 
accurate central differences; for example, at point $(i,j)$, $x_{\!.\xi} \approx (x_{i+1,j} - 
x_{i-1,j})/2h$, $x_{\!.\eta} \approx (x_{i,j+1} - x_{i,j-1})/2h$, $x_{\!.\xi\xi} \approx (x_{i+1,j} 
- 2x_{i,j} + x_{i-1,j})/h^2$ etc.\ where $h = 1/512$ is the grid spacing. The resulting system of 
nonlinear algebraic equations was solved using a Gauss-Seidel iterative method where in the 
equations of the $(i,j)$ node all terms are treated as known from their current values except for 
$x_{i,j}$ and $y_{i,j}$ which are solved for. The convergence of the method was accelerated using a 
minimal polynomial extrapolation technique \cite{Sidi_2012}, and iterations were carried out until 
machine precision was reached.

\begin{figure}[tb]
 \centering
 \subfigure[mean error, $\tau_{\mathrm{mean}}$]
   {\label{sfig: e_mean elliptic} \includegraphics[scale=1.0]{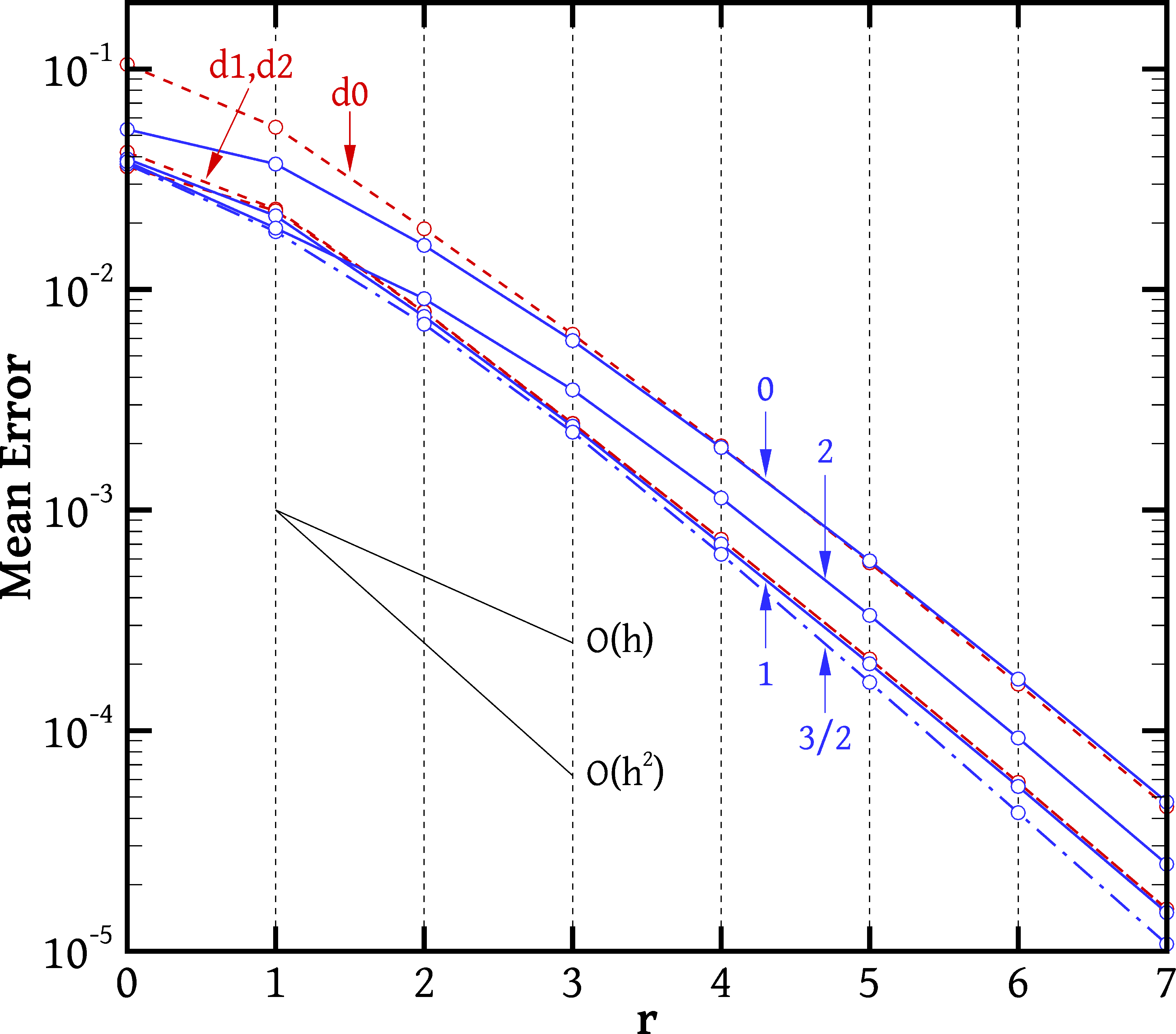}}
 \subfigure[maximum error, $\tau_{\mathrm{max}}$] 
   {\label{sfig: e_max elliptic}  \includegraphics[scale=1.0]{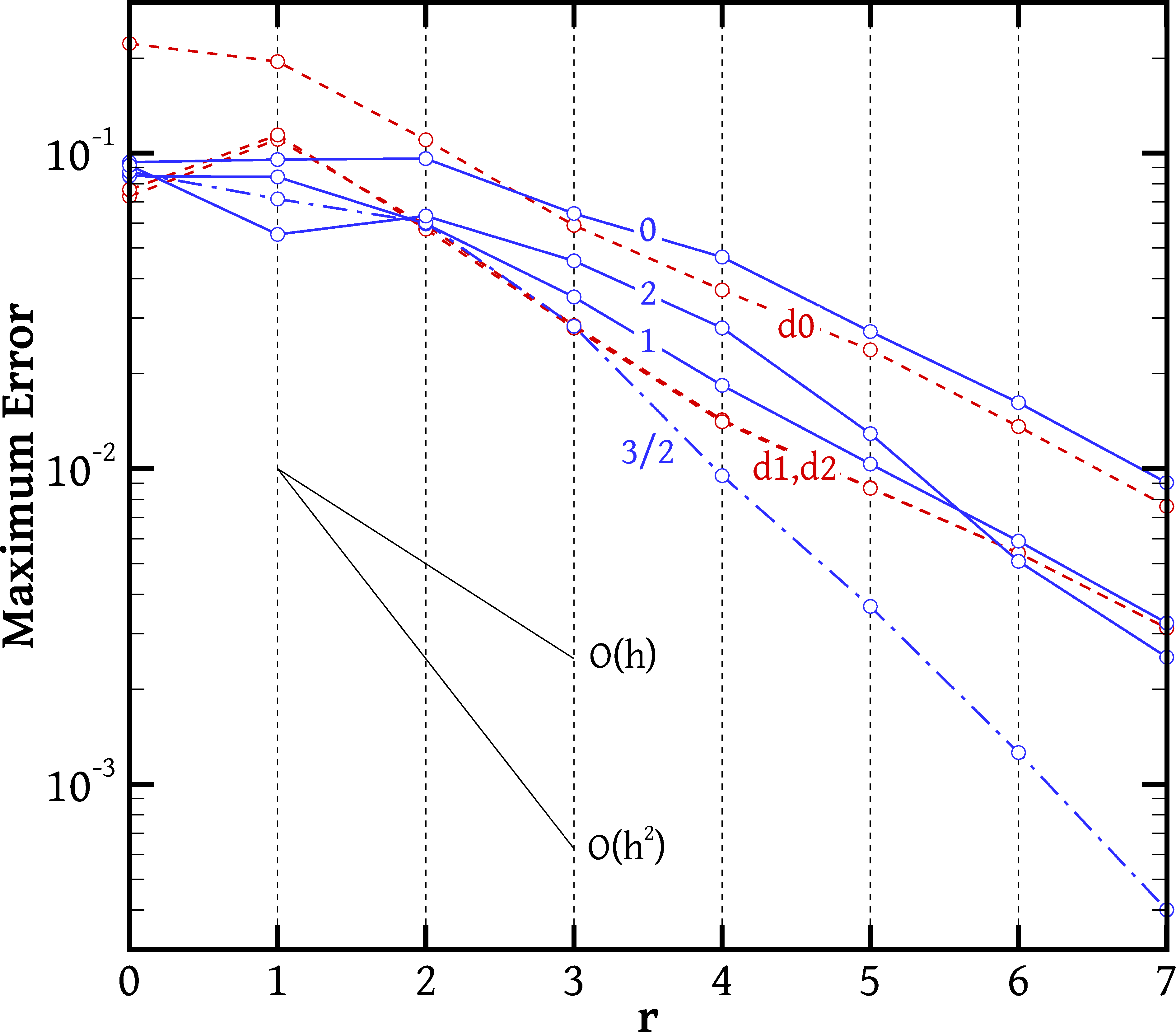}}
 \caption{The mean \subref{sfig: e_mean elliptic} and maximum \subref{sfig: e_max elliptic} errors 
(defined by Eqs.\ \eqref{eq: error mean} and \eqref{eq: error max}, respectively) of the gradient 
schemes applied on the function $\phi = \tanh(x) \tanh(y)$ on smooth curvilinear grids (Fig.\ 
\ref{fig: grids elliptic}). The abscissa $r$ designates the grid; $r = 0$ is the coarsest grid 
(Fig.\ \ref{sfig: grid elliptic 0}), and grid $r$ comes from subdividing every cell of grid $r-1$ 
into 4 child cells \textit{in the computational space} (see text). The blue solid lines correspond 
to the LS methods with $q$ = 0, 1 and 2 as indicated on each curve; the blue dash-dot line 
corresponds to the LS method with $q = 3/2$; and the red dashed lines correspond to the DT methods 
$\mathrm{d}c$ where $c$ is the number of corrector steps.}
 \label{fig: errors elliptic}
\end{figure}

After obtaining $x(\xi,\eta)$ and $y(\xi,\eta)$, a series of successively refined grids in the 
physical domain were constructed by drawing lines of constant $\xi$ and lines of constant $\eta$ at 
intervals $\Delta \xi = \Delta \eta = 0.25 / 2^r$, for $r = 0, 1, \ldots 7$. The first four grids 
are shown in Fig.\ \ref{fig: grids elliptic}. The gradient calculation methods were then applied on 
each of these grids. The errors of each method are depicted in Fig.\ \ref{fig: errors elliptic}.

These grids exhibit all kinds of grid irregularity, but the unevenness and the skewness diminish 
with grid refinement, as explained in  Section \ref{sec: preliminary}. In particular, Table 
\ref{table: exactness linear} shows that the measures of both these grid qualities have magnitudes 
of $O(h)$. These measures were defined in Sec.\ \ref{sec: preliminary}: $\|\vf{c}'_f-\vf{m}_f\| / 
\|\vf{N}_f-\vf{P}\|$ for unevenness and $\|\vf{c}_f - \vf{c}'_f\| / \|\vf{N}_f-\vf{P}\|$ for 
skewness (see Fig\ \ref{fig: unstructured grid} for definitions). The values listed in Table 
\ref{table: exactness linear} are the average values among all faces of each grid, excluding 
boundary faces. Therefore, it is expected that eventually the methods will behave as in the 
Cartesian case. Indeed, Fig.\ \ref{sfig: e_max elliptic} shows that as the grid is refined 
$\tau_{\mathrm{max}}$ tends to decrease at a first-order rate for all methods, except for the least 
squares method with $q = 3/2$, for which it decreases at a second-order rate. Accordingly, Fig.\ 
\ref{sfig: e_mean elliptic} shows that $\tau_{\mathrm{mean}}$ tends to decrease at a second-order 
rate for all methods, with the second-order rate attained earlier by the $q = 3/2$ method. 
Therefore, like on the Cartesian grids, the methods are second-order accurate at interior cells but 
revert to first-order accuracy at boundary cells, except for the $q = 3/2$ method.

Of the LS methods the least accurate is the unweighted method, followed by the $q = 2$ method. The 
undisputed champion is the $q = 3/2$ method because, as mentioned, it retains its second-order 
accuracy even at boundary cells. Concerning the DT methods, the method with no corrections (Eq.\ 
\ref{eq: gauss gradient 0}) performs similarly to the unweighted LS method. Application of a 
corrector step (Eq.\ \eqref{eq: gauss gradient 1}) now does make a difference, since skewness is 
present at any finite grid density, bringing the accuracy of the method on a par with the best 
weighted LS methods, except of course the $q = 3/2$ method. We also tried a second corrector step 
but it did not bring any noticeable improvement.

Next, on the same grids we also applied discrete gradient operators to calculate the gradient of a 
linear and of a quadratic function. The schemes tested are the DT gradient operators 
$\nabla^{\mathrm{d}0}$, $\nabla^{\mathrm{d}1}$ and $\nabla^{\mathrm{d}\infty}$ (the latter 
approximated with 100 correction steps), and the LS gradient operators with $q$ = 0, 1 and 3/2, 
denoted as $\nabla^{\mathrm{ls}0}$, $\nabla^{\mathrm{ls}1}$ and $\nabla^{\mathrm{ls}3/2}$, 
respectively, in Tables \ref{table: exactness linear} and \ref{table: exactness quadratic}. 
``Exactness'' is sometimes used as an aid to either determine the order of accuracy of a method or 
to design a method to achieve a desired order of accuracy (e.g.\ Appendix \ref{sec: appendix ls 
O(2)}): a first-order accurate gradient scheme would normally be exact for linear functions, and a 
second-order accurate gradient scheme would normally be exact for quadratic functions. However, the 
results listed in Tables \ref{table: exactness linear} and \ref{table: exactness quadratic} show 
that in the present case the DT gradient is not exact even for linear functions while the LS 
gradient is exact for linear functions but not for quadratic functions, despite both methods being 
second-order accurate.

In particular, Table \ref{table: exactness linear} shows that the DT gradient without corrector 
steps ($\nabla^{\mathrm{d}0}$) is not exact for the linear function (the errors are not zero) but 
converges to the exact gradient at a rate that approaches second-order as the grid is refined. 
Performing a corrector step ($\nabla^{\mathrm{d}1}$) brings a significant improvement in accuracy, 
with an observed convergence rate order of between 2 and 3, but still the operator is not exact. 
This inexactness is anticipated since the grids are skewed and the DT scheme cannot cope with 
skewness. However, grid refinement causes skewness to diminish and the DT accuracy to improve at a 
second-order rate. In the limit of many corrector steps ($\nabla^{\mathrm{d}\infty}$) the operator 
becomes exact, with the errors at machine precision levels even at the coarsest grid. All of the LS 
schemes are also exact for the linear function.

Concerning the quadratic function (Table \ref{table: exactness quadratic}), none of 
the schemes is exact but they are all second-order accurate, with the $q = 3/2$ LS scheme being the 
most accurate, and the $q = 0$ LS and zero-correction DT schemes being the least accurate. A single 
corrector step in the DT scheme ($\nabla^{\mathrm{d}1}$) brings the maximum attainable improvement, 
since the error levels of $\nabla^{\mathrm{d}1}$ and $\nabla^{\mathrm{d}\infty}$ are nearly 
identical. The second-order convergence rates are due to the improvement of grid quality with 
refinement, as explained in Sections \ref{sec: gauss} and \ref{sec: least squares}.

\begin{table}[thb]
\caption{Mean errors (Eq.\ \eqref{eq: error mean}) of various schemes for calculating the gradient 
$\nabla\phi = (1,2)$ of the linear function $\phi(x,y) = x + 2y + 0.5$ on the series of grids shown 
in Fig.\ \ref{fig: grids elliptic}. Also displayed are the measures of grid skewness and unevenness 
(defined in Sec.\ \ref{sec: preliminary}), averaged over all faces of each grid excluding boundary 
faces.}
\label{table: exactness linear}
\begin{center}
\begin{small}   
\renewcommand\arraystretch{1.25}   
{
\newcommand{\G}[2]{\;$\nabla^{\mathrm{#1}#2}$\;}
\newcommand{\e}[1]{$\cdot 10^{-#1}$}
\begin{tabular}{ c | c c | c c c c c c }
\toprule
 Grid $r$ & Skew.      & Unev.     &  \G{d}{0}    &  \G{d}{1}   & \G{d}{\infty} &  \G{ls}{0}  &  \G{ls}{1}  & \G{ls}{3/2}
\\ \midrule
 0        & 1.33\e{1}  & 1.10\e{1} &  4.16\e{1}   & 8.00\e{2}   &  2.19\e{15}   & 6.04\e{16}  & 7.39\e{16}  & 9.37\e{16}
\\
 1        & 5.68\e{2}  & 8.00\e{2} &  1.86\e{1}   & 1.90\e{2}   &  4.46\e{15}   & 1.34\e{15}  & 1.45\e{15}  & 1.68\e{15}
\\
 2        & 3.03\e{2}  & 4.83\e{2} &  6.93\e{2}   & 4.80\e{3}   &  9.53\e{15}   & 2.43\e{15}  & 2.41\e{15}  & 2.81\e{15}
\\
 3        & 1.64\e{2}  & 2.70\e{2} &  2.41\e{2}   & 8.91\e{4}   &  1.93\e{14}   & 4.42\e{15}  & 4.52\e{15}  & 4.84\e{15}
\\
 4        & 8.44\e{3}  & 1.43\e{2} &  7.72\e{3}   & 1.16\e{4}   &  3.97\e{14}   & 8.74\e{15}  & 8.77\e{15}  & 9.07\e{15}
\\
 5        & 4.29\e{3}  & 7.41\e{3} &  2.32\e{3}   & 1.43\e{5}   &  7.98\e{14}   & 1.70\e{14}  & 1.71\e{14}  & 1.74\e{14}
\\
 6        & 2.16\e{3}  & 3.78\e{3} &  6.64\e{4}   & 1.94\e{6}   &  1.60\e{13}   & 3.36\e{14}  & 3.38\e{14}  & 3.40\e{14}
\\
 7        & 1.09\e{3}  & 1.91\e{3} &  1.86\e{4}   & 3.25\e{7}   &  3.20\e{13}   & 6.69\e{14}  & 6.70\e{14}  & 6.72\e{14}
\\
\bottomrule
\end{tabular}
} 
\end{small}
\end{center}
\end{table}

\begin{table}[thb]
\caption{Mean errors (Eq.\ \eqref{eq: error mean}) of various schemes for calculating the gradient 
$\nabla\phi = (2x+2y,2x-2y)$ of the quadratic function $\phi(x,y) = x^2 + 2xy - y^2$ on the series 
of grids shown in Fig.\ \ref{fig: grids elliptic}.}
\label{table: exactness quadratic}
\begin{center}
\begin{small}   
\renewcommand\arraystretch{1.25}   
{
\newcommand{\G}[2]{\;$\nabla^{\mathrm{#1}#2}$\;}
\newcommand{\e}[1]{$\cdot 10^{-#1}$}
\begin{tabular}{ c | c c c c c c }
\toprule
 Grid $r$  &  \G{d}{0} &  \G{d}{1} & \G{d}{\infty} &  \G{ls}{0} &  \G{ls}{1} & \G{ls}{3/2}
\\ \midrule
 0         & 4.95\e{1} & 2.28\e{1} &  2.29\e{1}    & 3.54\e{1}  & 1.94\e{1}  & 1.47\e{1} 
\\
 1         & 2.20\e{1} & 1.27\e{1} &  1.29\e{1}    & 1.72\e{1}  & 8.16\e{2}  & 5.81\e{2} 
\\
 2         & 8.38\e{2} & 4.48\e{2} &  4.47\e{2}    & 8.09\e{2}  & 3.13\e{2}  & 2.24\e{2} 
\\
 3         & 2.84\e{2} & 1.46\e{2} &  1.46\e{2}    & 3.20\e{2}  & 1.12\e{2}  & 7.43\e{3} 
\\
 4         & 8.99\e{3} & 4.63\e{3} &  4.64\e{3}    & 1.14\e{2}  & 3.76\e{3}  & 2.05\e{3} 
\\
 5         & 2.66\e{3} & 1.40\e{3} &  1.40\e{3}    & 3.75\e{3}  & 1.19\e{3}  & 5.28\e{4} 
\\
 6         & 7.45\e{4} & 4.00\e{4} &  4.01\e{4}    & 1.15\e{3}  & 3.50\e{4}  & 1.35\e{4} 
\\
 7         & 2.07\e{4} & 1.10\e{4} &  1.10\e{4}    & 3.29\e{4}  & 9.74\e{5}  & 3.41\e{5} 
\\
\bottomrule
\end{tabular}
} 
\end{small}
\end{center}
\end{table}

\subsection{Grids of localised high distortion}
\label{ssec: results refined}

Structured grids that are constructed not by solving partial differential equations, as in Section 
\ref{ssec: results curvilinear}, but by algebraic methods may lack the property that unevenness and 
skewness diminish with grid refinement. This is especially true if the domain boundaries include 
sharp corners at points other than grid line endpoints. For example, the grid of Fig.\ \ref{fig: 
grid algebraic} is structured, consisting of piecewise straight lines. At the line joining the 
sharp corners, the intersecting grid lines change direction abruptly. This causes significant 
skewness which is unaffected by grid refinement.

\begin{figure}[tb]
 \centering
 \includegraphics[scale=0.5]{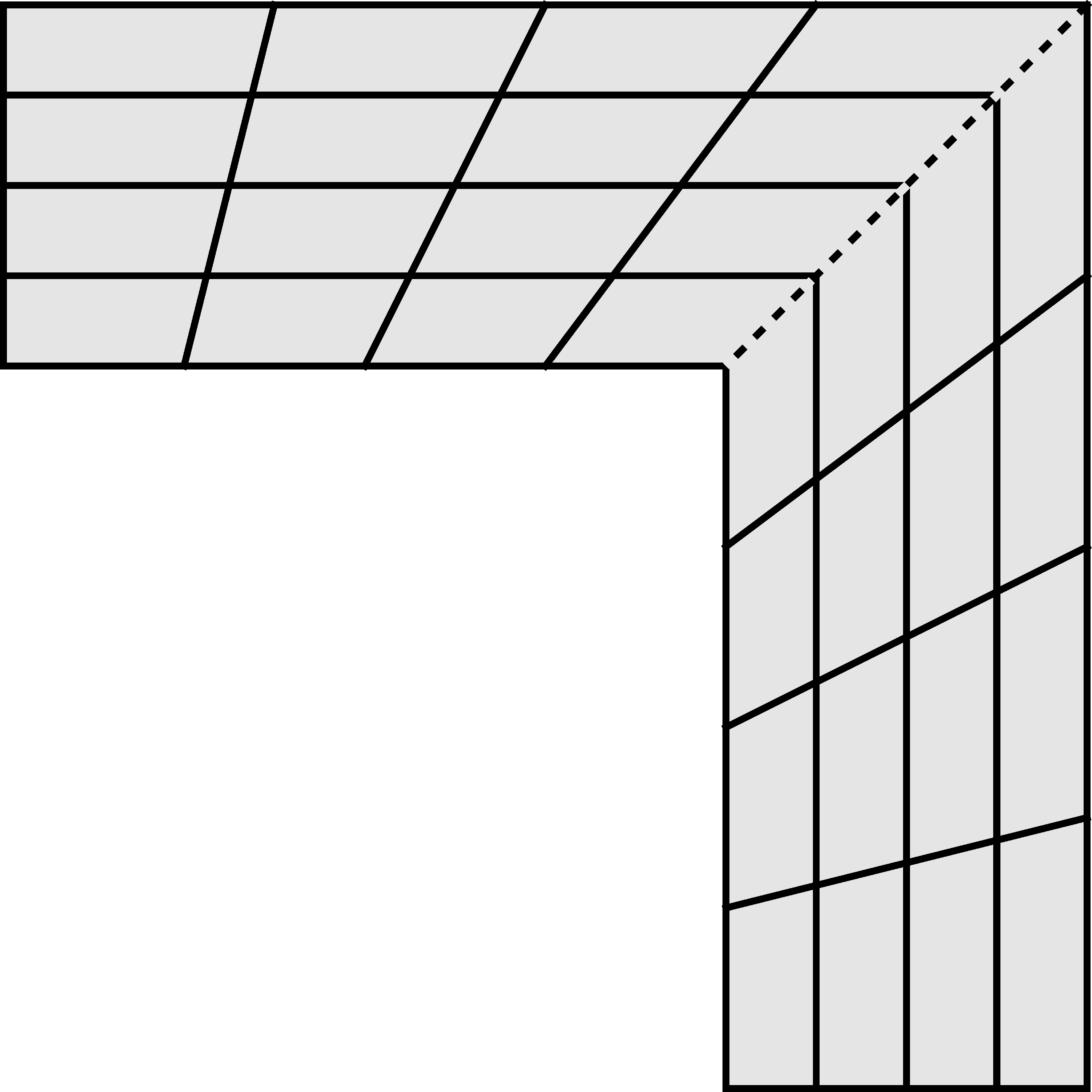}
 \caption{A structured grid where the grid lines belonging to one family change direction abruptly 
at the dashed line joining the pair of sharp corners, where grid skewness and not diminishing with 
grid refinement.}
 \label{fig: grid algebraic}
\end{figure}

\begin{figure}[tb]
 \centering
 \subfigure[$r=0$] {\label{sfig: grid refined 0} 
\includegraphics[scale=0.65]{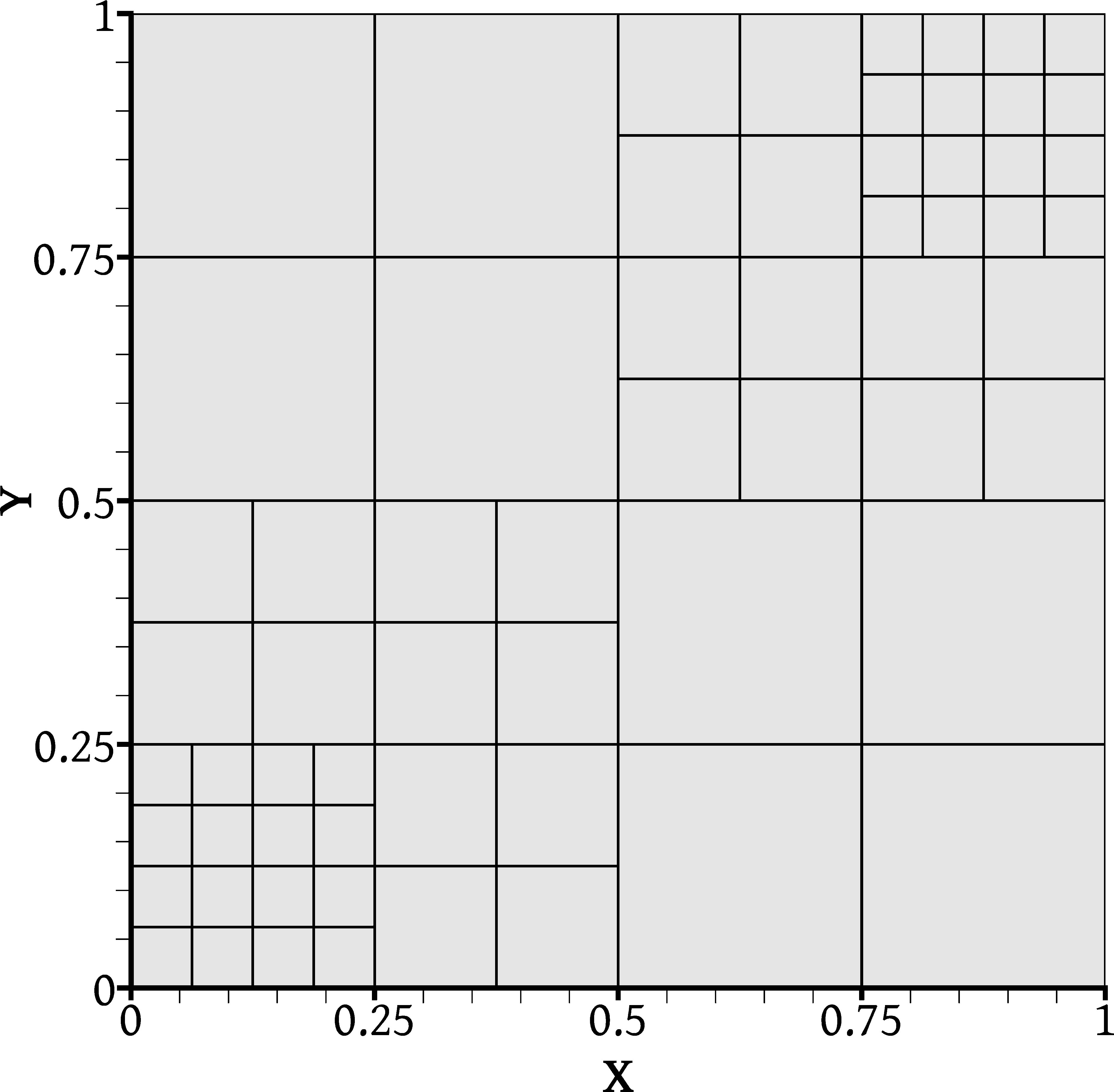}}
 \subfigure[$r=1$] {\label{sfig: grid refined 1} 
\includegraphics[scale=0.65]{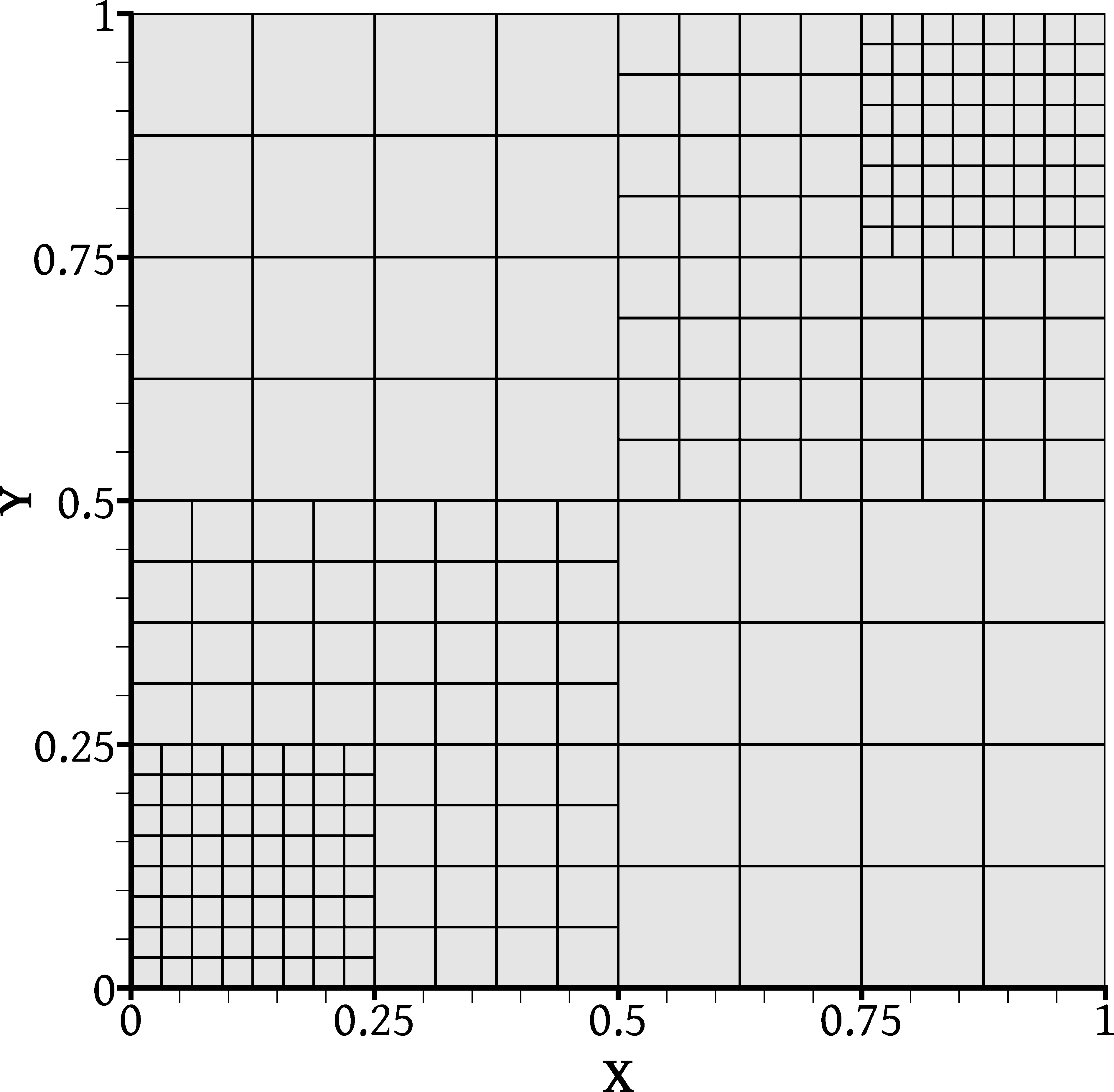}}
 \subfigure[$r=2$] {\label{sfig: grid refined 2} 
\includegraphics[scale=0.65]{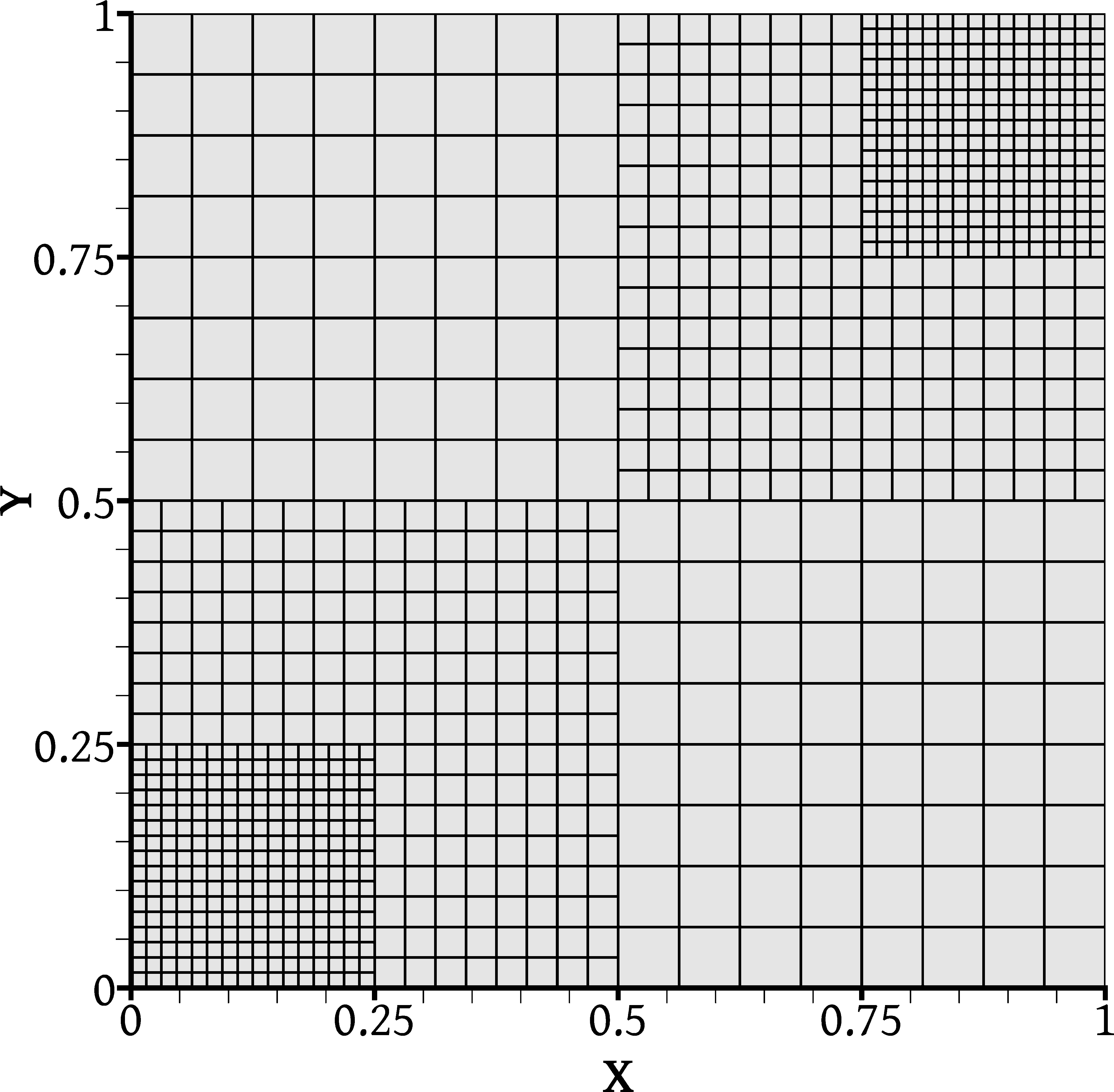}}
 \caption{A series of multi-level, or composite, grids. Each grid $r$ comes from the previous grid 
$r-1$ by evenly subdividing each cell into four child cells.}
 \label{fig: grids refined}
\end{figure}

\begin{figure}[!tb]
 \centering
 \subfigure[coarse cell] {\label{sfig: composite grid coarse cell} 
\includegraphics[scale=0.85]{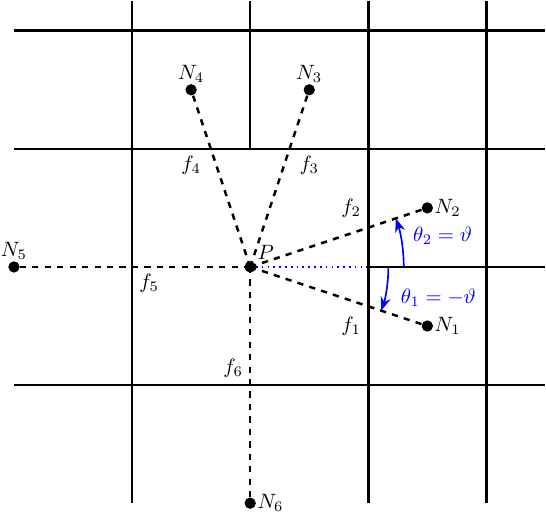}}
\quad
 \subfigure[fine cell] {\label{sfig: composite grid fine cell}
\includegraphics[scale=0.85]{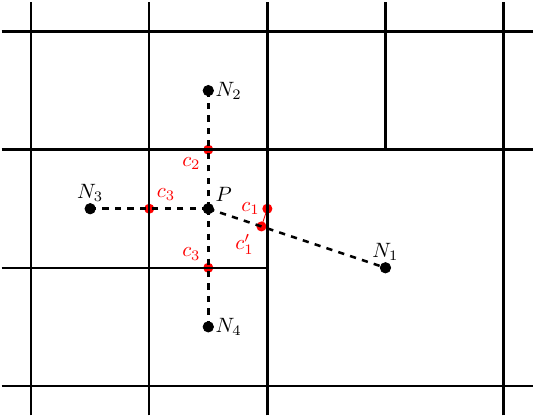}}
 \caption{Topological and geometrical characteristics of a coarse cell \subref{sfig: composite grid 
coarse cell} and of a fine cell \subref{sfig: composite grid fine cell} adjacent to a level 
interface, in a composite grid.}
 \label{fig: refined grid}
\end{figure}

A similar situation may occur when adaptive mesh refinement is used, depending on the treatment of 
the interaction between levels. Figure \ref{fig: grids refined} shows multi-level grids, which 
consist of regions of different fineness. Such grids are often called \textit{composite} grids 
\cite{Trottenberg_2001}. One possible strategy is to treat the cells at the level interfaces as 
topologically polygonal \cite{Muzaferija_1994, Jasak_1996, Syrakos_2012}. For example, cell $P$ of 
Fig.\ \ref{sfig: composite grid coarse cell} has 6 faces, each separating it from a single other 
cell. Its face $f_1$ separates it from cell $N_1$ which belongs to the finer level. Faces such as 
$f_1$, which lie on grid level interfaces, exhibit non-orthogonality, unevenness, and skewness. If 
the grid density is increased throughout the domain, as in the series of grids shown in Fig.\ 
\ref{fig: grids refined}, then these interface distortions remain insensitive to the grid fineness, 
like for the marked line in Fig.\ \ref{fig: grid algebraic}. Alternative schemes exist which avoid 
changing the topology of the cells by inserting a layer of transitional cells between the coarse and 
the fine part of the grid (e.g.\ \cite{Schneiders_1996, Chatzidai_2009}) but they also lead to high, 
non-diminishing grid distortions at the interface.

\begin{figure}[tb]
 \centering
 \subfigure[mean error, $\tau_{\mathrm{mean}}$]
   {\label{sfig: e_mean refined} \includegraphics[scale=1.0]{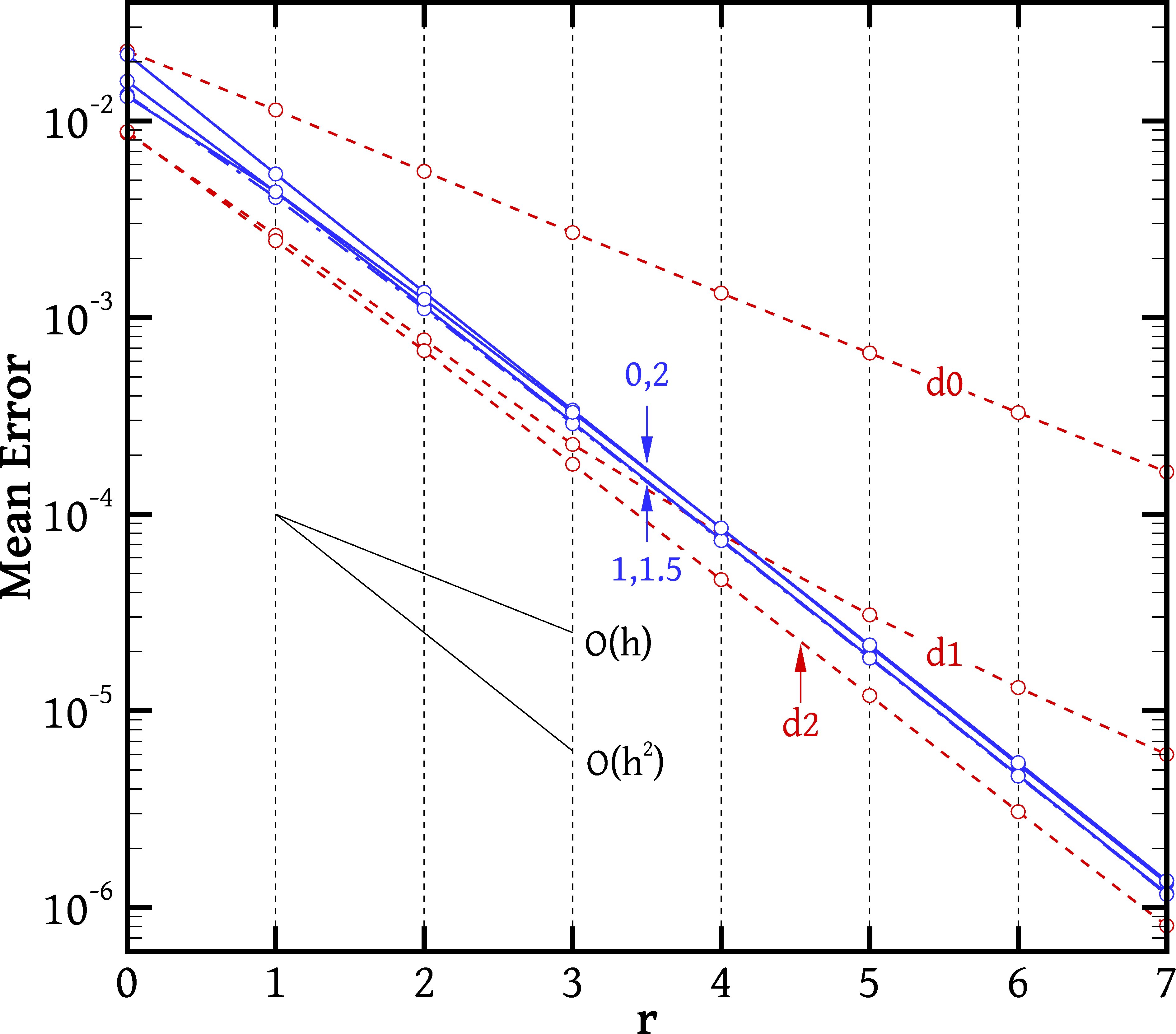}}
 \subfigure[maximum error, $\tau_{\mathrm{max}}$]
   {\label{sfig: e_max refined}  \includegraphics[scale=1.0]{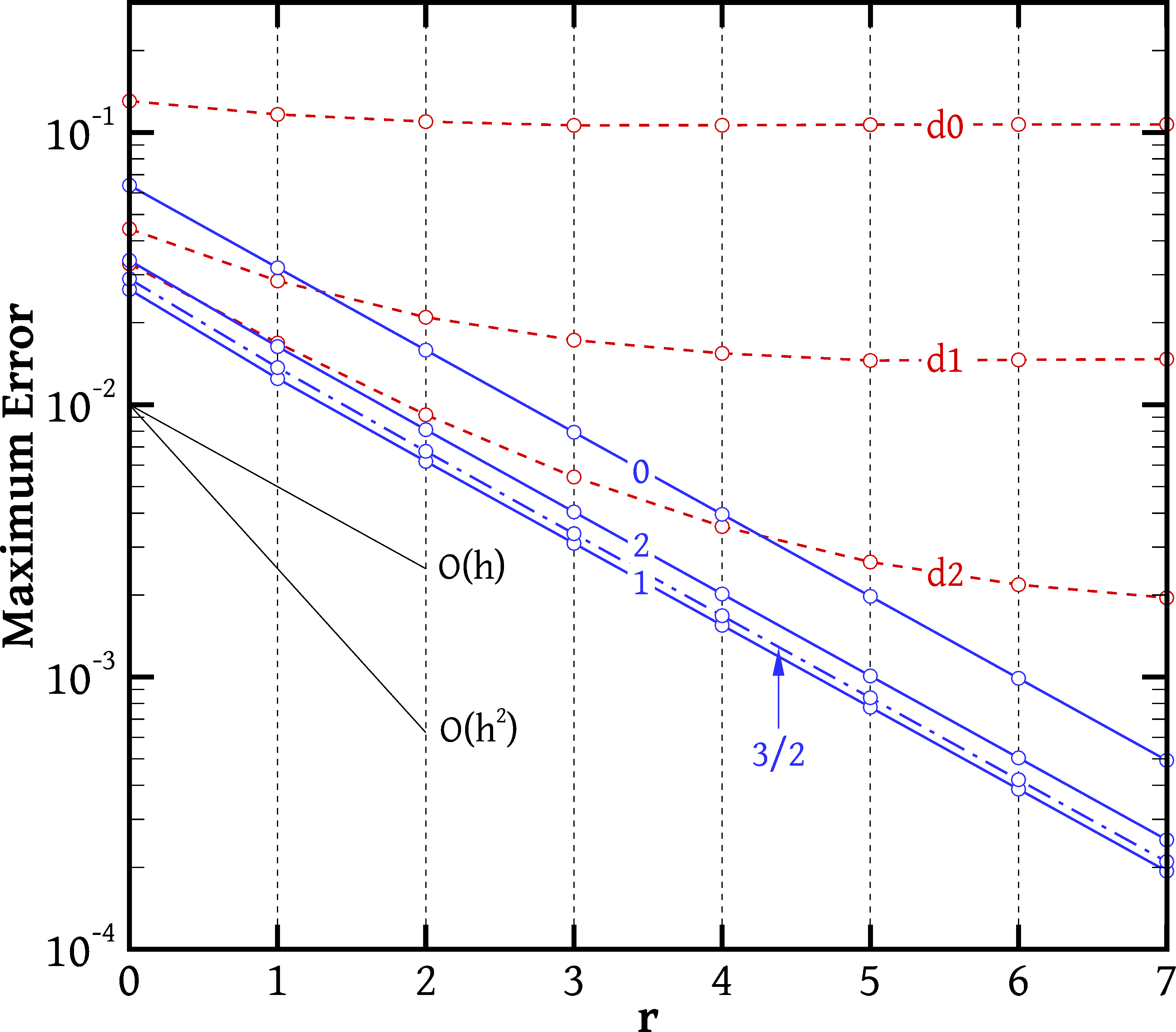}}
 \caption{The mean \subref{sfig: e_mean refined} and maximum \subref{sfig: e_max refined} errors 
(defined by Eqs.\ \eqref{eq: error mean by vol} and \eqref{eq: error max}, respectively) of the 
gradient calculation methods applied to the function $\phi = \tanh(x) \tanh(y)$ on locally refined 
grids (Fig.\ \ref{fig: grids refined}). The abscissa $r$ designates the grid; $r = 0$ is the 
coarsest grid (Fig.\ \ref{sfig: grid refined 0}), and grid $r$ comes from subdividing every cell of 
grid $r-1$ into 4 identical child cells. The blue solid lines correspond to the LS methods with 
weight exponents $q$ = 0, 1 and 2, which are indicated on each curve; the blue dash-dot line 
corresponds to the LS method with $q = 3/2$; and the red dashed lines correspond to the DT methods 
$\mathrm{d}c$ where $c$ is the number of corrector steps.}
 \label{fig: errors refined}
\end{figure}

We computed the gradient of the same function $\phi(x,y) = \tanh(x) \tanh(y)$ on a series of 
composite grids the first three of which are shown in Fig.\ \ref{fig: grids refined}. Figure 
\ref{fig: errors refined} shows how $\tau_{\mathrm{mean}}$ and $\tau_{\mathrm{max}}$ vary 
with grid refinement. This time, $\tau_{\mathrm{mean}}$ is defined a little differently than 
Eq.\ \eqref{eq: error mean} to account for the different grid levels: the error of each individual 
cell is weighted by the cell's volume (i.e.\ the area, in the present two-dimensional setting):
\begin{equation} \label{eq: error mean by vol}
 \tau_{\mathrm{mean}} \;\equiv\; \frac{1}{\Omega} \sum_{P=1}^{M_r} \Omega_P \| 
\nabla^{\mathrm{a}}\phi(\vf{P}) - \nabla 
\phi(\vf{P}) \|
\end{equation}
where $\Omega$ is the total volume of the domain and $\Omega_P$ is the volume of cell $P$.

We can identify three classes of cells that are topologically different. Apart from the familiar 
classes of interior and boundary cells, there is now also the class of cells that touch the level 
interfaces, which shall be called interface cells (these belong to two sub-classes, coarse- and 
fine-level cells, as in Figs.\ \ref{sfig: composite grid coarse cell} and \ref{sfig: composite grid 
fine cell}, respectively). Interface cells possess high skewness and unevenness that do not diminish 
with grid refinement. The behaviour of the gradient-calculation methods at interior and boundary 
cells has already been tested in Sections \ref{ssec: results cartesian} and \ref{ssec: results 
curvilinear}, so our interest now focuses on interface cells. Skewness, the most detrimental grid 
distortion, is encountered only at those and therefore it is there that the maximum errors (Fig.\ 
\ref{sfig: e_max refined}) occur.

In Fig.\ \ref{sfig: e_max refined} we observe that at the level interfaces all the LS methods 
converge to the correct solution at a first-order rate, $\tau_{\mathrm{max}}$ being lowest for 
the $q = 1$ method, followed closely by the $q = 3/2$ and $q = 2$ methods. On the other hand, none 
of the DT methods converge to the correct solution there, although nearly an order of magnitude 
accuracy improvement is obtained with each corrector step. We can determine the operator to whom 
$\nabla^{\mathrm{d}0}$ converges as follows. For the interface cell $P$ of Fig.\ \ref{sfig: 
composite grid fine cell}, formula \eqref{eq: gauss gradient 0} amounts to the following series of  
approximations:
\begin{align}
 \nonumber
 \phi_{\!.x}(\vf{P}) \;&\approx\; \frac{\phi(\vf{c}_1) - \phi(\vf{c}_3)}{h} \;\approx\; 
                             \frac{\phi(\vf{c}'_1) - \phi(\vf{c}_3)}{h}
\\[0.2cm]
 \label{eq: gauss derivative at refined grid}
                \;&\approx\; \frac{\left[ \alpha \, \phi(\vf{N}_1) + (1-\alpha) \, \phi(\vf{P}) 
\right] \;-\;
                                   \left[ 0.5 \, \phi(\vf{P}) + 0.5 \, \phi(\vf{N}_3) \right]}{h}
                \;\equiv\; \phi_{\!.x}^{\mathrm{d}0}(\vf{P})
\end{align}
where $h$ is the length of the side of cell $P$. In the last step, the values of $\phi$ at points 
$\vf{c}'_1$ and $\vf{c}_3$ were approximated with linear interpolation between points $\vf{N}_1$ and 
$\vf{P}$, and $\vf{P}$ and $\vf{N}_3$, respectively; $\alpha$ is an interpolation factor which 
equals $\alpha = 0.3$ for the present geometry. We then substitute in \eqref{eq: gauss derivative at 
refined grid} $\phi(\vf{N}_1)$ and $\phi(\vf{N}_3)$ with their two-dimensional Taylor series about 
$\vf{P}$, considering that if $\vf{P} = (x_0,\, y_0)$ then $\vf{N}_1 = (x_0\!+\!3h/2,\, 
y_0\!-\!h/2)$ and $\vf{N}_3 = (x_0\!-\!h,\, y_0)$ -- see Fig.\ \ref{sfig: composite grid fine cell}. 
The following result is obtained:
\begin{equation*}
 \phi_{\!.x}^{\mathrm{d}0}(\vf{P}) \;=\; 
 \frac{3\alpha + 1}{2} \, \phi_{\!.x}(\vf{P}) \;-\; \frac{\alpha}{2} \, \phi_{\!.y}(\vf{P}) \;+\; 
O(h)
\end{equation*}
Therefore, as $h \rightarrow 0$, $\phi_{\!.x}^{\mathrm{d}0}$ converges not to $\phi_{\!.x}$ but to 
an operator that involves both $\phi_{\!.x}$ and $\phi_{\!.y}$.

Next we examine the mean error in Fig.\ \ref{sfig: e_mean refined}. The plot can be interpreted by 
considering separately the error contributions of each class of cells. The contributions of 
interior and boundary cells to $\tau_{\mathrm{mean}}$ are both $O(h^2)$, as discussed in Section 
\ref{ssec: results cartesian} (for the $q = 3/2$ method the boundary cell contribution is $O(h^3)$).

At interface cells the LS methods behave similarly as on boundary cells, because they produce $O(h)$ 
errors there as well due to unevenness and skewness. The total length of the level interfaces is 
constant, $O(1)$. The number of interface cells is $O(h^{-1})$ because it equals this constant 
length divided by the cell size which is $O(h)$. Their contribution to the mean error in Eq.\ 
\eqref{eq: error mean by vol} is (number of cells) $\times$ (volume of one cell) $\times$ (error at 
a cell) = $O(h^{-1}) \times O(h^2) \times O(h) = O(h^2)$. This is confirmed by Fig.\ \ref{sfig: 
e_mean refined}, where all the LS methods converge to the exact solution at a second-order rate.

On the other hand, for the DT methods the contribution of interface cells to $\tau_{\mathrm{mean}}$ 
is (number of cells) $\times$ (volume of one cell) $\times$ (error at a cell) = $O(h^{-1}) \times 
O(h^2) \times O(1) = O(h)$. Figure \ref{sfig: e_mean refined} shows that for the d0 method (no 
corrector steps) this $O(h)$ component is so large that it dominates $\tau_{\mathrm{mean}}$ even at 
coarse grids. For the d1 method (one corrector step), at coarse grids this $O(h)$ component is 
initially small compared to the bulk $O(h^2)$ component that comes from all the other cells, so that 
$\tau_{\mathrm{mean}}$ appears to decrease at a second-order rate up to a refinement level of $r = 
3$; but eventually it becomes dominant and beyond $r = 5$ the d1 curve is parallel to the d0 curve, 
with a first-order slope. With two corrector steps, the $O(h)$ component is so small that up to $r = 
7$ it is completely masked by the $O(h^2)$ component and the method appears to be second-order 
accurate. More grid refinements are necessary to reveal its asymptotic first-order accuracy.

An observation that raises some concern in Fig.\ \ref{sfig: e_mean refined} is that the DT methods 
outperform the LS methods at those grids where they have not yet degraded to first order. This 
suggests that there may be some room for improvement in the latter. A potential source of the 
problem is suggested by Fig.\ \ref{sfig: composite grid coarse cell}. Along the horizontal 
direction, cell $P$ has one cell on its left side ($N_5$) and two cells on its right side ($N_1$ and 
$N_2$). Since points $\vf{N}_1$ and $\vf{N}_2$ are quite close to each other there is some overlap 
in the information they convey. Yet the weights of the LS method depend only on the distance of 
$\vf{N}_f$ from $\vf{P}$, while any clustering of the $\vf{N}_f$ points in some direction is not 
taken into account. Thus, points $\vf{N}_1$ and $\vf{N}_2$, being closer to $\vf{P}$ than $\vf{N}_5$ 
may individually contribute equally ($q=1$) or more ($q > 1$) to the calculation of the gradient at 
$\vf{P}$ than point $\vf{N}_5$ does. Combined they contribute much more. So, the horizontal 
component of the gradient is calculated using mostly information from the right of cell $P$, whereas 
information from its left is undervalued. The DT methods do not suffer from this deficiency because 
they weigh the contribution of each point by the area of the respective face; faces $f_1$ and $f_2$ 
are half in size than $f_5$, and so points $\vf{N}_1$ and $\vf{N}_2$ \textit{together} contribute to 
the gradient approximately as much as $\vf{N}_5$ alone does.

\begin{figure}[tb]
 \centering
 \subfigure[mean error, $\tau_{\mathrm{mean}}$]
   {\label{sfig: e_mean refined cor} \includegraphics[scale=1.0]{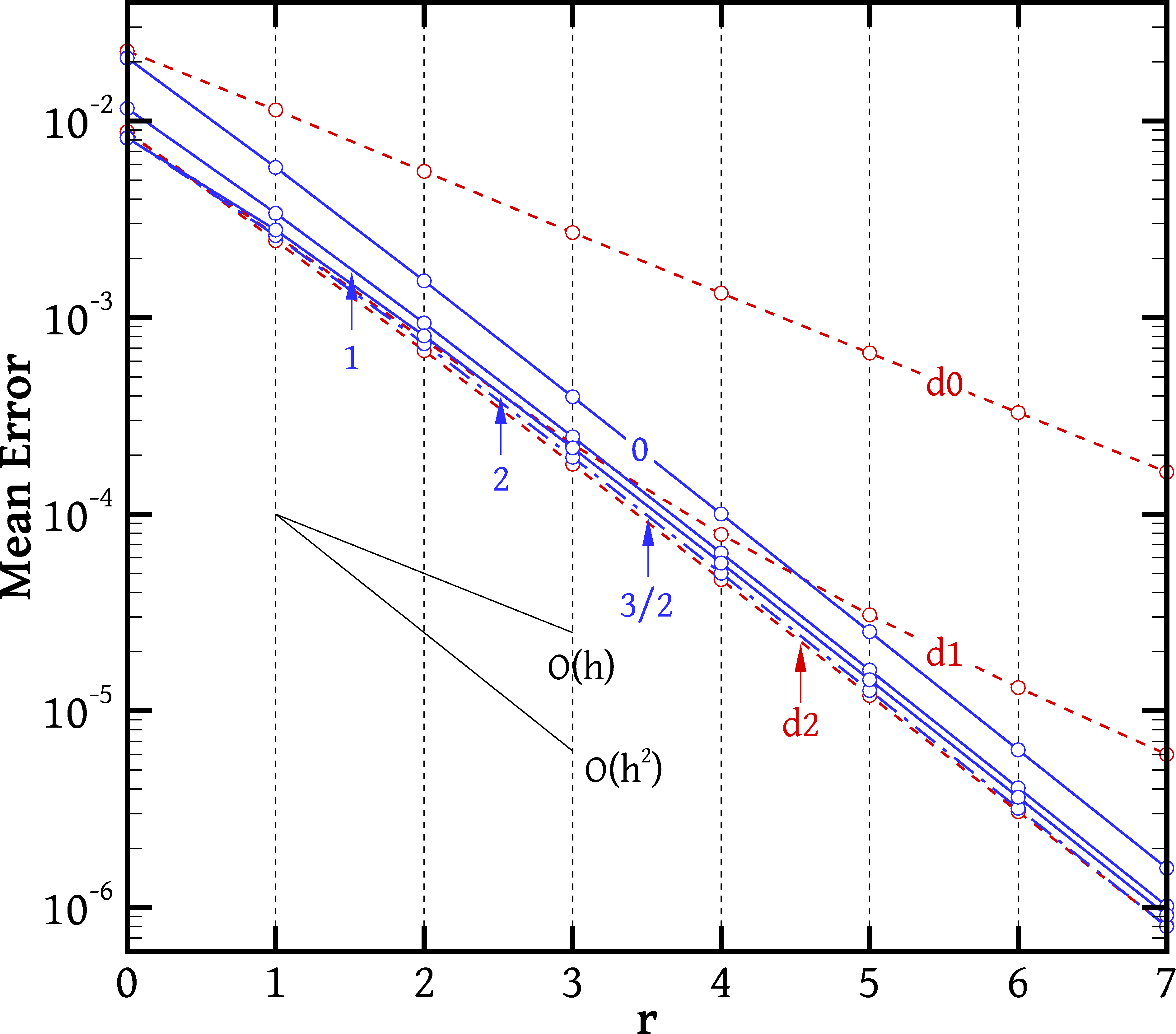}}
 \subfigure[maximum error, $\tau_{\mathrm{max}}$]
   {\label{sfig: e_max refined cor}  \includegraphics[scale=1.0]{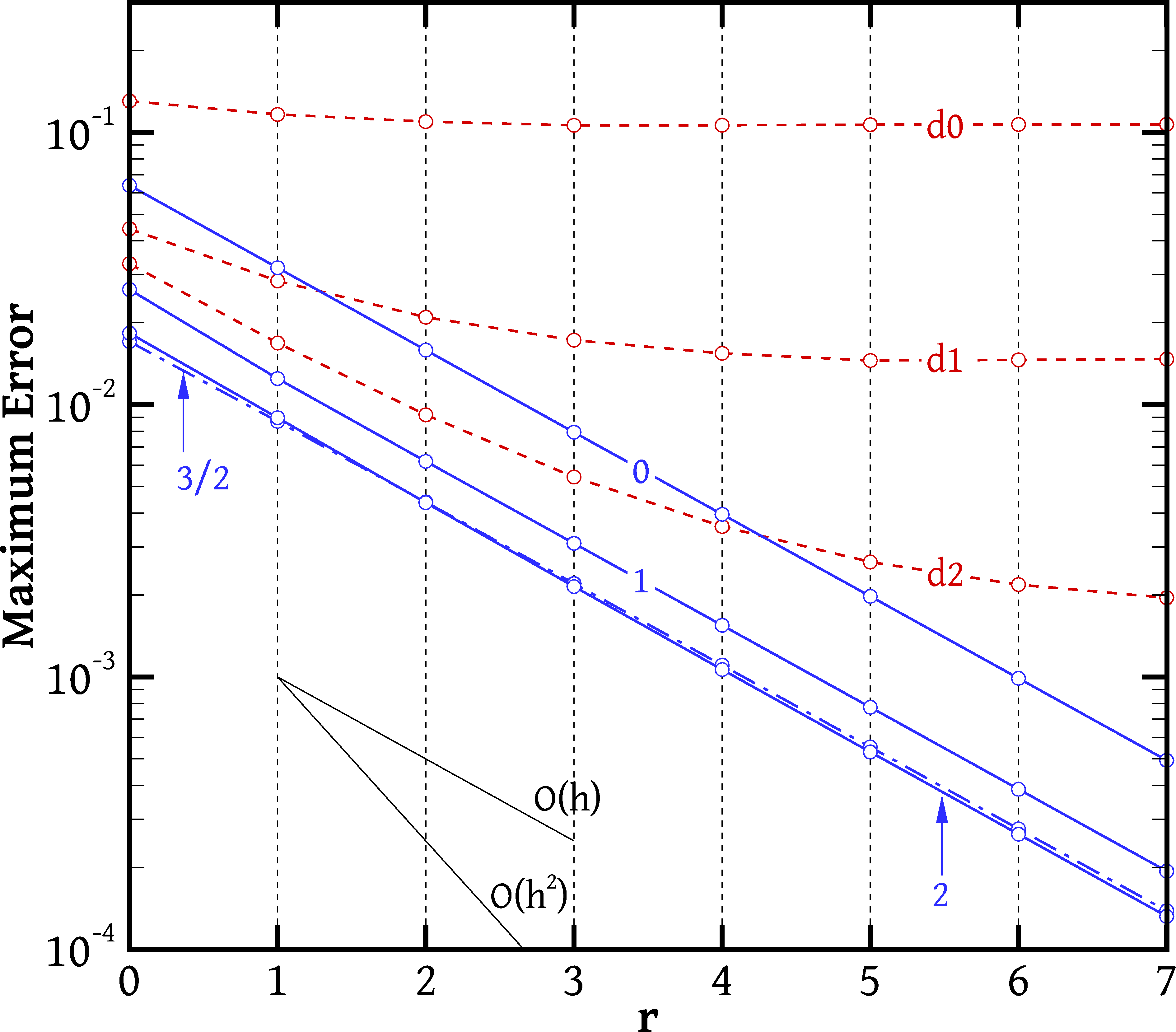}}
 \caption{As for Fig.\ \ref{fig: errors refined}, but with the weights \eqref{eq: weights refined} 
applied at level interfaces in the case of the LS methods.}
 \label{fig: errors refined cor}
\end{figure}

In order to seek a remedy, we investigate the contributions of points $\vf{N}_1$, $\vf{N}_2$ and 
$\vf{N}_5$ (Fig.\ \ref{sfig: composite grid coarse cell}) to the vector $\beta_{\varsigma}$ of the 
error expression \eqref{eq: WLS e}. As for Eq.\ \eqref{eq: beta_tau component 1}, we substitute for 
$\varsigma_f$ from Eq.\ \eqref{eq: varsigma_f} into the expression for $\beta_{\varsigma}$, but 
then proceed to substitute $\Delta x_f = \Delta r_f \cos \theta_f$, $\Delta y_f = \Delta r_f \sin 
\theta_f$. For the particular points under consideration we have, with reference to Fig.\ \ref{sfig: 
composite grid coarse cell}, $\theta_1 = -\vartheta$, $\theta_2 = \vartheta$ and $\theta_5 = \pi$, 
while $\Delta r_2 = \Delta r_1$. Finally, concerning the weights, since cells $P$ and $N_5$ both 
belong to the same grid level we choose not to tamper with the weight $w_5$ and simply use the 
possibly second-order accurate $q = 3/2$ scheme: $w_5 = (\Delta r_5)^{-3/2}$. Due to symmetry we set 
$w_2 = w_1$ and our goal is to suitably select $w_1$ to achieve a small error. Putting everything 
together, the joint contribution of these three points to $\beta_{\varsigma}$ is, neglecting higher 
order terms,
\begin{equation} \label{eq: beta_tau refined}
 \begin{bmatrix}
 \phi_{\!.xx} \left( \left( \cos \vartheta \right)^3 (\Delta r_1)^3 \, w_1^2  \,-\, \frac{1}{2} 
\right) \;+\; 
 \phi_{\!.yy} \left( \sin \vartheta \right)^2 \cos\vartheta \, (\Delta r_1)^3 \, w_1^2
\\[0.4cm]
 2 \phi_{\!.xy} \, \cos\vartheta \left( \sin\vartheta \right)^2 (\Delta r_1)^3 \, w_1^2
 \end{bmatrix}
\end{equation}
Since the values of the higher order derivatives can be arbitrary, it is obvious that the above 
contribution does not become zero for any choice of $w_1$; thus, second-order accuracy cannot be 
achieved. The best that can be done is to choose $w_1 = ( \cos\vartheta )^{-3/2} (\Delta r_1)^{-3/2} 
/ \sqrt{2}$ so that the term in parentheses multiplied by $\phi_{.xx}$ becomes zero. In fact, since 
this choice is not guaranteed to minimise the error and since $\cos \vartheta \approx 1$ for 
relatively small $\vartheta$, we chose to drop the $\cos \vartheta$ factor and applied the 
following scheme to all least squares methods:
\begin{equation} \label{eq: weights refined}
 w_f \;=\;
 \begin{cases}
  \| \vf{N}_f - \vf{P} \|^{-q}                     &
  \qquad  \text{face } f \text{ does not touch a finer level}
\\
  \frac{1}{\sqrt{2}} \| \vf{N}_f - \vf{P} \|^{-q}  &
  \qquad  \text{face } f \text{ touches a finer level}
 \end{cases}
\end{equation}
(In the three-dimensional case, cell $P$ of Fig.\ \ref{sfig: composite grid coarse cell} would have 
four fine-level neighbours on its right and the $1/\sqrt{2}$ factor in \eqref{eq: weights refined} 
would become $1/\sqrt{4} = 1/2$). Using this scheme, we obtained the results shown in Fig.\ 
\ref{fig: errors refined cor}, which can be seen to be better than those of Fig.\ \ref{fig: errors 
refined}, especially for the $q = 3/2$ and $q = 2$ methods which now rival the d2 method with 
respect to the mean error.

The modification \eqref{eq: weights refined} is easy and inexpensive. More elaborate methods could 
be devised, applicable to more general grid configurations, in order to properly deal with the issue 
of point clustering at certain angles. For example, a point $\vf{N}_f$ could be assigned to span a 
``sector'' of angle $\Delta \theta_f = |\theta_{f+1} - \theta_{f-1}| / 2$ (assuming that the points 
are numbered in either clockwise or anti-clockwise order); the sum of all these sectors would then 
equal $2\pi$, and they could be incorporated into the weights such that points with smaller sectors 
would have less influence over the solution.  Alternatively, the face areas could be incorporated 
into the weights, as in the DT method \cite{Shima_2010}.

\subsection{Grids with arbitrary distortion}
\label{ssec: results random}

\begin{figure}[tb]
 \centering
 \subfigure[$r=0$] {\label{sfig: grid random 0} 
\includegraphics[scale=0.65]{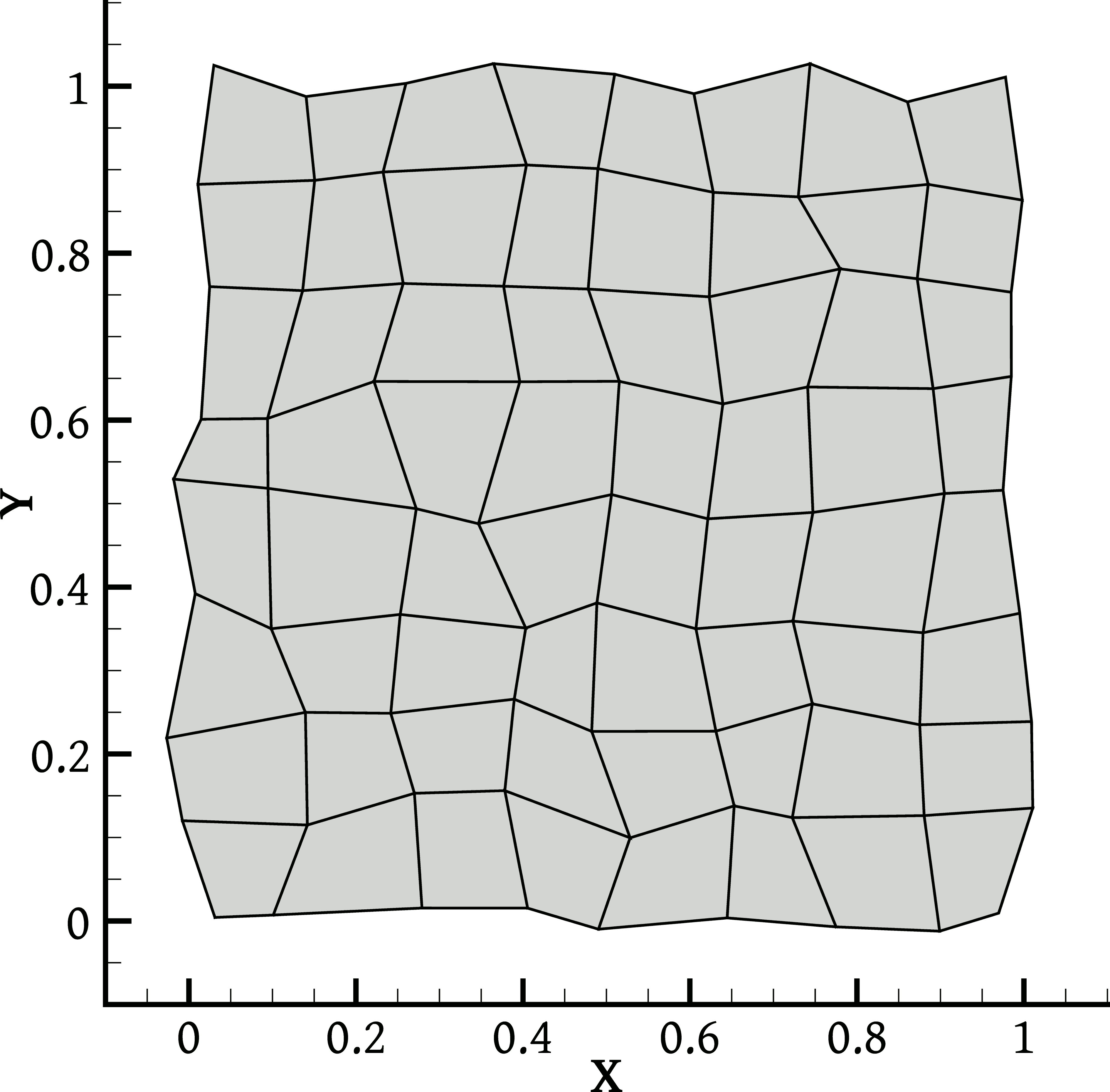}}
 \subfigure[$r=1$] {\label{sfig: grid random 1} 
\includegraphics[scale=0.65]{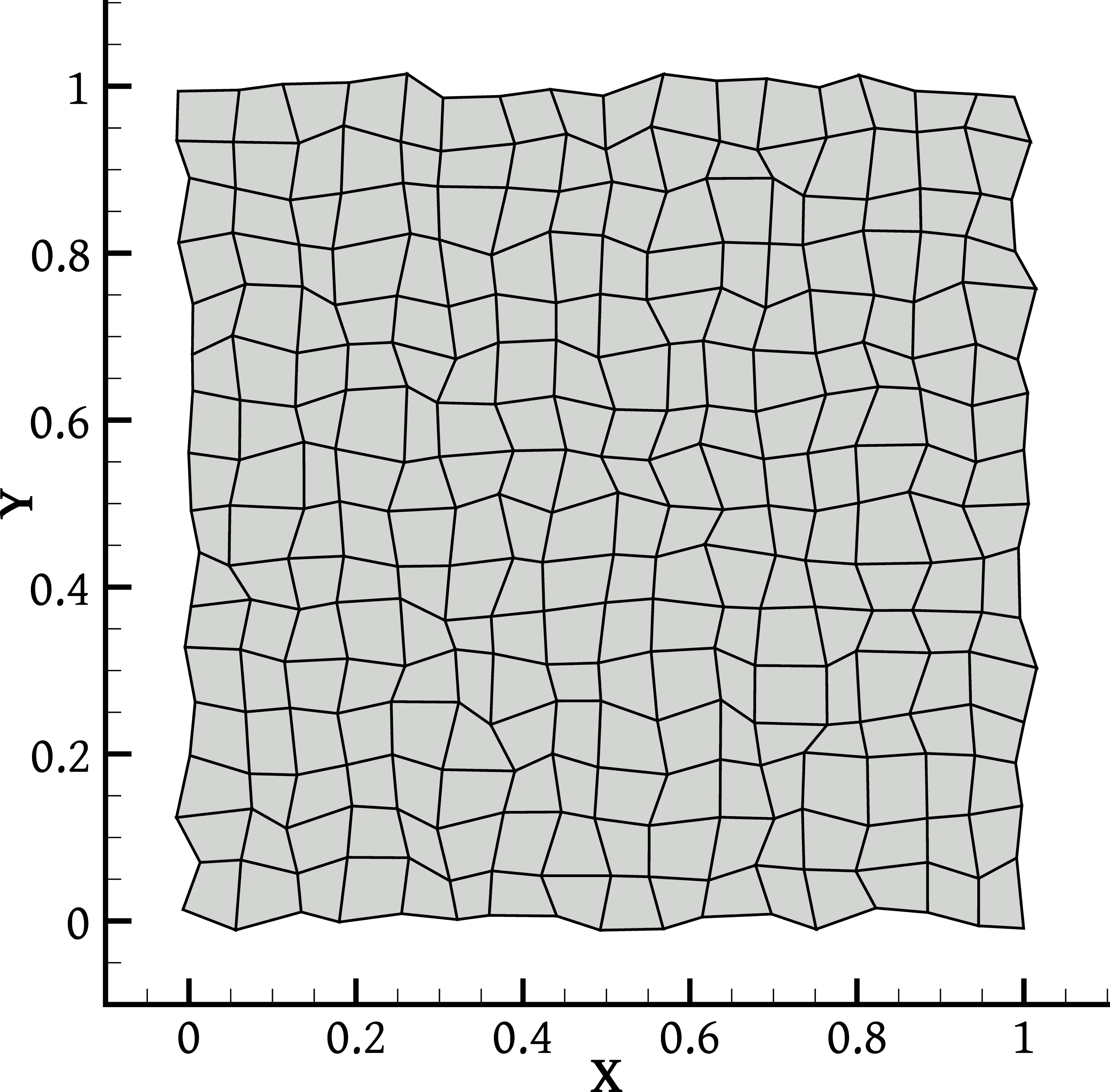}}
 \subfigure[$r=2$] {\label{sfig: grid random 2} 
\includegraphics[scale=0.65]{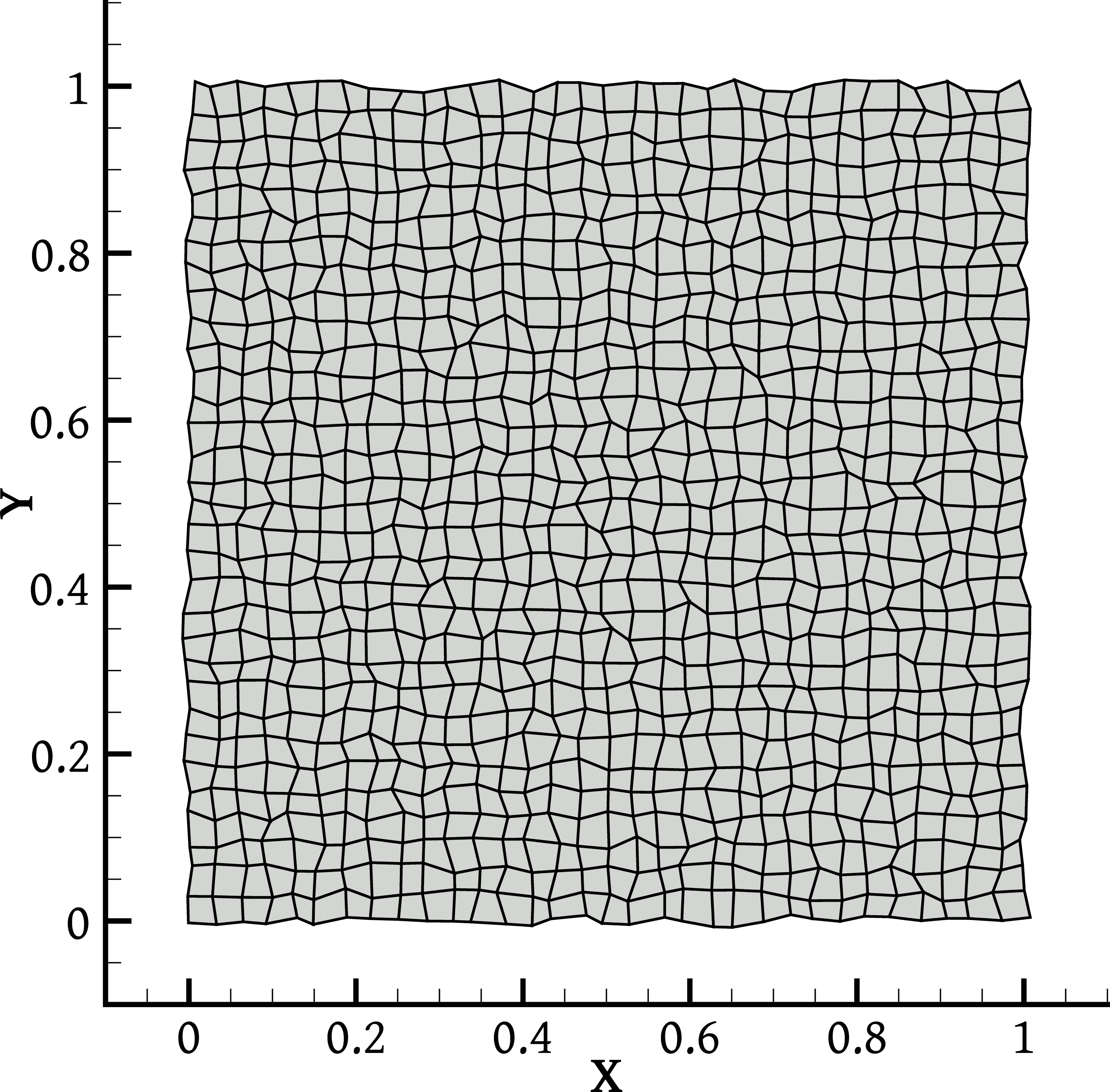}}
 \caption{A series of excessively distorted grids. Grid $r$ is constructed by random perturbation of 
the nodes of Cartesian grid $r+1$ of Section \ref{ssec: results cartesian}.}
 \label{fig: grids random}
\end{figure}

As mentioned in Section \ref{sec: gauss}, general-purpose unstructured grid generation methods 
result in unevenness and skewness that are insensitive to grid fineness throughout the domain, not 
just in isolated regions as for the grids of Section \ref{ssec: results refined}. Therefore, in the 
present Section the methods are tested under these conditions; grids of non-diminishing distortion 
were generated by randomly perturbing the vertices of a Cartesian grid. Using such a process we 
constructed a series of grids, the first three of which are shown in Fig.\ \ref{fig: grids random}. 
The perturbation procedure is applied as follows: Suppose a Cartesian grid with grid spacing $h$. 
If node $(i,j)$ has coordinates $(x_{ij},\, y_{ij})$ then the perturbation procedure moves it to a 
location $(x'_{ij},\, y'_{ij}) = (x_{ij} + \delta^x_{ij},\, y_{ij} + \delta^y_{ij})$ where 
$\delta^x_{ij}$ and $\delta^y_{ij}$ are random numbers in the interval $[-0.25 h, 0.25 h)$. Because 
all perturbations are smaller than $h/4$ in both $x$ and $y$, it is ensured that all grid cells 
remain simple convex quadrilaterals after all vertices have been perturbed. Grids based on triangles 
as well as three-dimensional cases will be considered in Section \ref{sec: PDE solution}.

\begin{figure}[tb]
 \centering
 \subfigure[mean error, $\tau_{\mathrm{mean}}$]
   {\label{sfig: e_mean random} \includegraphics[scale=1.0]{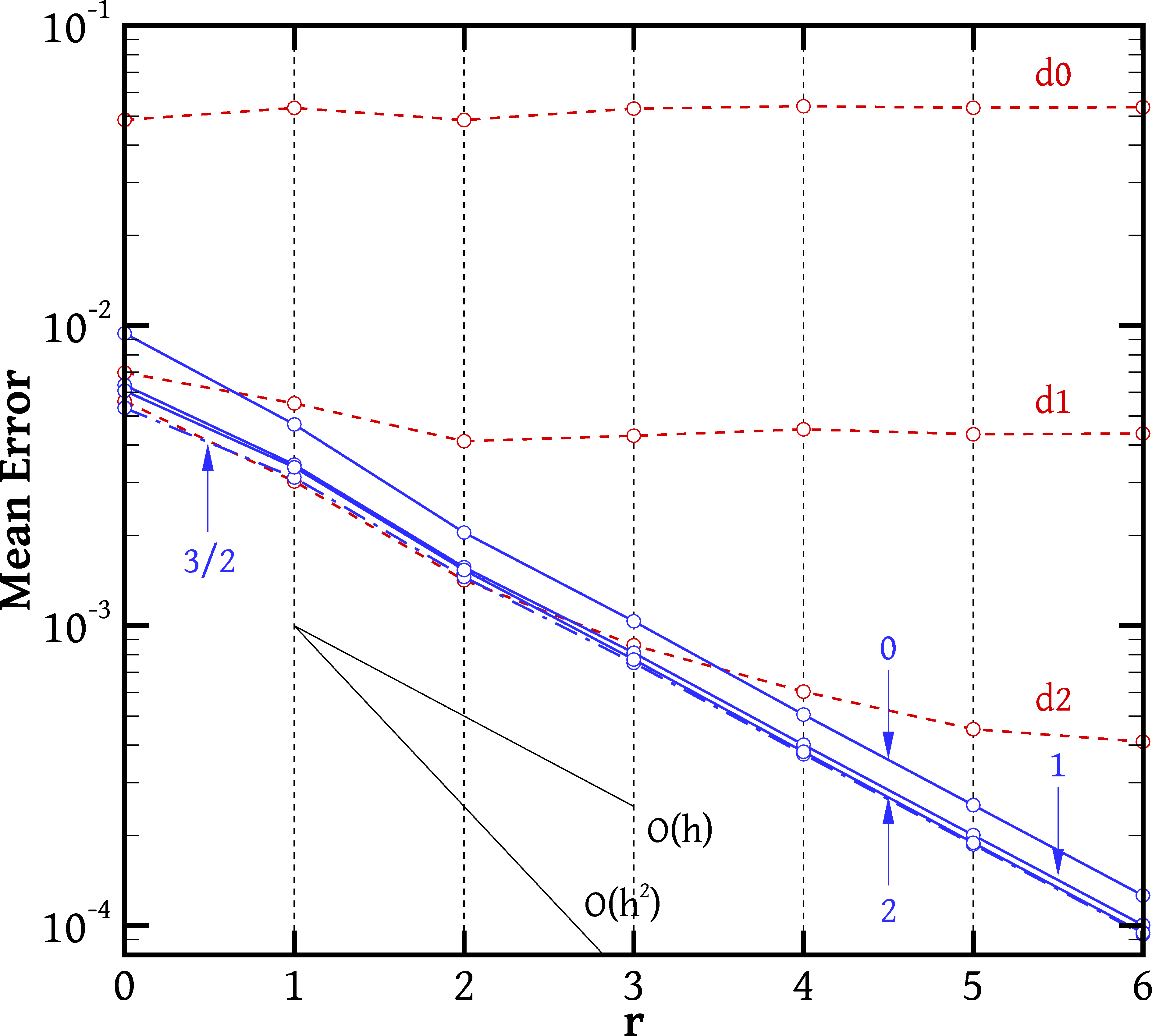}}
 \subfigure[maximum error, $\tau_{\mathrm{max}}$]
   {\label{sfig: e_max random}  \includegraphics[scale=1.0]{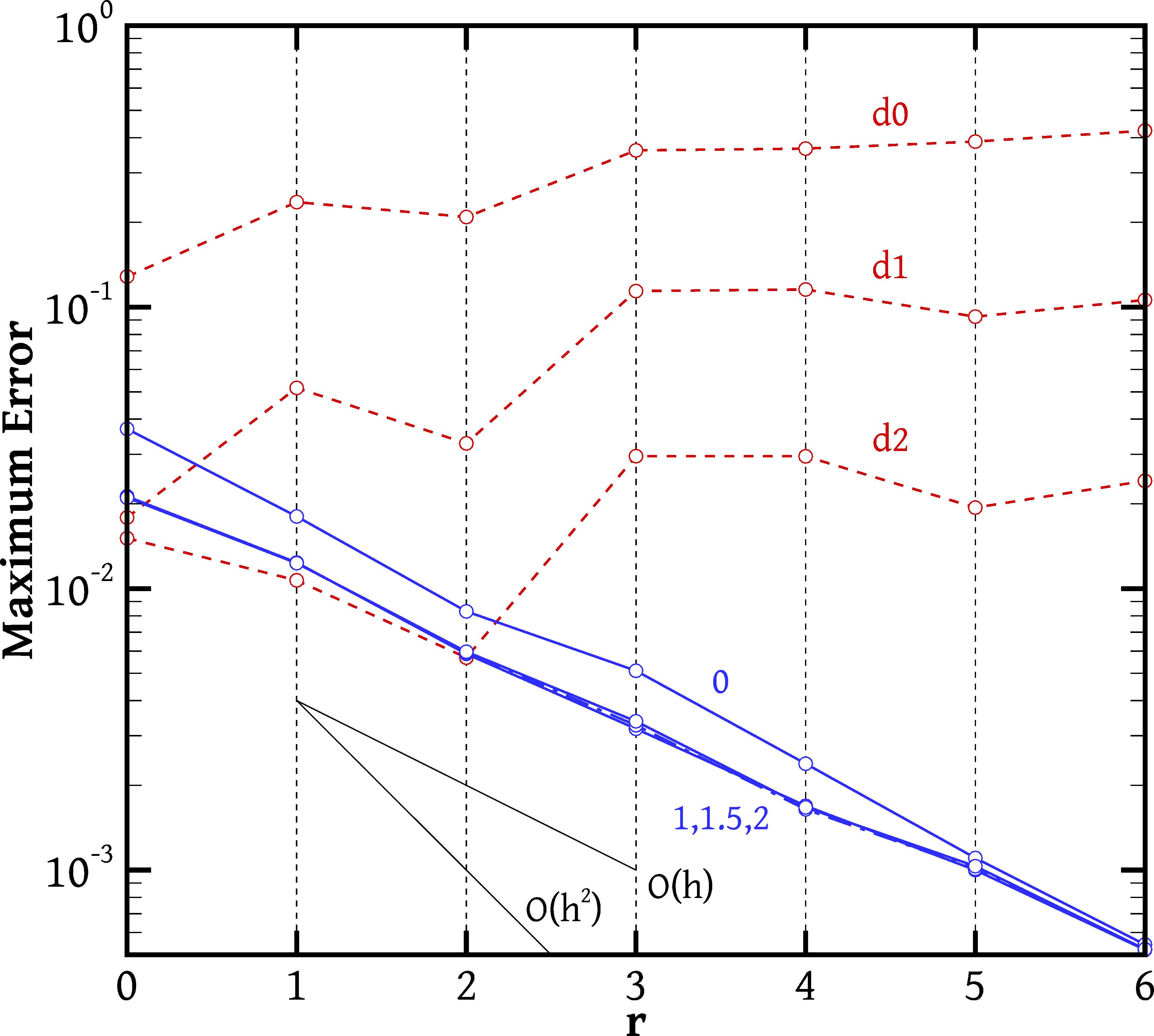}}
 \caption{The mean \subref{sfig: e_mean random} and maximum \subref{sfig: e_max random} errors 
(defined by Eqs.\ \eqref{eq: error mean} and \eqref{eq: error max}, respectively) of the gradient 
calculation methods applied on the function $\phi = \tanh(x) \tanh(y)$ on the series of globally 
distorted grids (Fig.\ \ref{fig: grids random}). The blue solid lines correspond to the least 
squares methods with weight exponents $q$ = 0, 1 and 2, which are indicated on each curve; the blue 
dash-dot line corresponds to the least squares method with $q = 3/2$; and the red dashed lines 
correspond to the divergence theorem methods $\mathrm{d}c$ where $c$ is the number of corrector 
steps.}
 \label{fig: errors random}
\end{figure}

The gradient of the same function $\phi = \tanh(x) \tanh(y)$ is calculated, and the errors are 
plotted in Fig.\ \ref{fig: errors random}. This time, all cells belong to a common category. 
Concerning the mean error, Fig. \ref{sfig: e_mean random} is in full agreement with the theory. The 
LS methods converge to the exact gradient at a first-order rate, with the unweighted method being, 
as usual, the least accurate, while the differences between the weighted methods are very slight. 
On the other hand, the DT methods, as expected, do not converge to the exact gradient. Performing 
corrector steps improves things, but in every case the convergence eventually stagnates at some grid 
fineness.

As far as the maximum errors are concerned, Fig.\ \ref{sfig: e_max random} shows that those of the 
LS methods decrease at a first order rate, with the unweighted method being the least accurate. On 
the other hand, the maximum errors of the DT methods actually increase with grid refinement. 
Presumably this is due to the fact that as the number of grid nodes is increased the probability of 
encountering higher degrees of skewness somewhere in the domain increases. Performing corrector 
steps reduces the error, but it is interesting to note that grid refinement causes a somewhat larger 
error increase when more corrector steps are performed. The deterioration of the 
$\nabla^{\mathrm{d}0}$ method with grid refinement propagates across the iterative correction 
procedure and eventually it is expected that the errors produced by the $\nabla^{\mathrm{d}c}$ 
operator used as a predictor will become so large that it will provide less (rather than more) 
accurate face centre values to the resulting ``corrected'' operator $\nabla^{\mathrm{d}(c+1)}$, 
making the latter worse than $\nabla^{\mathrm{d}c}$ itself.

Similarly to Sec.\ \ref{ssec: results curvilinear} we also tested the gradient schemes on a linear 
function. The results are listed in Table \ref{table: exactness linear random}, together with grid 
quality metrics. Table \ref{table: exactness linear random} confirms that in the present case grid 
skewness and unevenness are roughly independent of grid refinement, i.e.\ their measures are $O(1)$, 
whereas they were $O(h)$ in the structured grid case (Table \ref{table: exactness linear}). As a 
result, the DT gradient schemes with a finite number of corrector steps ($\nabla^{\mathrm{d}0}$ and 
$\nabla^{\mathrm{d}1}$ in Table \ref{table: exactness linear random}) now do not converge to the 
exact gradient. In the limit of infinite corrector steps the DT gradient 
($\nabla^{\mathrm{d}\infty}$) becomes exact for linear functions, as are all the LS gradient schemes 
($\nabla^{\mathrm{ls}0}$ and $\nabla^{\mathrm{ls}1}$ in Table \ref{table: exactness linear random}).

\begin{table}[thb]
\caption{Mean errors (Eq.\ \eqref{eq: error mean}) of various schemes for calculating the gradient 
$\nabla\phi = (1,2)$ of the linear function $\phi(x,y) = x + 2y + 0.5$ on the series of grids shown 
in Fig.\ \ref{fig: grids random}. Also displayed are the measures of grid skewness and unevenness 
(defined in Sec.\ \ref{sec: preliminary}), averaged over all faces of each grid excluding boundary 
faces.}
\label{table: exactness linear random}
\begin{center}
\begin{small}   
\renewcommand\arraystretch{1.25}   
{
\newcommand{\G}[2]{\;$\nabla^{\mathrm{#1}#2}$\;}
\newcommand{\e}[1]{$\cdot 10^{-#1}$}
\begin{tabular}{ c | c c | c c c c c }
\toprule
 Grid $r$ & Skew.      & Unev.     &  \G{d}{0}    &  \G{d}{1}   & \G{d}{\infty} &  \G{ls}{0}  &  \G{ls}{1}  
\\ \midrule
 0        & 4.79\e{2}  & 4.94\e{2} &  1.99\e{1}   & 1.55\e{2}   &  3.83\e{15}   & 9.84\e{16}  & 1.05\e{15}  
\\
 1        & 5.46\e{2}  & 5.32\e{2} &  2.43\e{1}   & 1.96\e{2}   &  7.29\e{15}   & 1.55\e{15}  & 1.59\e{15}  
\\
 2        & 4.97\e{2}  & 5.14\e{2} &  2.17\e{1}   & 1.61\e{2}   &  1.48\e{14}   & 3.07\e{15}  & 3.15\e{15}  
\\
 3        & 5.16\e{2}  & 5.16\e{2} &  2.29\e{1}   & 1.81\e{2}   &  2.94\e{14}   & 6.11\e{15}  & 6.18\e{15}  
\\
 4        & 5.22\e{2}  & 5.17\e{2} &  2.32\e{1}   & 1.89\e{2}   &  5.87\e{14}   & 1.22\e{14}  & 1.23\e{14}  
\\
 5        & 5.18\e{2}  & 5.17\e{2} &  2.32\e{1}   & 1.86\e{2}   &  1.18\e{13}   & 2.45\e{14}  & 2.45\e{14}  
\\
 6        & 5.18\e{2}  & 5.14\e{2} &  2.32\e{1}   & 1.87\e{2}   &  2.34\e{13}   & 4.86\e{14}  & 4.86\e{14}  
\\
\bottomrule
\end{tabular}
} 
\end{small}
\end{center}
\end{table}

\section{Use of the gradient schemes within finite volume PDE solvers}
\label{sec: PDE solution}

So far we have examined the gradient schemes \textit{per se}, examining their truncation error 
through mathematical tools and numerical experiments. Although there are examples of independent use 
of a gradient scheme such as in post-processing, the application of main interest is within finite 
volume methods (FVMs) for the solution of partial differential equations (PDEs). The gradient scheme 
is but a single component of the FVM and how it affects the overall accuracy depends also on the PDE 
solved as well as on the rest of the FVM discretisation. In the present section we provide some 
general comments and some simple demonstrations. The focus is on the effect of the gradient scheme 
on unstructured grids, where the DT scheme was shown to be inconsistent; on such grids the tests of 
Sec.\ \ref{ssec: results random} showed that the weights exponent of the LS method plays a minor 
role and so we only test the $q = 1$ LS variant.

%
%

We begin with the observation that the approximation formula \eqref{eq: gauss gradient 0} is very 
similar to the formulae used by FVMs for integrating convective terms of transport equations over 
grid cells. Therefore, according to the same reasoning as in Section \ref{sec: gauss}, such 
formulae also imply truncation errors of order $O(1)$ on arbitrary grids; this is true even if an 
interpolation scheme other than \eqref{eq: linear interpolation} is used to calculate 
$\bar{\phi}(\vf{c}'_f)$, or even if the exact values $\phi(\vf{c}'_f)$ are known and used. The 
order of the truncation error can be increased to $O(h)$ by accounting for skewness through a 
correction such as \eqref{eq: phi cor}, provided that the gradient used is at least first-order 
accurate. These observations may raise concern about the overall accuracy of the FVM; however, it 
is known that the order of reduction of the discretisation error with grid refinement is often 
greater than that of the truncation error. Thus, $O(1)$ \textit{truncation} errors do not 
necessarily imply $O(1)$ \textit{discretisation} errors; the latter can be of order $O(h)$ or even 
$O(h^2)$. This phenomenon has been observed by several authors (including the present ones 
\cite{Syrakos_2012}) and a literature review can be found in \cite{Diskin_2010}. The question then 
naturally arises of whether and how the accuracy of the gradient scheme would affect the overall 
accuracy of the FVM that uses it. A general answer to this question does not yet exist, and each 
combination of PDE / discretisation scheme / grid type must be examined separately. For simple cases 
such as one-dimensional ones the problem can be tackled using theoretical tools but on general 
unstructured grids this has not yet been achieved \cite{Diskin_2010}. Thus we will resort to some 
simple numerical experiments that amount to solving a Poisson equation.

\subsection{Tests with an in-house solver}
\label{ssec: in-house tests}

We solve the following Poisson equation:
\begin{align}
 & \nabla \cdot \left( -k \nabla \phi \right) \;=\; b(x,y)  & &
\mathrm{on}\; \Omega = [0,1] \!\times\! [0,1]
\label{eq: Poisson}
\\[0.2cm]
 & \phi \;=\; c(x,y) & & \mathrm{on} \;  S_{\Omega}
\label{eq: Poisson BC}
\end{align}
with $k = 1$, where $S_{\Omega}$ is the boundary of the domain $\Omega$ (the unit square), and
\begin{align}
 b(x,y) \;&=\; 2 \tanh(x) \tanh(y) \left[ 2 - (\tanh(x))^2 - (\tanh(y))^2 \right]
\label{eq: Poisson b tanh}
\\[0.2cm]
 c(x,y) \;&=\; \tanh(x) \tanh(y) 
\label{eq: Poisson c tanh}
\end{align}
This is a heat conduction equation with a heat source term and Dirichlet boundary conditions. The 
source term was chosen such that the exact solution of the above equation is precisely $\phi = 
c(x,y)$. The domain was discretised by a series of progressively finer grids of $32 \times 32$, $64 
\times 64$, \ldots, $512 \times 512$ cells denoted as grids 0, 1, \ldots 4, respectively. The grids 
were generated by the same randomised distortion procedure as those of Fig.\ \ref{fig: grids 
random}, only that their boundaries are straight (Fig. \ref{sfig: grid diffusion quads}). 
Corresponding undistorted Cartesian grids were also used for comparison. According to the FV 
methodology, we integrate Eq.\ \eqref{eq: Poisson} over each cell, apply the divergence theorem 
and use the midpoint integration rule to obtain for each cell $P$ a discrete equation of the 
following form:
\begin{equation} \label{eq: Poisson discrete}
 \sum_{f=1}^F D_f \;=\; b(\vf{P}) \, \Omega_P
\end{equation}
where $b(\vf{P})$ is the value of the source term at the centre of cell $P$ and
\begin{equation} \label{eq: D_f}
 D_f \;=\; \int_{S_f} -k \nabla \phi \cdot \vf{n}_f \, \mathrm{d}S
 \;\approx\;
 - S_f \, k \, \nabla \phi(\vf{c}_f) \cdot \vf{n}_f
\end{equation}
is the diffusive flux through face $f$ of the cell. We will test here two alternative 
discretisations of the fluxes $D_f$, both of which utilise some discrete gradient operator.

\begin{figure}[tb]
 \centering
 \subfigure[] {\label{sfig: diffusion 2} \includegraphics[scale=0.80]{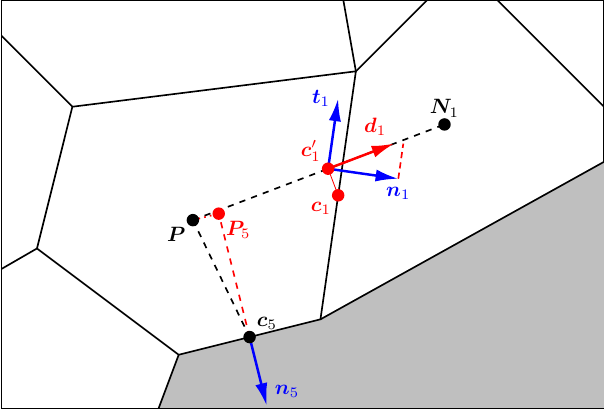}}
 \subfigure[] {\label{sfig: diffusion 1} \includegraphics[scale=0.80]{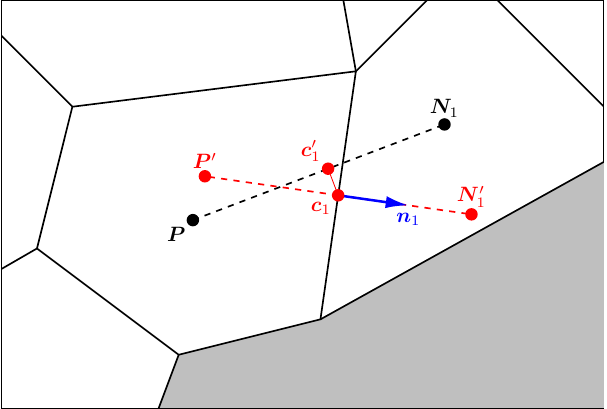}}
 \caption{\subref{sfig: diffusion 2} Adopted notation for the ``over-relaxed'' diffusion scheme 
(Eq.\ \eqref{eq: D_f overrelaxed}) and for the scheme for boundary faces (Eq.\ \eqref{eq: D_f 
boundary}); \subref{sfig: diffusion 1} Notation for the custom scheme \eqref{eq: D_f standard}. 
Face 1 of cell $P$ is an inner face, while its face 5 is a boundary face.}
 \label{fig: diffusion}
\end{figure}

The first is the ``over-relaxed'' scheme \cite{Jasak_1996, Traore_2009} which, according to 
\cite{Demirdzic_2015}, is very popular and is the method of choice in commercial and public-domain 
codes such as FLUENT, STAR-CD, STAR-CCM+ and OpenFOAM. This scheme splits the normal unit vector 
$\vf{n}_f$ in \eqref{eq: D_f} into two components, one in the direction $\vf{N}_f - \vf{P}$ 
and one tangent to the face. The unit vectors along these two directions are denoted as $\vf{d}_f$ 
and $\vf{t}_f$, respectively (Fig.\ \ref{sfig: diffusion 2}). Thus, if $\vf{n}_f = \alpha \vf{d}_f 
+ \beta \vf{t}_f$ then by taking the dot product of this expression with $\vf{n}_f$ itself and 
noting that $\vf{t}_f$ and $\vf{n}_f$ are perpendicular ($\vf{n}_f \cdot \vf{t}_f = 0$) we obtain 
$\alpha = (\vf{d}_f \cdot \vf{n}_f)^{-1}$. We can therefore split $\vf{n}_f = \vf{d}_f^* + 
\vf{t}_f^*$ where $\vf{d}_f^* = \vf{d}_f / (\vf{d}_f \cdot \vf{n}_f)$ and $\vf{t}_f^* = \beta 
\vf{t}_f$ calculated most easily as $\vf{t}_f^* = \vf{n}_f - \vf{d}_f^*$. Then the flux \eqref{eq: 
D_f} is approximated as
\begin{align}
 D_f \;&\approx\; -S_f \, k \, (\nabla \phi(\vf{c}'_f) \cdot
                  \overbrace{ \vf{d}_f) \frac{1}{(\vf{n}_f \cdot \vf{d}_f)} }^{\vf{d}_f^*}
                  \;-\; S_f \, k \, \nabla \phi(\vf{c}'_f) \cdot \vf{t}_f^* \nonumber
\\[0.25cm]
       &\approx\; -S_f \, k \, \frac{\phi(\vf{N}_f)-\phi(\vf{P})}{\|\vf{N}_f - \vf{P}\|} \, 
\frac{1}{(\vf{d}_f \cdot \vf{n}_f)} 
                  \;-\; S_f \, k \, \overline{\nabla^{\mathrm{a}} \phi}(\vf{c}'_f) \cdot \vf{t}_f^* 
\nonumber
\\[0.25cm]
       &=\;       -S_f \, k \, \frac{\phi(\vf{N}_f)-\phi(\vf{P})}{(\vf{N}_f - \vf{P}) \cdot 
\vf{n}_f}
                  \;-\; S_f \, k \, \overline{\nabla^{\mathrm{a}} \phi}(\vf{c}'_f) \cdot \vf{t}_f^*
\label{eq: D_f overrelaxed}
\end{align}
where $\nabla^{\mathrm{a}}$ is a discretised gradient calculated at the cell centres (either a DT 
gradient $\nabla^{\mathrm{d}}$ or the LS gradient $\nabla^{\mathrm{ls}}$) whose role we wish to 
investigate, and $\overline{\nabla^{\mathrm{a}} \phi}(\vf{c}'_f)$ is its value calculated at point 
$\vf{c}'_f$ by linear interpolation \eqref{eq: linear interpolation}. This scheme obviously 
deviates somewhat from the midpoint integration rule, substituting $\nabla \phi(\vf{c}'_f) \cdot 
\vf{n}_f$ instead of $\nabla \phi(\vf{c}_f) \cdot \vf{n}_f$ in \eqref{eq: D_f}, i.e.\ it does not 
account for skewness. 

We also try an alternative scheme which is the standard scheme in our code \cite{Syrakos_2006a} and 
is a slight modification of a scheme proposed in \cite{Ferziger_2002}:
\begin{equation} \label{eq: D_f standard}
 D_f \;\approx\; - S_f \, k \, \frac{\phi(\vf{N}'_f)-\phi(\vf{P}')}{\|\vf{N}'_f - \vf{P}'\|}
\end{equation}
where
\begin{align*}
 \phi(\vf{P}') \;&\approx\; \phi(\vf{P}) \;+\; \nabla^{\mathrm{a}} \phi(\vf{P}) \cdot 
(\vf{P}'-\vf{P})
 \\
 \phi(\vf{N}'_f) \;&\approx\; \phi(\vf{N}_f) \;+\; \nabla^{\mathrm{a}} \phi(\vf{N}_f) \cdot 
(\vf{N}'_f-\vf{N}_f)
\end{align*}
and points $\vf{P}'$ and $\vf{N}'_f$ (Fig.\ \ref{sfig: diffusion 1}) are such that the line segment 
joining these two points has length $\|\vf{N}_f - \vf{P}\|$, is perpendicular to face $f$, and its 
midpoint is $\vf{c}_f$. Thus this scheme tries to account for skewness. We will refer to this as 
the ``custom'' scheme.

Finally, irrespective of whether scheme \eqref{eq: D_f overrelaxed} or \eqref{eq: D_f standard} is 
used for the inner face fluxes, the boundary face fluxes are always calculated as
\begin{equation} \label{eq: D_f boundary}
 D_f \;\approx\; - S_f \, k \, \frac{\phi(\vf{c}_f) - \phi(\vf{P}_f)}{\|\vf{c}_f - \vf{P}_f \|}
\end{equation}
where $f$ is a boundary face, point $\vf{P}_f$ is the projection of $\vf{P}$ on the line which is 
perpendicular to face $f$ and passes through $\vf{c}_f$ (see Fig.\ \ref{sfig: diffusion 2}), and
\begin{equation*}
 \phi(\vf{P}_f) \;\approx\; \phi(\vf{P}) \;+\; \nabla^{\mathrm{a}} \phi(\vf{P}) \cdot 
(\vf{P}_f-\vf{P})
\end{equation*}

Grid non-orthogonality activates in all schemes, \eqref{eq: D_f overrelaxed}, \eqref{eq: D_f 
standard} and \eqref{eq: D_f boundary}, a component that involves the approximate gradient 
$\nabla^{\mathrm{a}}$ (it is activated in \eqref{eq: D_f standard} also by unevenness or skewness). 
It is not hard to show that this component contributes $O(1)$ to the truncation error if 
$\nabla^{\mathrm{a}}$ is first-order accurate, and $O(1/h)$ if it is zeroth-order accurate. On the 
undistorted Cartesian grids that were used for comparison the $\nabla^{\mathrm{a}}$ terms of the 
schemes are not activated; in fact, schemes \eqref{eq: D_f overrelaxed} and \eqref{eq: D_f standard} 
reduce to the same simple formula there. The reason why we do not simply substitute 
$\overline{\nabla^{\mathrm{a}}\phi}(\vf{c}_f)$ for $\nabla \phi(\vf{c}_f)$ in Eq.\ \eqref{eq: D_f} 
instead of using schemes such as \eqref{eq: D_f overrelaxed} and \eqref{eq: D_f standard} is that 
this would allow spurious high-frequency (cell-to-cell) oscillations in the numerical solution 
$\phi$. This occurs because the $\nabla^{\mathrm{a}}$ calculation is mostly based on differences 
of $\phi$ values that are two cells apart, and so the gradient field $\nabla^{\mathrm{a}}\phi$ can 
be smooth even though $\phi$ oscillates from cell to cell. Schemes such as \eqref{eq: D_f 
overrelaxed} and \eqref{eq: D_f standard} express the diffusive flux mostly as a function of the 
direct difference of $\phi$ across the face, $\phi(\vf{N}_f) - \phi(\vf{P})$, and thus do not allow 
such oscillations to go undetected. The gradients $\nabla^{\mathrm{a}}$ are used in an auxiliary 
fashion, in order to make the diffusive fluxes consistent. However, it has recently been shown 
\cite{Nishikawa_2011, Jalali_2014, Nishikawa_2017} that schemes such as \eqref{eq: D_f overrelaxed} 
and \eqref{eq: D_f standard} can also be interpreted as equivalent to substituting  an interpolated 
value of $\nabla^{\mathrm{a}}\phi$ for $\nabla\phi(\vf{c}_f)$ in \eqref{eq: D_f} and adding a 
damping term which is a function of the direct difference $\phi(\vf{N}_f) - \phi(\vf{P})$ and of the 
gradients $\nabla^{\mathrm{a}}\phi(\vf{P})$ and $\nabla^{\mathrm{a}}\phi(\vf{N}_f)$ and tends to 
zero with grid refinement.

\begin{table}[b]
\caption{Number of outer iterations for convergence (residual per unit volume below $10^{-8}$) of 
the over-relaxed \eqref{eq: D_f overrelaxed} and custom \eqref{eq: D_f standard} schemes for 
solving the Poisson equation \eqref{eq: Poisson} -- \eqref{eq: Poisson c tanh} on the $512 \times 
512$ distorted grids, employing various gradient schemes. Also shown is the percentage of CPU time 
spent on computing the gradient.}
\label{table: FVM convergence}
\begin{center}
\begin{small}   
\renewcommand\arraystretch{1.25}   
{
\newcommand{\G}[2]{\;$\nabla^{\mathrm{#1}#2}$\;}
\begin{tabular}{ l l | c c c c c c }
\toprule
 Scheme       &             &  \G{d}{0} & \G{d}{1} & \G{d}{2} & \G{d}{\infty} & \G{dx}{} & \G{ls}{}
\\ \midrule
 Over-relaxed & iterations  &       901 &      617 &      601 &           786 &      573 &   585
\\
              & CPU \%      &       6.8 &     14.9 &     22.1 &           9.7 &     17.2 &  14.3
\\ \midrule
 Custom       & iterations  &       792 &      874 &      914 &           752 &      919 &   919
\\
              & CPU \%      &       7.4 &     16.5 &     23.6 &          10.5 &     17.9 &  15.2
\\
\bottomrule
\end{tabular}
} 
\end{small}
\end{center}
\end{table}

The system of discrete equations \eqref{eq: Poisson discrete} is linear, but for convenience it was 
solved with a fixed-point iteration procedure where in each outer iteration a linear system is 
solved (by a few inner iterations of a preconditioned conjugate gradient solver) whose equations 
involve only the unknowns at a cell and at its immediate neighbours, thus avoiding extended 
stencils. The matrix of this linear system is assembled only from the parts of Eq.\ \eqref{eq: D_f 
overrelaxed} or \eqref{eq: D_f standard} that directly involve $\phi(\vf{P})$ and $\phi(\vf{N}_f)$, 
while the terms involving the gradients are calculated using the estimate of $\phi$ from the 
previous outer iteration and incorporated into the right-hand side of the linear system, as is 
customary. Outer iterations were carried out until the magnitude of the residual per unit volume 
had fallen below $10^{-8}$ in every grid cell, and the number of required such iterations for each 
diffusion flux scheme and gradient scheme are listed  in Table \ref{table: FVM convergence}, where 
the operator $\nabla^{\mathrm{d}\infty}$ is obtained with the scheme \eqref{eq: gauss outer 
iterative scheme}, while the operator $\nabla^{\mathrm{dx}}$ is the one obtained by applying the 
divergence theorem to the auxiliary cell of Fig.\ \ref{fig: dual cell}. Table \ref{table: FVM 
convergence} includes the percentage of the total calculation time consumed in calculating the 
gradient. It may be seen that the cost of corrector steps is quite significant, with the 
$\nabla^{\mathrm{d}2}$ gradient requiring more than 22\% of the total computational effort. On the 
other hand, the scheme \eqref{eq: gauss outer iterative scheme} is quite efficient, obtaining the 
$\nabla^{\mathrm{d}\infty}$ gradient at a cost that is almost as low as that of 
$\nabla^{\mathrm{d}0}$; however, in the case of the over-relaxed scheme, it also requires somewhat 
more outer iterations for convergence. The LS gradient costs about the same as the 
$\nabla^{\mathrm{d}1}$ gradient. The cost of the $\nabla^{\mathrm{dx}}$ gradient appears somewhat 
inflated due to a quick but not very efficient implementation. Overall, it seems that the iterative 
convergence rate depends mostly on the chosen flux discretisation scheme rather than on the choice 
of gradient scheme: in most cases the over-relaxed scheme converges faster than the custom scheme, 
which is not surprising since the former is known for its robustness \cite{Demirdzic_2015}. 
Nevertheless, the choice of gradient scheme may significantly affect the iterative convergence rate 
in the FVM solution of other kinds of problems (depending of course also on the choice of FVM 
discretisation and iterative solver). In \cite{Diskin_2011} it is reported that for convection 
problems unweighted LS gradients can result in much better iterative convergence than weighted ones, 
and that including more distant neighbours in the LS gradient computation also improves convergence.

\begin{figure}[tb]
 \centering
 \subfigure[mean error, $\epsilon_{\mathrm{mean}}$]
   {\label{sfig: e_mean diffusion} \includegraphics[scale=1.0]{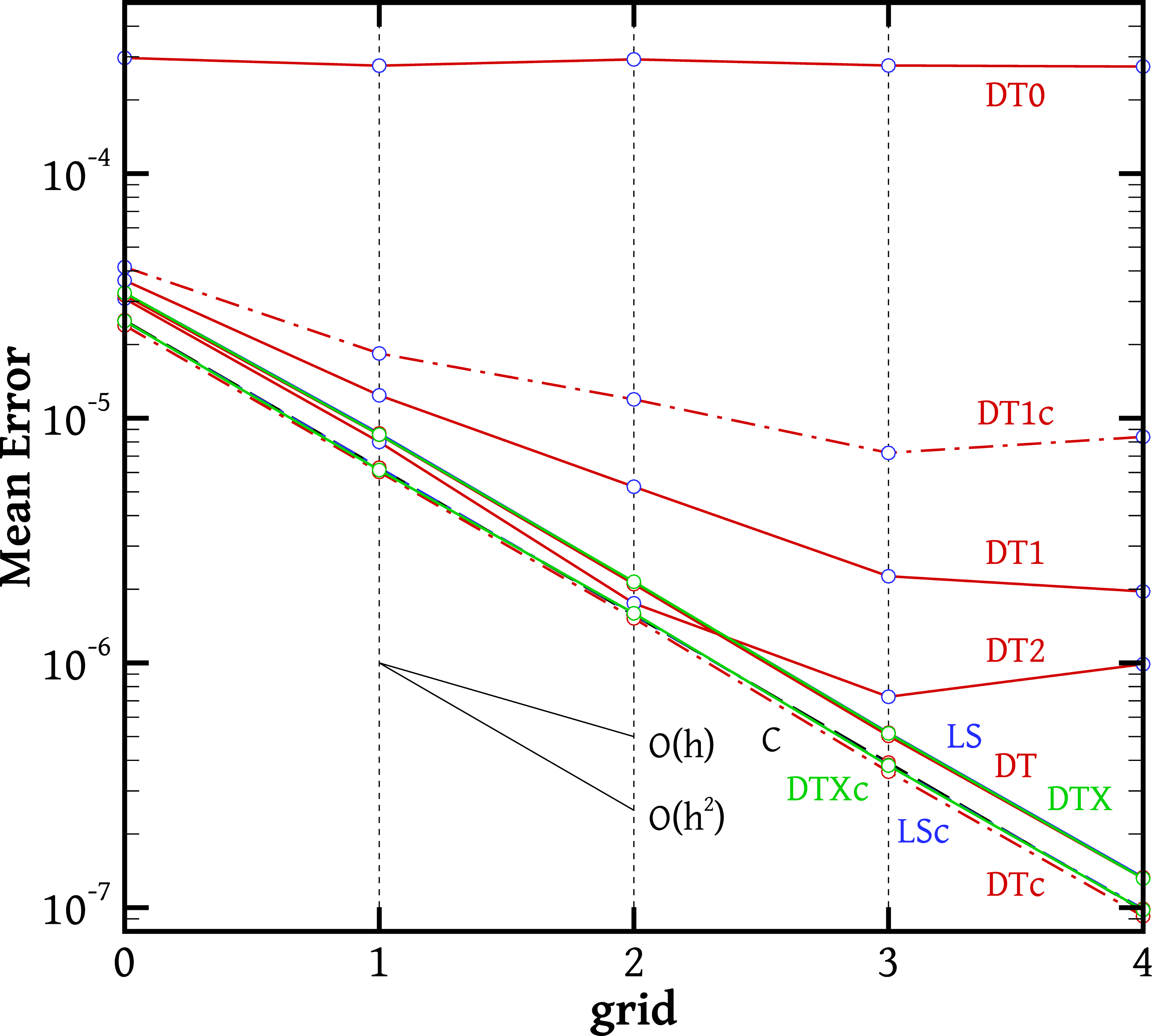}}
 \subfigure[maximum error, $\epsilon_{\mathrm{max}}$]
   {\label{sfig: e_max diffusion}  \includegraphics[scale=1.0]{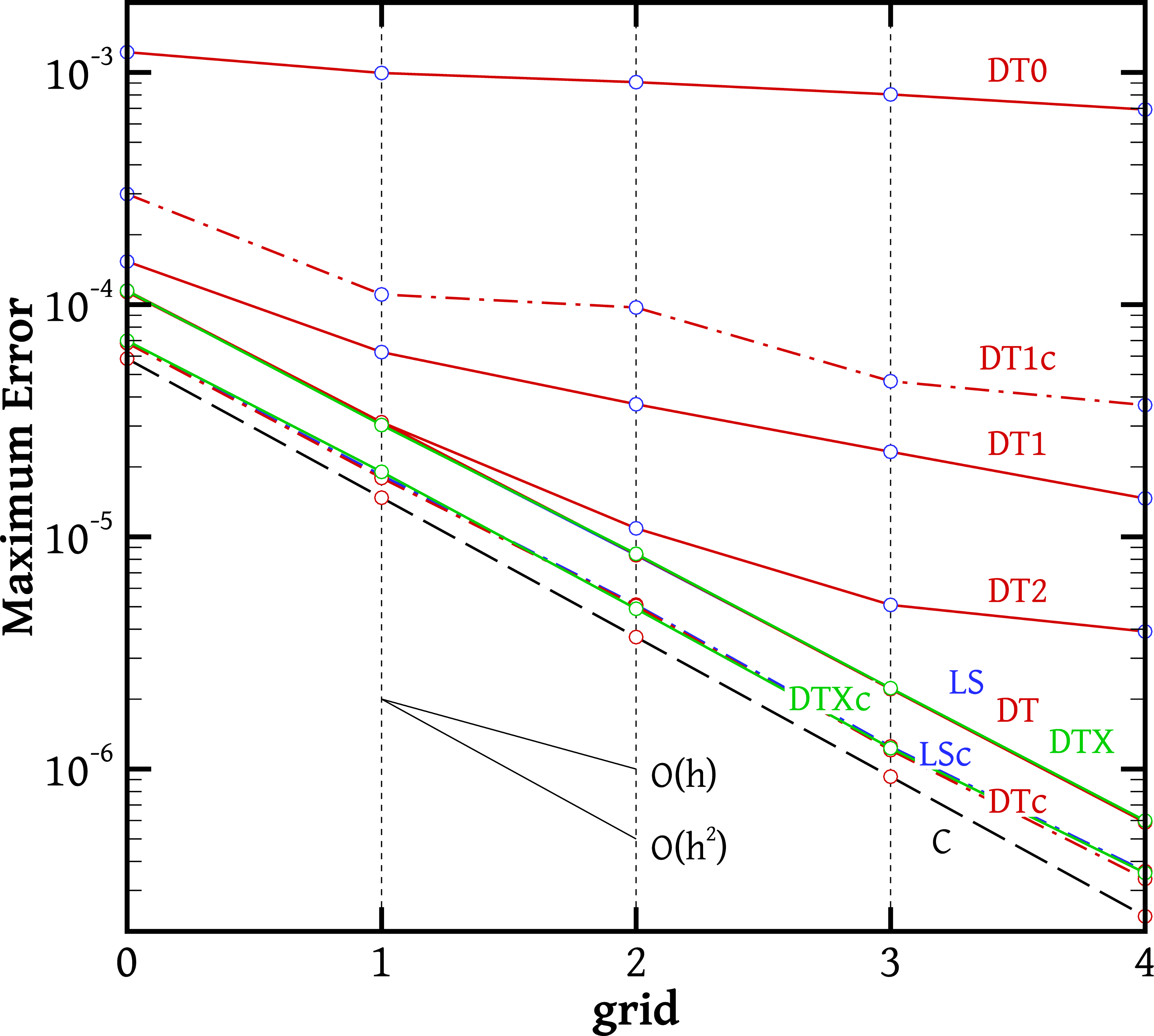}}
 \caption{The mean \subref{sfig: e_mean diffusion} and maximum \subref{sfig: e_max diffusion} 
discretisation errors of various FVM schemes to solve the diffusion problem \eqref{eq: Poisson} -- 
\eqref{eq: Poisson c tanh}, on a series of highly distorted quadrilateral grids (Fig.\ \ref{sfig: 
grid diffusion quads}). Grids 0, 1, \ldots, 4 have $32 \times 32$, $64 \times 64$, \ldots, $512 
\times 512$ volumes, respectively. DT0, DT1, DT2, DT, DTX and LS (solid lines) denote the 
``over-relaxed'' scheme \eqref{eq: D_f overrelaxed} with the $\nabla^{\mathrm{d}0}$, 
$\nabla^{\mathrm{d}1}$, $\nabla^{\mathrm{d}2}$, $\nabla^{\mathrm{d}\infty}$, $\nabla^{\mathrm{dx}}$ 
and $\nabla^{\mathrm{ls}}$ gradient schemes, respectively. DT1c, DTc, DTXc and LSc (dash-dot lines) 
denote the ``custom'' scheme \eqref{eq: D_f standard} with the respective gradient schemes. C 
denotes results on Cartesian grids, where all methods are identical.}
 \label{fig: errors diffusion}
\end{figure}

The discretisation errors with respect to grid refinement are plotted in Fig.\ \ref{fig: errors 
diffusion} for various flux and gradient scheme combinations. The discretisation error is defined 
as the difference between the exact solution $\phi$ of the problem \eqref{eq: Poisson} -- 
\eqref{eq: Poisson BC} and the numerical (FVM) solution; Fig.\ \ref{sfig: e_mean diffusion} plots 
the mean absolute value of the discretisation error over all grid cells, and Fig.\ \ref{sfig: e_max 
diffusion} plots the maximum absolute value.

The first thing that one can notice from Figure \ref{fig: errors diffusion} is that the common DT 
gradient ($\GRAD{d0}$, $\GRAD{d1}$, $\GRAD{d2}$) leads to zeroth-order accuracy with respect to 
both the mean and maximum discretisation errors, irrespective of whether the over-relaxed or custom 
flux scheme is used. Similar trends as those of Fig.\ \ref{fig: errors random} are observed: 
corrector steps reduce the error but are eventually unable to converge to the exact solution. 
Although on coarse grids the DT1 and DT2 lines give the impression that they converge, eventually 
grid refinement cannot reduce the error below a certain point. In fact, an increase of the mean 
error can be observed on grid 4. Things are worse when the DT gradient is used in combination with 
the custom scheme (line DT1c), presumably because the gradient $\nabla^{\mathrm{a}}$ plays a more 
significant role in that scheme than in the over-relaxed scheme.

On the other hand, the $\GRAD{d\infty}$ gradient obtained through the iterative procedure 
\eqref{eq: gauss outer iterative scheme}, the ``auxiliary cell'' gradient $\GRAD{dx}$, and the LS 
gradient $\GRAD{ls}$, all lead to second-order convergence to the exact solution. As long as the 
gradient scheme is consistent, the FVM accuracy seems to depend not on the gradient scheme but on 
the flux scheme, i.e.\ lines DT, LS and DTX (over-relaxed scheme \eqref{eq: D_f overrelaxed}) 
nearly coincide, as do lines DTc, LSc and DTXc (custom scheme \eqref{eq: D_f standard}). The mean 
discretisation error of the custom scheme is about the same as that obtained on Cartesian grids 
(line C), or even marginally lower when used in combination with the $\GRAD{d\infty}$ gradient! The 
mean discretisation error of the over-relaxed scheme is only about 40\% higher, which is a very 
good performance given that it does not account for skewness and that its iterative convergence is 
better (Table \ref{table: FVM convergence}).

\subsection{Tests with OpenFOAM}
\label{ssec: openfoam tests}

\begin{figure}[tb]
 \centering
 \subfigure[] {\label{sfig: grid diffusion quads} 
\includegraphics[scale=0.7]{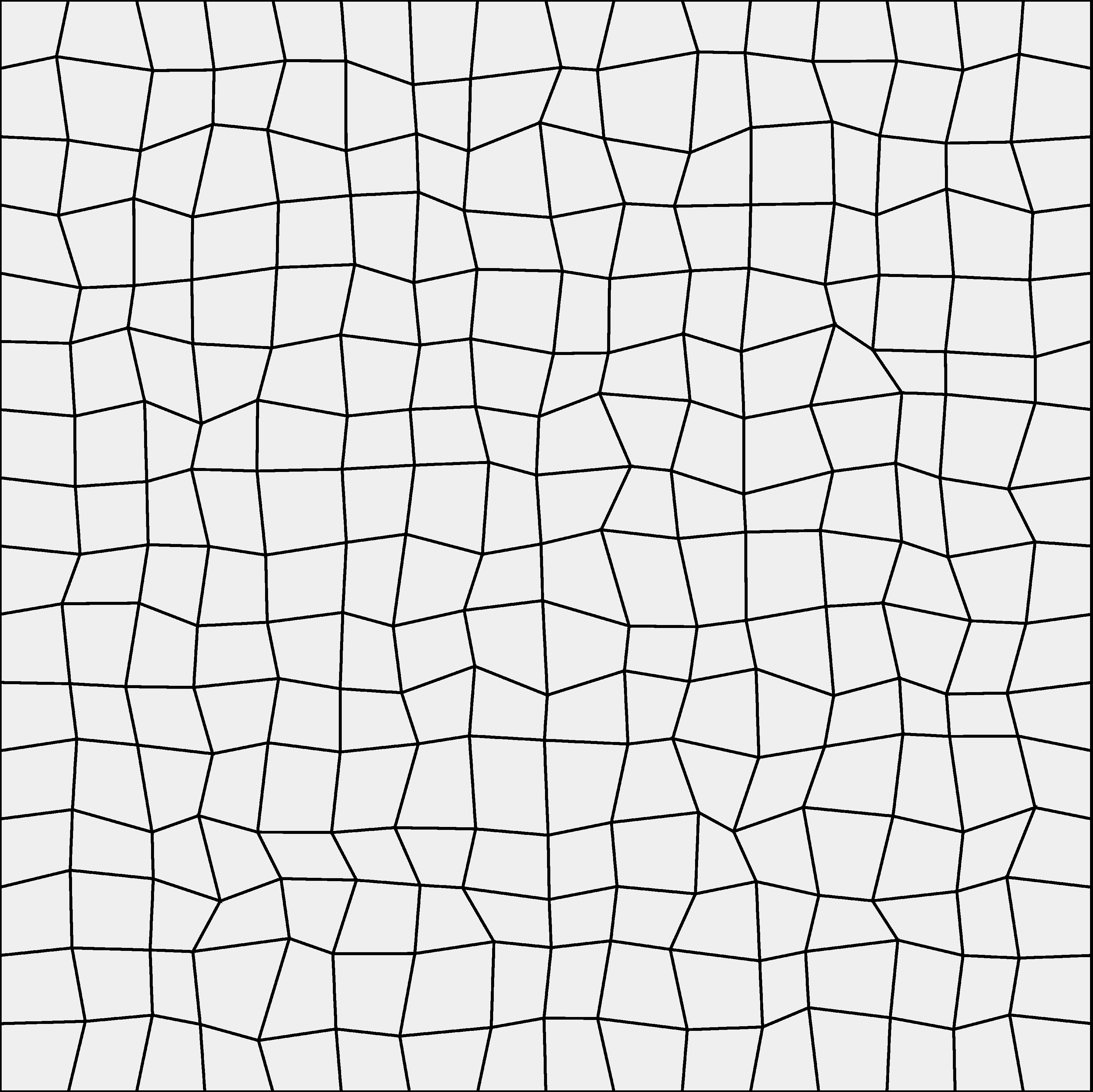}}
 \quad
 \subfigure[] {\label{sfig: grid diffusion netgen}   
\includegraphics[scale=0.7]{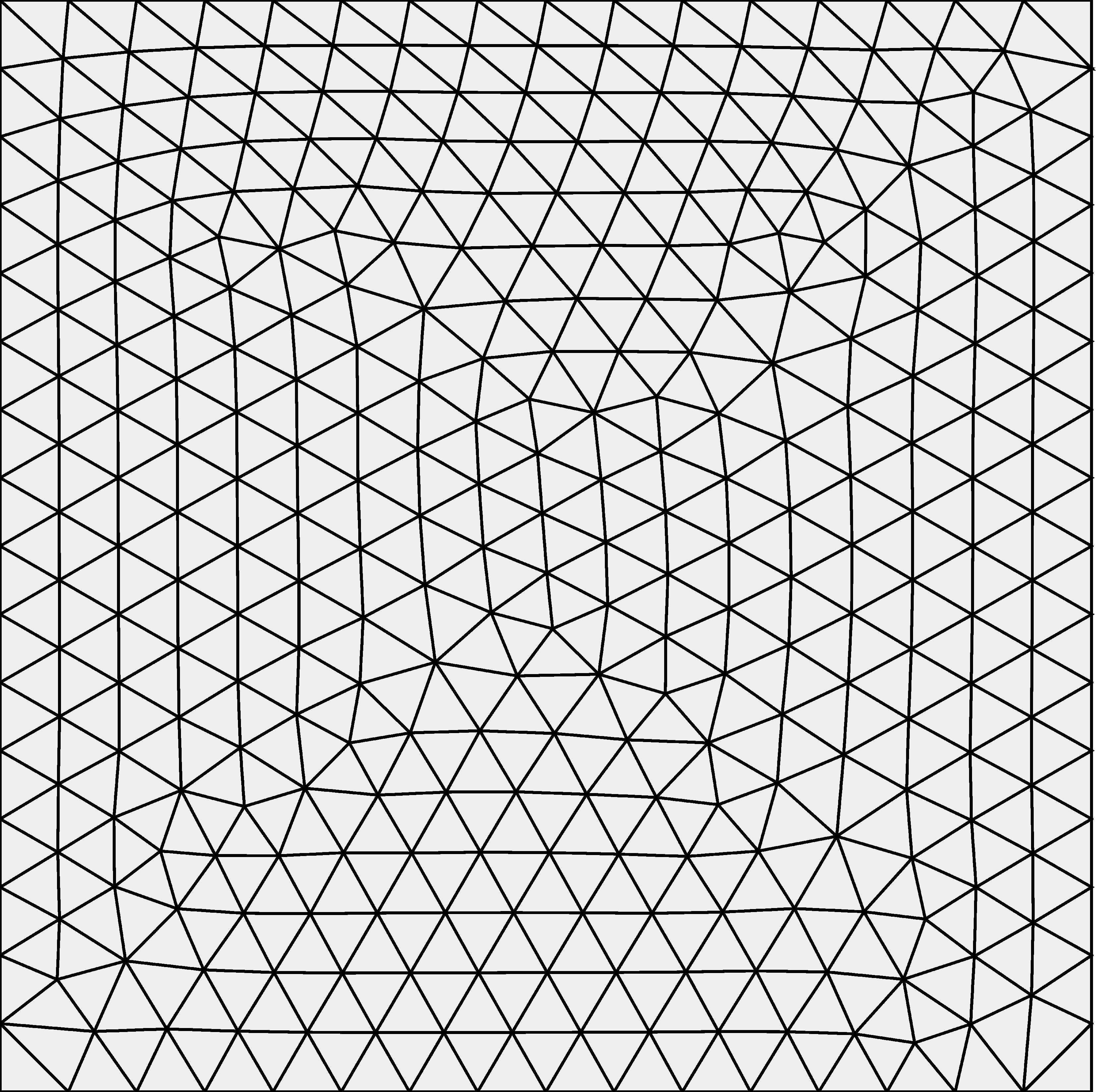}}
 \quad
 \subfigure[] {\label{sfig: grid diffusion gmsh} 
\includegraphics[scale=0.7]{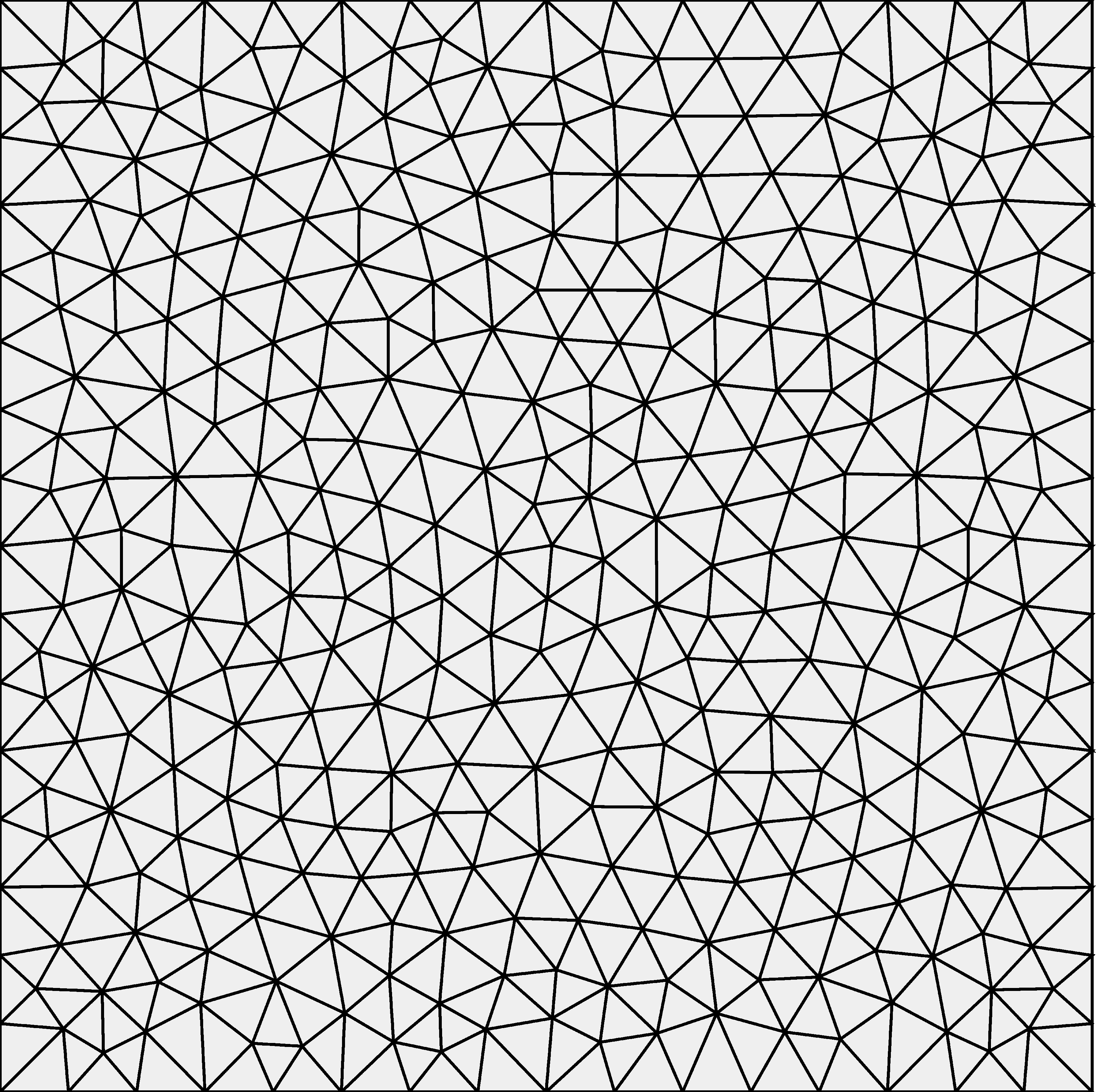}}
 \caption{Samples from each grid category employed to solve the 2D Poisson problems using OpenFOAM: 
\subref{sfig: grid diffusion quads} a grid of distorted quadrilaterals; \subref{sfig: grid 
diffusion 
netgen} a grid constructed with the Netgen algorithm; \subref{sfig: grid diffusion gmsh} a grid 
constructed with the Gmsh algorithm. All domains have unit length along each dimension.}
 \label{fig: grids diffusion}
\end{figure}

In order to investigate further the severity of the problem associated with the use of the common 
DT gradient we performed also a set of experiments using the popular public domain CFD solver 
OpenFOAM (\url{https://openfoam.org}), version 4.0, (see e.g. \cite{Robertson_2015}) which we have 
been using recently in our laboratory \cite{Lampropoulos_2016, Karapetsas_2017}. We solved again a 
Poisson problem \eqref{eq: Poisson} -- \eqref{eq: Poisson BC} but instead of \eqref{eq: Poisson c 
tanh} and \eqref{eq: Poisson b tanh} we have $c(x,y) = \sin(\pi x) \sin(\pi y)$ and $b(x,y) = 2\pi^2 
c(x,y)$, respectively. The problem was solved on the same set of grids composed of distorted 
quadrilaterals (Fig.\ \ref{sfig: grid diffusion quads}), but, as these are artificially distorted 
grids, we also repeated the calculations using a couple of popular grid generation procedures which 
are more typical of real-life engineering applications, namely the Netgen \cite{Schoberl_1997} 
(\url{https://ngsolve.org}) and Gmsh Delaunay \cite{Geuzaine_2009} (\url{http://gmsh.info}) 
algorithms implemented as plugins in the SALOME preprocessor (\url{www.salome-platform.org}). Both 
generate grids of triangular cells, coarse samples of which are shown in Fig.\ \ref{fig: grids 
diffusion}. The Netgen algorithm can be seen to be effective in producing smooth grid, reminiscent 
of that shown in Fig.\ \ref{fig: triangular structured grid}, in four distinct areas within the 
domain that meet at an ``X''-shaped interface where grid distortion is high. The Gmsh grids are less 
regular. The ``LaplacianFoam'' component of OpenFOAM was used to solve the equations with default 
options, which include the DT gradient \eqref{eq: gauss gradient 0} as the chosen gradient scheme 
(option ``gradSchemes'' is set to ``Gauss linear''). We repeated the calculations with the gradient 
option set to LS (gradSchemes: leastSquares).

\begin{figure}[!tb]
 \centering
 \subfigure[mean error, 2D problem]
   {\label{sfig: e_mean OpenFOAM 2D} \includegraphics[scale=1.0]{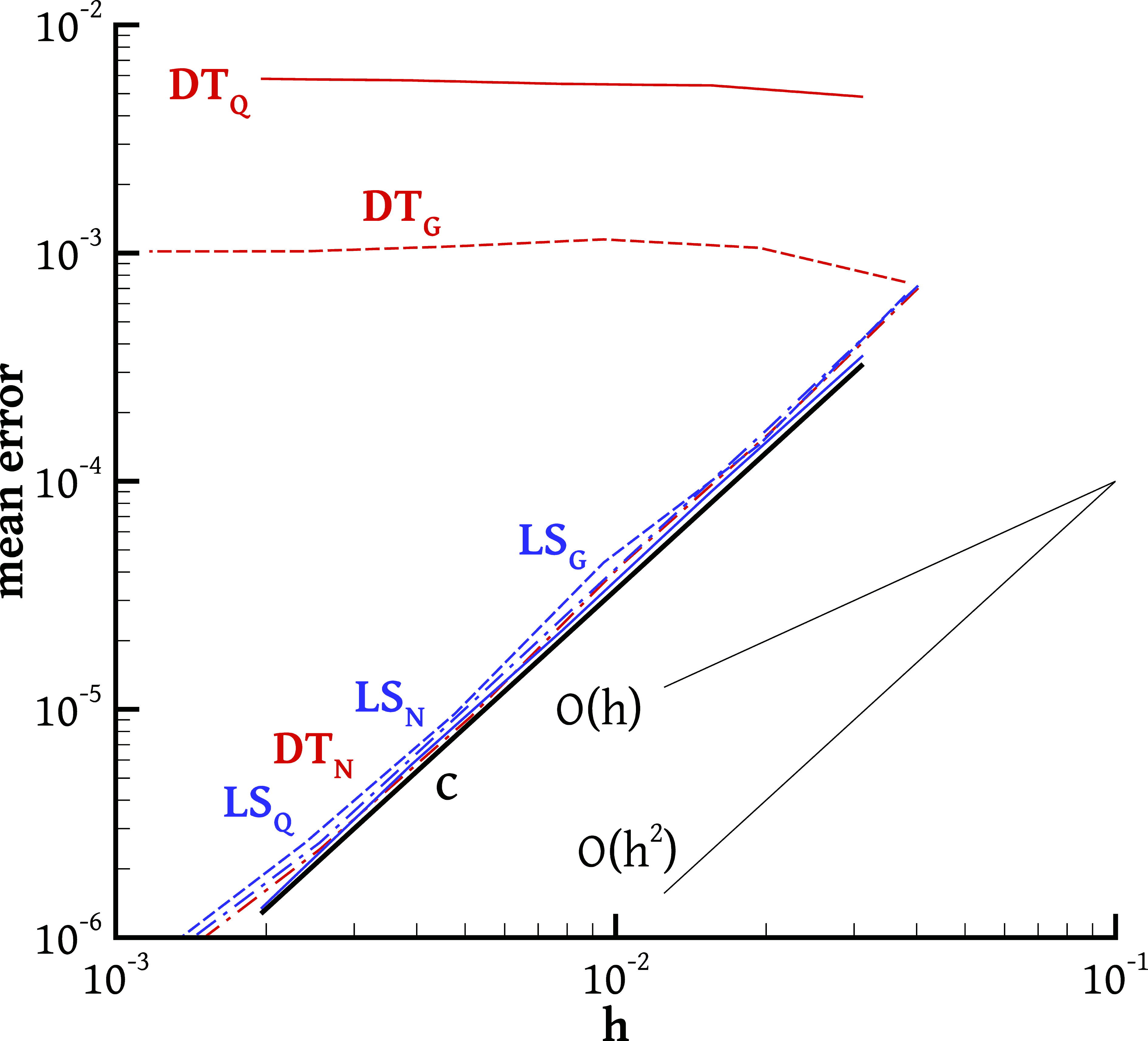}}
 \subfigure[mean error, 3D problem]
   {\label{sfig: e_mean OpenFOAM 3D}  \includegraphics[scale=1.0]{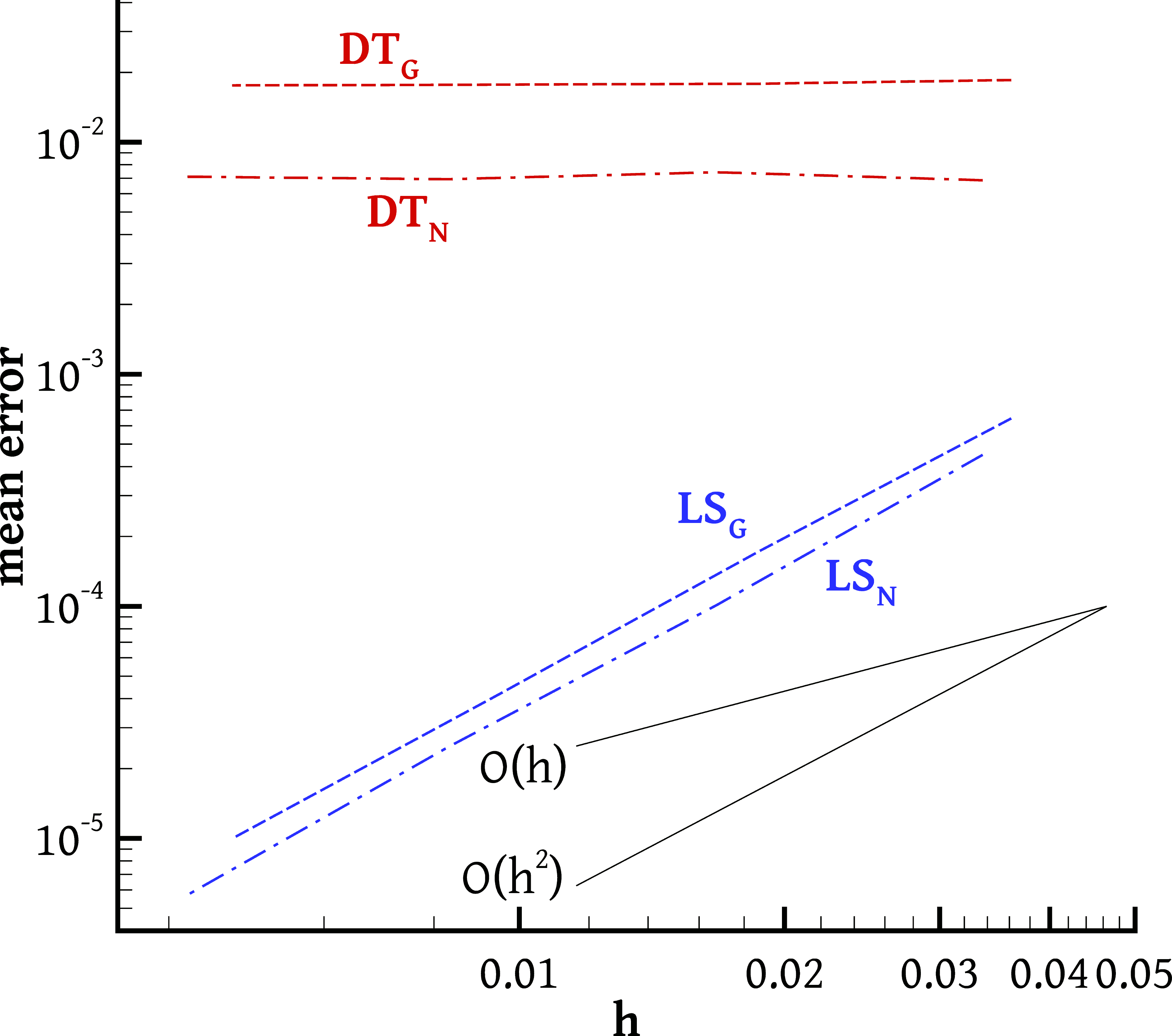}}
 \caption{The mean errors of the OpenFOAM solutions of the 2D \subref{sfig: e_mean OpenFOAM 2D} and 
3D \subref{sfig: e_mean OpenFOAM 3D} Poisson equations, for various configurations. The abscissa 
$h$ is the mean cell size, $h = (\Omega/M)^{1/d}$ where $\Omega = 1$ is the domain volume, $M$ is 
the total number of grid cells, and $d = 2$ (for 2D) or $3$ (for 3D). The black solid line (C) 
corresponds to the 2D problem solved on undistorted Cartesian grids. Red lines (DT) and blue lines 
(LS) correspond to results obtained with the DT and LS gradients, respectively. Subscript Q 
($\mathrm{DT}_{\mathrm{Q}}$, $\mathrm{LS}_{\mathrm{Q}}$ -- solid lines) corresponds to grids of 
distorted quadrilaterals (Fig.\ \ref{sfig: grid diffusion quads}). Subscript N 
($\mathrm{DT}_{\mathrm{N}}$, $\mathrm{LS}_{\mathrm{N}}$ -- dash-dot lines) corresponds to Netgen 
grids (Fig.\ \ref{sfig: grid diffusion netgen}). Subscript G  ($\mathrm{DT}_{\mathrm{G}}$, 
$\mathrm{LS}_{\mathrm{G}}$ -- dashed lines) corresponds to Gmsh Delaunay grids (Fig.\ \ref{sfig: 
grid diffusion gmsh}).}
 \label{fig: errors OpenFOAM}
\end{figure}

The discretisation errors are plotted in Fig.\ \ref{sfig: e_mean OpenFOAM 2D}, against the mean 
cell size\footnote{Unlike in previous plots we use the mean cell size instead of the grid level in 
Fig.\ \ref{fig: errors OpenFOAM} because in the case of the Netgen and Gmsh grids each level does 
not have precisely 4 (in 2D) or 8 (in 3D) times as many cells as the previous level.} $h$. We note 
that on the distorted quadrilateral grids and on the Gmsh grids the DT gradient (lines 
$\mathrm{DT}_{\mathrm{Q}}$ and $\mathrm{DT}_{\mathrm{G}}$) engenders zeroth-order accuracy, like 
the DT0 case of Fig.\ \ref{fig: errors diffusion}. On the other hand, the Netgen grids are quite 
smooth (Fig.\ \ref{sfig: grid diffusion netgen}), resembling grids that come from triangulation of 
structured grids in most of the domain, and this has the consequence that even with the DT gradient 
(lines $\mathrm{DT}_{\mathrm{N}}$) the mean error decreases at a second-order rate. With the LS 
gradient second-order accuracy is exhibited on all grids, and furthermore the accuracy is nearly at 
the same level as on the Cartesian grids.

Three-dimensional simulations were also performed. The governing equations are again \eqref{eq: 
Poisson} -- \eqref{eq: Poisson BC} where we set $c(x,y,z) = \sin(\pi x) \sin(\pi y) \sin(\pi z)$ 
and $b(x,y,z) = 3\pi^2 c(x,y,z)$. The domain is the unit cube $z,y,z \in [0,1]$ and the grids were 
generated with the Netgen and Gmsh algorithms, where in the latter we chose the Delaunay method to 
construct the grid at the boundaries and the frontal Delaunay algorithm to construct the grid in 
the interior. The discretisation errors are plotted in Fig.\ \ref{sfig: e_mean OpenFOAM 3D}. This 
time it is observed that the DT gradient engenders zeroth-order accuracy even on the Netgen grids. 
Presumably, in the 3D case although the Netgen algorithm produces relatively smooth surface grids 
on the bounding surfaces of the cube, in the bulk of the domain it packs the tetrahedra in a way 
that the skewness is large throughout. On the other hand, the LS gradient again always results in 
second-order accuracy.

\section{Final remarks and conclusions}
\label{sec: conclusions}

The previous sections showed that on arbitrary grids the DT and LS gradients are zeroth- and 
first-order accurate, respectively, but higher orders of accuracy, up to second, are obtained on 
particular kinds of grids -- a summary of the present findings is given in Table \ref{table: 
accuracy vs grid}. 

Unfortunately, general-purpose unstructured grid generation algorithms, such as the popular Netgen 
and Gmsh algorithms tested in Sec.\ \ref{sec: PDE solution}, retain high skewness at all levels of 
refinement (except in the 2D Netgen case) resulting in zeroth-order accuracy for both the DT 
gradient and the FVM solver that employs it. In our tests, OpenFOAM using the DT gradient was unable 
to reduce the discretisation errors by grid refinement even though the problems solved were very 
fundamental, namely linear Poisson problems with sinusoidal solutions on the simplest possible 
domains, the unit square and cube. Zeroth-order accuracy is very undesirable because it deprives the 
modeller of his/her main tool for estimating the accuracy of the solution, i.e.\ the grid 
convergence study. Therefore the common DT gradient should be avoided unless the grid refinement 
algorithm is known to reduce the skewness towards zero. This is especially important nowadays that
automatic unstructured grid generation algorithms are preferred over structured grids. On the other 
hand, if the LS gradient is employed instead then the FVM for the problems of Sec.\ \ref{sec: PDE 
solution} proved to be second-order accurate on all grids, irrespective of whether their skewness 
diminishes or not with refinement. The same holds for consistent DT gradient schemes, such as the 
proposed iterative scheme \eqref{eq: gauss outer iterative scheme}, the scheme employing the 
auxiliary cell of Fig.\ \ref{fig: dual cell}, or any other consistent scheme mentioned in Sec.\ 
\ref{sec: gauss}.

\begin{table}[t]
\small
\caption{Order of the mean and maximum truncation errors of the gradient schemes on various grid 
types.}
\label{table: accuracy vs grid}
\vspace{-0.6cm}
{%
\newcommand{\mc}[3]{\multicolumn{#1}{#2}{#3}}
\begin{center}
\begin{tabular}{l||cc|cc|cc}
\hline
       & \mc{2}{c|}{Common DT}  & \mc{2}{c|}{LS $q \neq 3/2$,} &  \mc{2}{c}{LS $q = 3/2$}
\\
       & \mc{2}{c|}{}           & \mc{2}{c|}{$\GRAD{d\infty}$, $\GRAD{dx}$}   & \mc{2}{c}{}
\\
                   & \textit{mean} & \textit{max.} & \textit{mean} & \textit{max.} & \textit{mean} 
& \textit{max.}\\
\hline \hline
Cartesian (\S\ref{ssec: results cartesian}) &
  2             & 1             & 2             & 1             & 2             & 2            \\
Smooth Structured (\S\ref{ssec: results curvilinear}) &
  2             & 1             & 2             & 1             & 2             & 2            \\
Locally Distorted (\S\ref{ssec: results refined}) &
  1             & 0             & 2             & 1             & 2             & 1            \\
Arbitrary (\S\ref{ssec: results random}) &
  0             & 0             & 1             & 1             & 1             & 1            \\
\hline
\end{tabular}
\end{center}
}%
\end{table}

On grids where both methods exhibit the same order of accuracy the optimal method is 
problem-dependent. Usually schemes of the same order of accuracy will produce comparable errors on 
any given grid but in some cases the errors can differ significantly. For example, in aerodynamics 
problems concerning the boundary layer flow over a curved solid boundary, when structured grids 
consisting of cells of very high aspect ratio are employed the DT method is known to significantly 
outperform the LS method \cite{Mavriplis_2003, Sozer_2014} despite both being second-order accurate 
(or first-order, on triangulated grids). This is because the DT method benefits from the alignment 
of the cell faces with the flow (see Appendix \ref{sec: appendix aerodynamics}). 

The present study has been restricted to the FVM solution of a Poisson problem. The effect of the 
choice of gradient discretisation scheme for the FVM solution of flow problems, including 
non-Newtonian or turbulent ones where the role of the gradient is more significant, will be 
investigated in a future study.


\section*{Acknowledgements}

This research has been funded by the LIMMAT Foundation under the Project ``MuSiComPS''.

\begin{appendices}
\renewcommand\theequation{\thesection.\arabic{equation}}

\section{Accuracy of the least squares gradient when the neighbouring points are arranged at equal 
angles}
\label{sec: appendix ls angles}
\setcounter{equation}{0}

Consider the case that all angles between two successive neighbours are equal; this means that 
they are equal to $2\pi/F$ ($F$ being the total number of neighbours). Therefore,
\begin{equation} \label{eq: appendix theta_f}
 \theta_f \;=\; \theta_1 \;+\; (f-1) \, \frac{2\pi}{F}
\end{equation}
As per Eq.\ \eqref{eq: beta_tau component 1}, we will express $\beta_{\varsigma}$ in the error 
expression \eqref{eq: WLS e} substituting for $\varsigma_f$ from Eq.\ \eqref{eq: varsigma_f}. We 
also substitute $\Delta x_f = \Delta r_f \cos \theta_f$ and $\Delta y_f = \Delta r_f \sin \theta_f$. 
In the special case that all neighbours are equidistant from $\vf{P}$, i.e.\ all $\Delta r_f$'s are 
equal, the weights $w_f$ are also equal and the products $(\Delta r_f)^3 w_f^2$ can be factored out 
of the sums:

\begin{equation} \label{eq: appendix beta_tau}
 \beta_{\tau}
 \;=\;
 (\Delta r_f)^3 w_f^2
 \begin{bmatrix}
  \frac{1}{2}\phi_{\!.xx} \sum\limits_{f=1}^F (\cos\theta_f )^3
  \;+\; \phi_{\!.xy} \sum\limits_{f=1}^F (\cos\theta_f )^2 \sin\theta_f
  \;+\; \frac{1}{2}\phi_{\!.yy} \sum\limits_{f=1}^F \cos\theta_f (\sin\theta_f)^2
\\[0.4cm]
  \frac{1}{2}\phi_{\!.xx} \sum\limits_{f=1}^F (\cos\theta_f )^2 \sin\theta_f
  \;+\; \phi_{\!.xy} \sum\limits_{f=1}^F \cos\theta_f (\sin\theta_f)^2
  \;+\; \frac{1}{2}\phi_{\!.yy} \sum\limits_{f=1}^F (\sin\theta_f)^3
 \end{bmatrix}
\end{equation}
where higher-order terms are omitted. If the neighbours are not equidistant but the weight scheme 
$w_f = (\Delta r_f)^{-3/2}$ is used, then all products $(\Delta r_f)^3 w_f^2 = 1$ are again equal 
and \eqref{eq: appendix beta_tau} still holds.

To proceed further, it will be convenient to use complex arithmetic: $\cos\theta_f = (e^{i\theta_f} 
+ e^{-i\theta_f})/2$ and $\sin\theta_f = (e^{i\theta_f} - e^{-i\theta_f})/2i$, where $i \equiv 
\sqrt{-1}$. Then it is straightforward to show that the sums that appear in the right hand 
side of Eq.\ \eqref{eq: appendix beta_tau} are equal to 

\begin{align}
 \sum\limits_{f=1}^F (\cos\theta_f)^3 \;&=\;
 \frac{1}{8} \sum\limits_{f=1}^F \left( e^{i3\theta_f} \;+\; 3e^{i\theta_f} \;+\; 3e^{-i\theta_f} 
\;+\; e^{-i3\theta_f} \right)
 \label{eq: appendix Sc3}
\\
 \sum\limits_{f=1}^F (\cos\theta_f)^2 \sin\theta_f \;&=\;
 \frac{1}{8i} \sum\limits_{f=1}^F \left( e^{i3\theta_f} \;+\; e^{i\theta_f} \;-\; e^{-i\theta_f} 
\;-\; e^{-i3\theta_f} \right)
 \label{eq: appendix Sc2s}
\\
 \sum\limits_{f=1}^F \cos\theta_f (\sin\theta_f)^2 \;&=\;
 -\frac{1}{8} \sum\limits_{f=1}^F \left( e^{i3\theta_f} \;-\; e^{i\theta_f} \;-\; e^{-i\theta_f} 
\;+\; e^{-i3\theta_f} \right)
 \label{eq: appendix Scs2}
\\
 \sum\limits_{f=1}^F (\sin\theta_f)^3 \;&=\;
 -\frac{1}{8i} \sum\limits_{f=1}^F \left( e^{i3\theta_f} \;-\; 3e^{i\theta_f} \;+\; 3e^{-i\theta_f} 
\;-\; e^{-i3\theta_f} \right)
 \label{eq: appendix Ss3}
\end{align}
The above expressions involve the sums $\sum e^{i\theta_f}$, $\sum e^{-i\theta_f}$, $\sum 
e^{i3\theta_f}$ and $\sum e^{-i3\theta_f}$; if these sums are zero then the sums on the left hand 
sides of \eqref{eq: appendix Sc3} -- \eqref{eq: appendix Ss3} are also zero and it follows, in the 
same way as for Eq.\ \eqref{eq: beta_tau component 1}, that the LS gradient is second-order 
accurate. So, we consider the first of these sums, substituting for $\theta_f$ from Eq.\ \eqref{eq: 
appendix theta_f}.

\begin{equation} \label{eq: appendix Seiu}
 \sum\limits_{f=1}^F e^{i\theta_f} \;=\; 
 \sum\limits_{f=1}^F e^{i\left(\theta_1 + (f-1) \frac{2\pi}{F} \right)} \;=\;
 e^{i\theta_1} \sum\limits_{f=1}^F e^{i(f-1) \frac{2\pi}{F}} 
\end{equation}
We denote $\zeta = e^{i2\pi/F}$ and note that it is the first $F$-th root of unity, i.e.\ $\zeta^F 
= e^{i2\pi F/F} = e^{i2\pi} = 1$. Then the sum that appears at the end of Eq.\ \eqref{eq: appendix 
Seiu} is

\begin{equation} \label{eq: appendix sum of roots of unity}
 \sum\limits_{n=0}^{F-1} \zeta^n \;=\; 1 \;+\; \zeta \;+\; \zeta^2 \;+\; \cdots \;+\; \zeta^{F-1} 
\;=\; 0
\end{equation}
Equation \eqref{eq: appendix sum of roots of unity} holds for any $F$-th root of unity except $\zeta 
= 1$, as follows from the identity $1 - \zeta^F = (1 - \zeta)(1 + \zeta + \zeta^2 + \cdots + 
\zeta^{F-1})$: the left hand side is zero because $\zeta^F = 1$, so the right hand side must be zero 
as well, and since $\zeta \neq 1$, it follows that $(1 + \zeta + \zeta^2 + \cdots + \zeta^{F-1}) = 
0$, i.e.\ Eq.\ \eqref{eq: appendix sum of roots of unity} \cite{Strang_2006}. Thus $\sum 
e^{i\theta_f} = 0$. In exactly the same manner it can be shown that $\sum e^{-i\theta_f} = 0$ as 
well.

We proceed in the same way for the sum $\sum e^{i3\theta_f}$:

\begin{equation} \label{eq: appendix Se3iu}
 \sum\limits_{f=1}^F e^{i 3 \theta_f} \;=\; 
 \sum\limits_{f=1}^F e^{i 3 \left(\theta_1 + (f-1) \frac{2\pi}{F} \right)} \;=\;
 e^{i 3 \theta_1} \sum\limits_{f=1}^F e^{i 3 (f-1) \frac{2\pi}{F}} 
\end{equation}
This time we denote $\zeta = e^{i 3 (2\pi) / F}$. It is also an $F$-th root of unity, since $\zeta^F 
= e^{i 3 (2\pi) F/F} = e^{i6\pi} = 1$. So, it will also satisfy Eq.\ \eqref{eq: appendix sum of 
roots of unity}, \textit{unless} $\zeta = 1$. Now, $\zeta = e^{i 2\pi (3/F)}$ will be $\neq 1$ 
precisely if $3/F$ is not an integer, i.e.\ if $F$ is neither $1$ nor $3$. Thus, for $F \neq 3$ 
Eq.\ \eqref{eq: appendix sum of roots of unity} still holds, and Eq.\ \eqref{eq: appendix Se3iu} 
says that $\sum e^{i3\theta_f}$ is zero. But if $F = 3$ then $\sum e^{i3\theta_f}$ is not zero. The 
same holds for the last sum, $\sum e^{-i3\theta_f}$, as can be shown in exactly the same way.

So, to summarise, if $F > 3$ then all the terms in \eqref{eq: appendix Sc3} -- \eqref{eq: appendix 
Ss3} become zero and the leading term of the truncation error vanishes, granting second-order 
accuracy to the LS gradient, provided that the $q = 3/2$ weights scheme is used (or that the 
neighbouring points are equidistant from $\vf{P}$). But if $F = 3$ then not all terms vanish and the 
leading part of the truncation error does not reduce to zero -- the gradient remains first-order 
accurate. This is unfortunate, as the $F = 3$ case is very common, corresponding to triangular 
cells.

\section{Using a LS weights matrix with off-diagonal elements}
\label{sec: appendix ls O(2)}
\setcounter{equation}{0}

Consider the one-dimensional case. For a general matrix $W$ the weighted LS method returns the value 
of $\phi_{\!.x}^{\mathrm{ls}}(\vf{P})$ which minimises the quantity

\begin{equation} \label{eq: NDW W(b-Az)}
 \| W(b-Az) \|^2 \;=\; \sum_i
                       \left[ \sum_{j=1}^F w_{ij} \left( \Delta \phi_j 
                                           - \Delta x_j \, \phi_{\!.x}^{\mathrm{ls}}(\vf{P}) 
\right) \right]^2
\end{equation}
It can be shown that if the method returns the exact derivative of any quadratic function $\phi$ 
then it is second-order accurate. The Taylor expansion of a quadratic function $\phi$ is just

\begin{equation} \label{eq: NDW Taylor 1d}
 \Delta \phi_f \;=\; \phi_{\!.x}(\vf{P}) \, \Delta x_f \;+\; \frac{1}{2} \phi_{\!.xx}(\vf{P}) \, 
(\Delta x_f)^2
\end{equation}
We can use this equation to substitute for $\Delta \phi_j$ into Eq.\ \eqref{eq: NDW W(b-Az)} to get

\begin{equation*}
 \| W(b-Az) \|^2 \;=\; \sum_i \left[
                       \left( \phi_{\!.x}(\vf{P}) - \phi_{\!.x}^{\mathrm{ls}}(\vf{P}) \right) 
\sum_{j=1}^F w_{ij} \Delta x_j
                       \;+\; \frac{1}{2} \phi_{\!.xx}(\vf{P}) \sum_{j=1}^F w_{ij} (\Delta x_j)^2
                       \right]^2
\end{equation*}
Apparently, if we can choose the weights so that $\sum_j w_{ij} (\Delta x_j)^2 = 0$ for each row 
$i$ of $W$, then the terms involving $\phi_{\!.xx}(\vf{P})$ will vanish and the remaining quantity 
will be minimised (actually, made zero) by the choice $\phi_{\!.x}^{\mathrm{ls}}(\vf{P}) = 
\phi_{\!.x}(\vf{P})$. Thus the procedure will produce the exact result. The equation $\sum_j w_{ij} 
(\Delta x_j)^2 = 0$ can be written in matrix notation as $X \, w_i = 0$, where $X$ is the $1 \times 
F$ matrix with elements $X_{1j} = (\Delta x_j)^2$ and $w_i$ is the $F \times 1$ column vector with 
elements $w_{ij}$. We can therefore set each row $w_i$ of $W$ equal to a basis vector of the null 
space of $X$. A simple choice is one where in row $f$ the first element is $(\Delta x_1)^{-2}$, the 
$f\!+\!1$ element is $-(\Delta x_{f+1})^{-2}$, and the rest of the elements are zero. $W$ is no 
longer square but has size $(F\!-\!1) \times F$, so the weighted system $WAz = Wb$ has one less 
equation than the original system $Az = b$. Equation $i$ of the system $W A z = W b$ has the form

\begin{equation} \label{eq: NDW eqn i}
 \phi_{\!.x}^{\mathrm{ls}}(\vf{P}) \, \left( \frac{1}{\Delta x_{i+1}} - \frac{1}{\Delta x_1} \right)
 \;=\;
 \frac{\Delta \phi_{i+1}}{(\Delta x_{i+1})^2} \;-\; \frac{\Delta \phi_1}{(\Delta x_1)^2}
\end{equation}
for $i = 1, \ldots, F-1$. A little manipulation shows that Eq.\ \eqref{eq: NDW eqn i} is just Eq.\ 
\eqref{eq: NDW Taylor 1d} for $f=i+1$, where $\phi_{\!.xx}(\vf{P})$ has been eliminated by using 
again Eq.\ \eqref{eq: NDW Taylor 1d}, but for $f=1$, to express it as a function of $\Delta \phi_1$ 
and $\Delta x_1$. Therefore, this method is nearly equivalent to the unweighted LS solution of the 
system of equations \eqref{eq: NDW Taylor 1d}. Extension to the 2D or 3D cases follows the same 
lines.

\section{The gradient schemes in aerodynamics problems}
\label{sec: appendix aerodynamics}
\setcounter{equation}{0}

It is known in aerodynamics that the DT gradient is more accurate than the LS gradient for the 
computation of boundary layer flow over a curved solid boundary when grid cells of very high aspect 
ratio are employed (typically in excess of 1000) \cite{Mavriplis_2003, Sozer_2014}. This is due to a 
fundamental difference between the methods: the LS gradient uses the directions in which the 
neighbours lie, represented by the unit vectors $\vf{d}_f = (\vf{N}_f - \vf{P}) / \| \vf{N}_f - 
\vf{P} \|$, while the DT gradient uses the directions normal to the faces, represented by the unit 
vectors $\vf{n}_f$. This allows the possibility of assisting the DT gradient by aligning the cell 
faces with the direction of the exact gradient, when this direction is known a priori.

\begin{figure}[thb]
 \centering
 \includegraphics[scale=0.80]{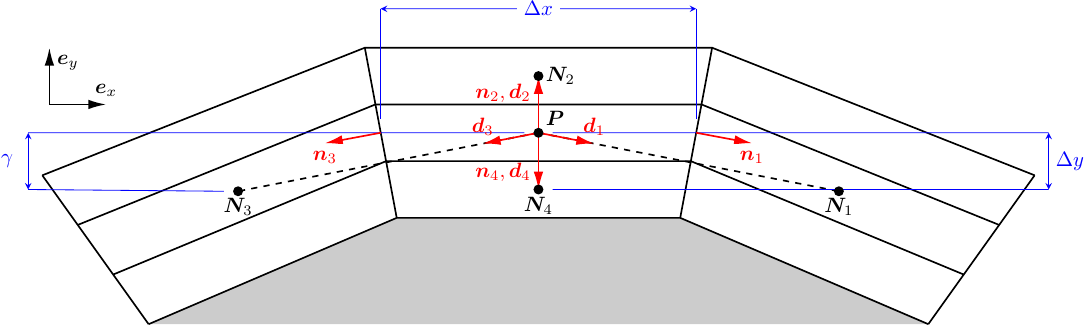}
 \caption{Structured grid of high aspect ratio cells over a curved boundary (the aspect ratio 
$\Delta x / \Delta y$ is greatly downplayed for clarity).}
 \label{fig: aerodynamics grid}
\end{figure}

So, consider a structured grid of high aspect ratio cells over a curved boundary as in Fig.\ 
\ref{fig: aerodynamics grid} and assume a boundary layer flow where the contours of the dependent 
variable $\phi$ follow the shape of the boundary so that $\phi(\vf{P}) = \phi(\vf{N}_1) = 
\phi(\vf{N}_3) \Rightarrow \Delta \phi_1 = \Delta \phi_3 = 0$ and the exact $\nabla \phi$ is 
directed normal to the boundary with approximate magnitude $(\phi(\vf{N}_2) - \phi(\vf{N}_4)) / 
2\Delta y$. The grid being structured all of the examined methods are second order accurate, but the 
unweighted LS gradient is known to perform very poorly: From expression \eqref{eq: minimised q=0} it 
follows that because $\Delta r_1$ and $\Delta r_3$ are much larger than $\Delta r_2$ and $\Delta 
r_4$ this scheme relies on information mostly in the directions $\vf{d}_1$ and $\vf{d}_3$, and since 
the directional derivatives $\Delta \phi_1 / \Delta r_1$ and $\Delta \phi_3 / \Delta r_3$ in these 
directions appear to be zero, the predicted gradient is also near zero although the true gradient 
may be very large. This gross underestimation is known to occur when the ratio of $y$-displacements 
$\gamma / \Delta y$ of neighbours $N_1$ and $N_3$ to neighbours $N_2$ and $N_4$ (Fig.\ \ref{fig: 
aerodynamics grid}) is larger than one (i.e.\ when the terms $(\partial \phi / \partial 
y)^{\mathrm{ls}} \gamma$ are larger than the terms $(\partial \phi / \partial y)^{\mathrm{ls}} 
\Delta y$ in the minimisation expression), and in practical applications it is typically around 50 
\cite{Mavriplis_2003}. As this ratio is proportional to $(\Delta x)^2 / (R \, \Delta y)$ 
\cite{Mavriplis_2003}, reducing the size $\Delta x$ of the cells in the streamwise direction rapidly 
reduces the error. Using a weighted LS method dramatically improves the accuracy, but this remedy 
does not work if the grid is triangulated as in Fig.\ \ref{fig: triangular structured grid}.

The inaccuracy of the LS gradients in this case is due to the invalidation of their fundamental 
assumption that the variation of $\phi$ in the neighbourhood of the cell is linear, as the cell size 
is large enough for the contours of $\phi$ to curve significantly across it.  On the other hand, the 
DT gradient benefits from the alignment of the normals $\vf{n}_2$ and $\vf{n}_4$ of the long faces 
(which dominate the calculation due to weighting by face area) with the direction of the exact 
gradient. Thus, for this particular $\phi$ distribution the given grid alignment assists this 
gradient scheme to achieve good accuracy. A similar (somewhat less favourable) situation holds in 
the triangulated grid case.

Another issue is related to numerical stability. Certain discretisation schemes use the values of 
the dependent variable at face centres, evaluated as $\phi(\vf{c}_f) \approx \phi(\vf{P}) + 
\nabla^{\mathrm{a}}\phi(\vf{P}) \cdot (\vf{c} - \vf{P})$. Consider the same boundary layer flow as 
before, but for the moment assume a straight boundary so that the grid looks like that of Fig.\ 
\ref{fig: parallelogram grid}, on which all the examined gradient schemes reduce to 
$\nabla^{\mathrm{s}}$, Eq.\ \eqref{eq: grad parallelogram}. Since $\phi(\vf{N}_1) = \phi(\vf{N}_3) = 
\phi(\vf{P})$, Eq.\ \eqref{eq: grad parallelogram} returns a gradient in the direction 
$(-\delta_y^{\xi}, \delta_x^{\xi})$ i.e.\ perpendicular to the vector $\vf{\delta}^{\xi}$. Thus 
$\phi(\vf{c}_1) \approx \phi(\vf{P}) + \nabla^{\mathrm{s}}\phi(\vf{P}) \cdot (\vf{c}_1 - \vf{P}) = 
\phi(\vf{P}) = \phi(\vf{N}_1)$, which is reasonable. Now consider a curved boundary case where cell 
$P$ has exactly the same shape as before but its neighbours $N_1$ and $N_3$ are tilted to follow the 
curvature as in Fig.\ \ref{fig: aerodynamics grid}. The DT gradient depends only on the face normals 
$\vf{n}_f$ which have not changed for cell $P$ and thus it gives the same result for 
$\phi(\vf{c}_1)$ as before. On the other hand, the vectors $\vf{d}_1$ and $\vf{d}_3$ have now 
changed and thus the LS gradient is affected by the curvature and may no longer be perpendicular to 
the direction $\vf{c}_1 - \vf{P}$; this leads to a predicted value for $\phi(\vf{c}_1)$ that is 
either larger or smaller than both $\phi(\vf{P})$ and $\phi(\vf{N}_1)$, introducing an oscillatory 
variation in the direction $\vf{d}_1$ that may be detrimental to numerical stability 
\cite{Shima_2010}. In other cases this issue can occur with both the LS and the DT gradients, and 
so some solvers such as OpenFOAM offer versions of these gradient schemes that use limiters in order 
to avoid this problem.

\end{appendices}

\bibliographystyle{ieeetr}
\bibliography{gradient_calculation}

\end{document}